%% file: QuadForms_V13.tex
\documentclass{article}
\input{format}
\usepackage{amsmath,amssymb,enumerate,amsbsy}
\usepackage{color}

\numberwithin{equation}{section}
	\newcommand{\xtra}[1]{}
	
\begin{document}
\title{Simultaneous Integer Values of \\Pairs of Quadratic Forms}

\author{D. R. Heath-Brown\footnote{Mathematical Institute, Radcliffe Observatory Quarter, Woodstock Road, Oxford
OX2~6GG, \emph{rhb@maths.ox.ac.uk}} \\
L. B. Pierce\footnote{Hausdorff Center for Mathematics, 60 Endenicher
  Allee, 53115 Bonn, \emph{pierce@math.uni-bonn.de}} }
  \date{}
\maketitle

\begin{abstract}
We prove that a pair of integral quadratic forms in 5 or more
variables will simultaneously represent ``almost all'' pairs of integers
that satisfy the necessary local conditions, provided that the forms
satisfy a suitable nonsingularity condition.  In particular such
forms simultaneously attain prime values if the obvious local
conditions hold. The proof uses the circle method, and
in particular pioneers a two-dimensional version of a Kloosterman
refinement. 

\end{abstract}
\noindent
\nonumfootnote{\emph{2010 Mathematics Subject Classification:} 11P55, 11D85}

\tableofcontents

\section{Introduction}\label{sec_intro}
Given two quadratic forms $Q_1, Q_2 \in \Z[x_1,\ldots, x_k]$, we would
like to understand which pairs of integers $n_1,n_2$ are represented
simultaneously by $Q_1$ and $Q_2$.  The situation for a single form is
fairly well understood, but less is known for pairs of
forms. Naturally there may be local obstructions to representability.
However one might expect that if the variety $V$ defined by the simultaneous
equations \begin{equation}\label{V_dfn}
\left\{\begin{array}{ccc}Q_1(\xbf) &=& 0 \\
Q_2(\xbf) &=& 0 \end{array}\right.
\end{equation}
is nonsingular, and if $k$ is large enough, then every pair $n_1,n_2$
satisfying the necessary local conditions should be representable, or
at least that all but finitely many such pairs are representable.

The method of Birch \cite{Bir61} can be adapted to obtain an asymptotic for the number of
representations of a pair of integers $n_1,n_2$ by quadratic forms
$Q_1,Q_2$ for $k \geq 13 +s$ variables, where  $s$ denotes the
dimension of the set \[\{ \xbf \in \Qbar^k: \rank \left( \begin{array}{c} \nabla Q_1(\xbf)
 \\ \nabla Q_2(\xbf)  \end{array} \right) <2 \};\] in particular, if the variety (\ref{V_dfn}) is nonsingular then $s=1$.
In a recent preprint, Munshi \cite{Mun13} employs a new `nested
$\del$-method' version of the circle method to reduce the number of
variables to $k \geq 11$ if (\ref{V_dfn}) is nonsingular.

In order to handle as small a value for $k$ as possible, this paper asks
only that ``almost all'' suitable pairs are representable.  By this we
mean that the number of suitable pairs for which
$|n_1|,|n_2|\le N$  but $n_1,n_2$ are not simultaneously
representable, should be $o(N^2)$ as $N$ tends to infinity.

As a by-product of our investigation we will be able to say something
about pairs of quadratic forms which simultaneously take prime values
infinitely often. Certainly such pairs of forms exist: for example, $Q_1$ and
$Q_2$ defined by \begin{eqnarray*}
Q_1(x_1,x_2,x_3) &=& x_1^2 + 4x_2^2\\
Q_2(x_1,x_2,x_3) &=& 4x_1^2 + x_3^2,
\end{eqnarray*}
simultaneously attain prime values infinitely often. Indeed
one can use the result of Fouvry and Iwaniec \cite {FI} to show that
for ``almost all'' odd integers $x_1$ there is a prime of the form
$x_1^2+4x_2^2$, and equally that for almost all such $x_1$ there is a
prime of the form $4x_1^2+x_3^2$.  The claim then follows easily.
\label{note_FI}

The geometric condition we must impose on the pair of quadratic forms
is as follows: \begin{condition}\label{cond1}
The projective variety defined over $\overline{\Q}$ by
\beq\label{res_var}
V : Q_1(\xbf)=Q_2(\xbf)=0
\eeq
is nonsingular, by which we mean for every $\xbf \in \overline{\Q}^k$, if
$Q_1(\xbf)=Q_2(\xbf)=0$ with $\xbf \neq 0$ then \beq\label{res_rk}
\rk \left( \begin{array}{c}
	 \nabla Q_1(\xbf) \\
	 \nabla Q_2(\xbf)
	 \end{array}
	 \right)
	 =2.
	 \eeq
\end{condition}

We will count solutions using a smooth non-negative weight $w: \R^k \maps
\R$ of compact support. We then define the
normalized weight function \beq\label{wt_dfn}
w_B({\bf x}) = w(B^{-1}{\bf x})
\eeq
and the weighted representation function \beq\label{RBw}
R_B(n_1,n_2)=\sum_{\substack{\xbf \in \Z^k \\ \u{Q}(\xbf)=\u{n}}}w_B(\xbf),
\eeq
where we use an underscore to denote variables in 2-dimensional spaces,
so that $\u{Q}(\xbf)=(Q_1(\xbf),Q_2(\xbf))$, for example.

By an ``integral quadratic form'' we will always mean a form with
integral matrix, so that the off-diagonal terms of the form
have even coefficients.
Our principal result is then the following, in which we set
\[ F(x,y) = \det(xQ_1 + yQ_2).\]
\begin{thm}\label{thm_r1}
Suppose we have two integral quadratic forms $Q_1, Q_2$ satisfying Condition 1, in $k \geq 5$ variables.
Then if $N=B^2$ we have
\beq\label{tt}
\sum_{\substack{\max\{|n_1|,|n_2|\}\le N\\ F(n_2,-n_1)\not=0}} |
R_B(\u{n})-
\Jcal_w(B^{-2}\u{n})\Sf(\u{n})B^{k-4}|^2\ll_{Q_1,Q_2,w} B^{2k-4-1/16},
\eeq
where $\Jcal_w(\u{\mu})$ and $\Sf(\u{n})$ are the singular integral and
singular series, given respectively by (\ref{sing_int_dfn}) and (\ref{sing_ser_dfn}). \end{thm}

Of course this is of little use without some information about
$\Jcal_w(\u{\mu})$ and $\Sf(\u{n})$.  This will be provided in our second result.
\begin{thm}\label{thm_r2}
Let $Q_1$ and $Q_2$ be as in Theorem \ref{thm_r1} and suppose that $k\ge 5$
and $F(n_2,-n_1)\not=0$. Then
$\Sf(\u{n})\ll_{\ep,Q_1,Q_2} \max(|n_1|,|n_2|)^{\ep}$ for any $\ep>0$.
Moreover if the system of equations
$Q_1(\xbf)=n_1$, $Q_2(\xbf)=n_2$ is solvable in every $p$-adic ring
$\Z_p$, then $\Sf(\u{n})$ is real and positive.  Indeed there is then a
constant $p_0$ depending only on $Q_1$ and $Q_2$ such that
\[\Sf(\u{n})\gg_{\ep,Q_1,Q_2} \max(|n_1|,|n_2|)^{-\ep}\prod_{p\le p_0}|F(n_2,-n_1)|_p^{k-2},\]
for any fixed $\ep>0$,
where $|\cdot|_p$ denotes the standard $p$-adic valuation.

Similarly, under the same conditions on $\u{Q}$, for any smooth weight
$w$ of compact support we have
$\Jcal_w(\u{\mu})\ll_{Q_1,Q_2,w} 1$.  Moreover there is a constant $C$,
depending on $Q_1$ and $Q_2$, with the following property. Suppose
that $w(\xbf)>0$ for $|\xbf|\le C$.  Then we have
\[\Jcal_w(\u{\mu})\gg_{Q_1,Q_2,w} |F(\mu_2,-\mu_1)|^{k-2}\]
for any $\u{\mu}$ in the region $1/2\le\max(|\mu_1|,|\mu_2|)\le 1$ for which
the system of equations
$Q_1(\xbf)=\mu_1$, $Q_2(\xbf)=\mu_2$ has a solution $\xbf\in\R^k$.
\end{thm}

As a corollary we have the following result.
\begin{thm}\label{thm_gen}
Suppose we have two integral quadratic forms $Q_1, Q_2$ satisfying Condition 1, in $k \geq 5$ variables.
Let $\Ecal(N)$ denote the number of integer pairs $(n_1,n_2)$ with $|n_1|,|n_2|\le N$ for which the system
$Q_1(\xbf)=n_1$, $Q_2(\xbf)=n_2$ has a real solution and solutions in
each $\Z_p$, but for which there is no solution $\xbf\in \Z^k$.  Then
\[ \Ecal(N) \ll_{Q_1,Q_2,\varpi} N^{2-\varpi},\]
with $\varpi=1/(8k^3)$.
\end{thm}
We have made no effort to get the best possible exponent $\varpi$
here, but note that our value is dependent only on $k$.

In particular, we may derive from Theorem \ref{thm_gen} the following
result on simultaneous prime values: \begin{thm}\label{thm_prime}
Suppose two integral quadratic forms $Q_1(\xbf), Q_2(\xbf)$
satisfy Condition 1 with $k \geq 5$. 
Suppose further that there is an $\xbf_0 \in \R^k$ such that
$Q_1(\xbf_0), Q_2(\xbf_0)>0$, and that for every prime $q$ there is an
$\xbf_q \in \Z^k$ for which $q \ndiv Q_1(\xbf_q)Q_2(\xbf_q)$. Then there are infinitely many pairs of primes simultaneously
representable by $Q_1(\xbf)$ and $Q_2(\xbf)$. \end{thm}

We will prove in Section \ref{sec_geom_cond} that for diagonal
quadratic forms, say $Q_1 = \sum a_ix_i^2$ and $Q_2 = \sum b_ix_i^2$,
Condition 1 is equivalent to the condition that the ratios $a_i/b_i$
are all distinct, for $i = 1,\ldots, k$. As a result, it is simple to find examples of pairs of forms that
satisfy the requirements of Theorem \ref{thm_prime}. For instance, one
may take \begin{eqnarray*}
Q_1(\xbf) &=& x_1^2 + x_2^2 + x_3^2 + x_4^2 + x_5^2 \\
Q_2(\xbf) & = & x_1^2 + 2x_2^2 +3 x_3^2 + 4x_4^2 + 5x_5^2,
\end{eqnarray*}
for which the choices $\xbf_0 = (1,0,\ldots,0)$, $\xbf_q =
(1,0,\ldots,0)$ clearly suffice.

If Condition \ref{cond1} is dropped from the hypotheses of Theorem
\ref{thm_prime}, then the result of the theorem can fail to hold. For
example, the pencil generated by $Q_1(\xbf) = x_1^2,$ $Q_2(\xbf) =
x_1^2 + \cdots + x_k^2$ is singular as soon as $k \geq 3$, and
certainly $Q_1$ will never attain prime values. In fact, it is
reasonable to conjecture that this is representative of the only type
of exception to arise. Motivated by Schinzel's Hypothesis, one may
conjecture that the result of Theorem \ref{thm_prime} should continue
to hold if the assumption of Condition \ref{cond1} is replaced by the
assumption that neither $Q_1(\xbf)$ nor $Q_2(\xbf)$ factors over $\Z$;
this is discussed further in Section \ref{sec_final}. \label{note_Schinzel}

We should remark at this point that one can handle ``target sets'' other
than the primes in much the same way.  Provided that the target set
has a counting function that grows faster than $N^{1-\varpi/2}$, all
that is necessary is that
one should be able to deal satisfactorily with the local conditions
that arise. In particular we can handle the case in which the target
set consists of the integers that are sums of two squares.  This
leads to an analytic proof of the following example of the Hasse Principle
(which is a special case of a far more general result due to
Colliot-Th\'el\`ene, Sansuc and Swinnerton-Dyer \cite[Theorem A, case
(i)(a)]{CTSSD}).
\label{note_CTSSD}

\begin{cor}
Let $q_1(\xbf,\ybf) = Q_1(\xbf)  - y_1^2  - y_2^2,$ $q_2(\xbf,\ybf) =
Q_2(\xbf)  - y_3^2  - y_4^2,$ where $Q_1,Q_2 \in \Z[x_1,\ldots, x_k]$
are integral quadratic forms satisfying Condition 1. Then
the Hasse Principle holds for the intersection $q_1(\xbf,\ybf) =
q_2(\xbf,\ybf)=0$ as soon as $k+4\ge 9$.
\end{cor}

We leave the details of the proof to the reader. Here $k+4$ is the
 total number of variables in the system $q_1(\xbf,\ybf) = q_2(\xbf,\ybf)=0$.  This result may be compared to
 the recent work of Browning and Munshi \cite{BrMu13} on the Hasse
 Principle for the intersection $q_1(\xbf,\ybf) = q_2(\xbf,\ybf)=0$.
 Browning and Munshi assume that the forms $q_1$ and $q_2$ take
the shape $q_1(\xbf,\ybf) = Q_1(\xbf) - y_1^2 - y_2^2$,
 $q_2(\xbf,\ybf) = Q_2(\xbf)$ in $k+2 \geq 9$ variables, where the
 intersection $Q_1(\xbf) = Q_2(\xbf) =0$ is assumed to be
 nonsingular.

We observe that the range
$k\ge 5$ in Theorem \ref{thm_gen} is best possible in the sense that
the theorem statement would be false with $k=4$.  To justify this we first note that for any
integer $L\ge 2$ the forms
\beq\label{Qpair1}
Q_1(\xbf)=x_1^2+x_3^2+x_4^2,\;\;\; Q_2(\xbf)=x_2^2+x_3^2+Lx_4^2
\eeq
satisfy Condition 1. It can be shown (by a method similar to the
 derivation of Theorem \ref{thm_prime} in Section \ref{sec_final})
 that there is a value $p_0$ depending only on $Q_1, Q_2$ such that
 the equations \beq\label{Q1Q2n}
Q_1(\xbf)=n_1,\;\;\; Q_2(\xbf)=n_2
\eeq
will have a nonsingular solution in $\F_p$, and hence a solution in
$\Z_p$, for any $p\ge p_0$ not dividing both $L$ and $n_2$. Moreover, there exists a modulus $M$ and
residue classes $a_1,a_2$ dependent only on $Q_1,Q_2$ and the primes
$p<p_0$ such that the system (\ref{Q1Q2n}) has a solution over $\Z_p$
for $p < p_0$ as long as \[n_1 \con a_1, n_2 \con a_2 \modd{M}. \]
For
the pair (\ref{Qpair1}), computation shows that $p_0=11, a_1=a_2=1,$
and $M=840$ suffice, so that taking $L$ as a large prime, one may
conclude that the system has solutions over every $\Z_p$ provided that
$n_1\equiv n_2\equiv 1\modd{840}$ and $L\nmid n_2$.  Moreover there
will be a real solution (with $x_3=x_4=0$) whenever $n_1$ and $n_2$
are positive.  Thus there are $\gg N^2$ pairs of positive integers
$n_1,n_2\le N$ for which the local conditions are everywhere
satisfied, with an implied constant independent of $L$. However for a global solution one clearly has
$|x_1|,|x_2|,|x_3|\le N^{1/2}$ and $|x_4|\le N^{1/2}L^{-1/2}$ so that
there are $O(N^2L^{-1/2})$ possible 4-tuples $(x_1,\ldots,x_4)$.  Thus
if $L$ is chosen sufficiently large there will be $\gg N^2$ pairs
$n_1,n_2$ for which there are local solutions but no global solution.
\label{note_k4_ex}

Thus far we have interpreted the results of Theorems \ref{thm_r1} and
\ref{thm_r2} as providing ``almost-every'' results. In particular,
these theorems verify that for a pair of integral quadratic forms satisfying
Condition 1 in $k \geq 5$ variables, for ``almost every'' pair
$(n_1,n_2)$ the counting function $R_B(\u{n})$ is  asymptotically
equal to $\Jcal_w(B^{-2}\u{n})\mathfrak{S}(\u{n})B^{k-4}$. Precisely,
if $N=B^2$ and for each $0 < \theta < 1/32$ we let $\Ecal_\theta(N)$ denote the
number of $|n_1|, |n_2| \leq N$ for which either $F(n_2,-n_1) =
0$ or the difference \[R_B(\u{n}) - \Jcal_w(B^{-2}\u{n})\mathfrak{S}(\u{n})B^{k-4}\]
fails to be $O(B^{k-4-\theta})$, then Theorem \ref{thm_r1} shows that
$\Ecal_\theta(N) \ll N^{2 - (\frac{1}{32}- \theta)}$.

Our results also imply universal results for certain forms in $k \geq
10$ variables. Suppose that $Q_1, Q_2$ are integral forms in $k_1+k_2$
variables that split, so that we may write \begin{eqnarray}\label{Qpair}
\begin{array}{lll}Q_1(\xbf,\ybf) &=& q_1(\xbf) + q_2(\ybf) \\
Q_2(\xbf,\ybf) & = & q_3(\xbf) + q_4(\ybf), \end{array}
\end{eqnarray}
where $\xbf \in \Z^{k_1}, \ybf \in \Z^{k_2}$, $q_1,q_3$ are integral quadratic
forms in $k_1$ variables and $q_2,q_4$ are integral quadratic forms in $k_2$
variables. Note that if the pair $Q_1,Q_2$ satisfies Condition 1 then each of the
pairs  $q_1,q_3$ and $q_2,q_4$ also satisfies Condition 1 (but the
reverse implication need not hold); this is visible by applying
Condition 2 (see section \ref{sec_geom_cond}) to the factorization \[ \det(xQ_1 +yQ_2) = \det(xq_1+yq_3)\det(xq_2 + yq_4).\]
Our results (in particular Propositions \ref{prop_major_int} and
\ref{prop_minor_arcs}) imply immediately the following, by a classical
application of the circle method: \begin{thm}\label{thm_s1}
Suppose that $k_1, k_2\geq 5$ and that we have two integral quadratic forms $Q_1,Q_2 \in \Z[X_1,\ldots, X_{k_1+k_2}]$
 satisfying Condition 1. Assume further that $Q_1$ and $Q_2$ are
split as in
(\ref{Qpair}).   Then for any $(n_1,n_2)$ such that $F(n_2,-n_1) \neq 0$, \[ R_B(n_1,n_2)  = \Sf(\u{n}) \Jcal_w(B^{-2}\u{n}) B^{k_1+k_2-4} + O(B^{k_1+k_2 -4-1/32}),\]
where as usual $R_B(n_1,n_2)$ is defined by (\ref{RBw}) and the
singular series $\Sf(\u{n})$ and singular integral $\Jcal_w(\u{n})$ are given respectively by (\ref{sing_ser_dfn}) and (\ref{sing_int_dfn}), and satisfy the properties of
Theorem \ref{thm_r2}. \end{thm}

Theorem \ref{thm_s1} does not apply to the case $(n_1,n_2)=(0,0)$,
 but we are able to modify our argument to establish the following
 results.
\begin{thm}\label{thm_s3}
Let $k_1, k_2\geq 5$ and suppose that we have two integral quadratic forms
$Q_1,Q_2 \in \Z[X_1,\ldots, X_{k_1+k_2}]$ satisfying Condition
1. Assume further that $Q_1$ and $Q_2$ are
split as in
(\ref{Qpair}).  Then
\[ R_B(0,0)  = \Sf(\u{0}) \Jcal_w(\u{0}) B^{k_1+k_2-4} + O(B^{k_1+k_2 -4-1/32}),\]
where the  singular series $\Sf(\u{0})$ and singular  integral
$\Jcal_w(\u{0})$ are given respectively by (\ref{sing_ser_dfn}) and
(\ref{sing_int_dfn}), and satisfy the properties of Theorem \ref{thm_s2} below.
\end{thm}
\begin{thm}\label{thm_s2}
Let $Q_1$ and $Q_2$ be as in Theorem \ref{thm_s3}. If the system of equations
$Q_1(\xbf)=Q_2(\xbf)=0$ has a nonzero solution in every $p$-adic ring
$\Z_p$, then $\Sf(\u{0})$ is real and positive. Suppose further that the weight $w(\x)$ is supported on a neighbourhood
of the origin.   Then if the system of equations
$Q_1(\xbf)=Q_2(\xbf)=0$ has a nonzero solution $\xbf\in\R^k$, we have
$\Jcal_w(\u{0})>0$.
\end{thm}

\subsection{Method of proof}\label{sec_method}
Our approach uses the two-dimensional circle method. The novel idea that
allows us to reduce to $k \geq 5$ is a two-dimensional Kloosterman
refinement applied to the contribution of the minor arcs.

The proofs of our results are
long and involved.  In order to aid the reader, we begin by sketching the
method of proof schematically, as follows.  With the notations
introduced above we set
\beq\label{S_dfn1}
S(\al_1,\al_2) = \sum_{\x \in \Z^k} e(\al_1 Q_1(\xbf) + \al_2
Q_2(\xbf))w_B(\xbf),
\eeq
for any $(\al_1, \al_2) \in \R^2$, so that \[ R_B(n_1,n_2) = \int_0^1 \int_0^1 S(\al_1,\al_2) e(-\al_1 n_1 - \al_2 n_2)
d\al_1 d\al_2 . \] With a suitable definition of the major arcs $\Mf$ and minor arcs
$\mf$ in two dimensions, we may break $R_B$ into the two pieces \[ R_B(n_1,n_2) = \iint_{\Mf} + \iint_{\mf}\]
so that
\begin{eqnarray}\label{RMm}
\lefteqn{\sum_{\u{n}\in\Z^2} \left|R_B(n_1,n_2) - \iint_\Mf
  S(\al_1,\al_2)e(-\al_1n_1 - \al_2n_2) d\al_1 d\al_2 \right|^2 }
\hspace{2cm}\nonumber	\\ &=& \sum_{\u{n}\in\Z^2} \left|\iint_{\mf}
 S(\al_1,\al_2)e(-\al_1n_1 - \al_2n_2) d\al_1 d\al_2\right|^2. \end{eqnarray}
Temporarily set $f(\al_1,\al_2) = S(\al_1,\al_2)
\chi_{\mf}(\al_1,\al_2)$. Then the right hand side of (\ref{RMm}) is
\begin{eqnarray}\label{Parseval}
\lefteqn{\sum_{\u{n}\in \Z^2} \left| \iint_{[0,1]^2} f(\al_1,\al_2)
e(-\al_1n_1 - \al_2n_2) d\al_1 d\al_2\right|^2}\hspace{2cm}\nonumber\\
&=& \iint_{[0,1]^2} |f(\al_1,\al_2)|^2 d\al_1 d\al_2 \nonumber\\
&=& \iint_{\mf} |S(\al_1,\al_2)|^2 d\al_1 d\al_2, \end{eqnarray}
where in the first equality we have applied Parseval's identity.

We will show that if $\max\{|n_1|,|n_2|\}\le N=O(B^2)$ and
$F(n_2,-n_1)\not=0$, the contribution of the major arcs may be approximated as
\beq\label{maj_ME}
\iint_{\Mf} S(\al_1,\al_2) e(-\al_1n_1 - \al_2n_2)d\al_1d\al_2=
M(n_1,n_2) + E(n_1,n_2),
\eeq
where  $M(n_1,n_2)$ is the expected main term, including the singular
series and the singular integral, and $E(n_1,n_2)$ is an acceptable
error term. Thus it will follow that \begin{eqnarray}\label{min_E}
\lefteqn{ \sum_{\substack{\max\{|n_1|,|n_2|\}\le N\\ F(n_2,-n_1)\not=0}} |R_B(n_1,n_2) - M(n_1,n_2)|^2}\hspace{2cm}\nonumber\\
& \ll&
\iint_{\mf} |S(\al_1,\al_2)|^2 d\al_1 d\al_2 + E(N),
\end{eqnarray}
where \[ E(N) = \sum_{\u{n}\in\Z^2} |E(n_1,n_2)|^2.\]

When $\max(|n_1|,|n_2|)$ is of order $B^2$ the expected size of $M(n_1,n_2)$ is $B^{k-4}$, since roughly
speaking, this is the probability that as $\xbf$ ranges over $B^k$
possible choices, both of the values $Q_1(\xbf)-n_1$ and
$Q_2(\xbf)-n_2$ (each of size up to $O(B^2)$), is zero. Imagine for the moment that $M(n_1,n_2)\gg B^{k-4}$; this implicitly assumes that the singular
series and singular integral have suitable lower bounds. (In fact,
proving such lower bounds is a significant source of complication.)
Suppose furthermore that the right hand side of (\ref{min_E}) is
bounded by $B^\be$ for some $\be>0$. Under these two significant
assumptions, (\ref{min_E}) would imply that \[ \# \{(n_1,n_2)\in\Z^2:\, \max(|n_1|,|n_2|)\ll B^2: R_B(n_1,n_2)=0 \} \cdot (B^{k-4})^2
	\ll B^\be,\]
	whence
\beq\label{RB1}
\# \{(n_1,n_2)\in\Z^2:\, \max(|n_1|,|n_2|)\le N \,: R_B(n_1,n_2)=0 \}
	\leq B^{\be - (2k-8)}
	\eeq
for $N= O(B^2)$. If $\be$ is such that the right hand side of (\ref{RB1}) is $o(N^2)$, then we may
conclude that almost all pairs of integers $n_1, n_2$ of size $N$
are represented simultaneously by $Q_1, Q_2$.

Thus we seek to bound (\ref{min_E}) by $B^\be$ with $\be < 2k-4$.
The mean square argument above has appeared before in settings in
which the circle method is used to show that almost all integers (in
some suitable sense) are represented by a particular form.  It dates
from the work of Hardy and Littlewood \cite{HL} who showed, for example, that almost all positive integers are the sum of 5
non-negative cubes.

Our main goals now are therefore a representation of the contribution of
the major arcs in the form (\ref{maj_ME}), and a bound for the right
hand side of (\ref{min_E}) of the form
$o(B^{2k-4})$. In particular, this will require an estimate for the
integral over the minor arcs of the shape
\beq\label{minor_intro}
\iint_{\mf} |S(\al_1,\al_2)|^2 d\al_1 d\al_2 = o(B^{2k-4}).
\eeq
It is here that the main innovation of our work lies. Note that the application of Parseval's identity in (\ref{Parseval})
has achieved several things: first, it has raised the weighted
exponential sum $S(\al_1,\al_2)$ to a higher power and homogenized the
problem, thus passing the problem of counting simultaneous solutions
of $Q_1(\xbf)=n_1$, $Q_2(\xbf)=n_2$ to the problem of counting
simultaneous solutions of \[Q_1(\xbf_1) - Q_1(\xbf_2) =0, \quad Q_2(\xbf_1) - Q_2(\xbf_2)=0,\]
where $\xbf_1, \xbf_2 \in \Z^k$.
Moreover, it has passed the power of $S(\al_1,\al_2)$ inside the
relevant integral, thus increasing the possibility for advantageous
averaging. In order to achieve this, we must treat the double integral
over $\al_1, \al_2$, and the possible interactions between $\al_1$
and $\al_2$, nontrivially. We accomplish this by developing a
two-dimensional Kloosterman refinement.

Historically, the Kloosterman refinement applies to a generating
function of a single variable, say \[ F(\al) = \sum_n r(n) e(\al n),\]
in which case one may be interested in computing the representation number
\[ r(0) = \int_0^1 F(\al) d\al.\]
To do so, one would traditionally apply the circle method to estimate
portions of the integral with $\al\approx a/q$ for certain rational
numbers $a/q$. Kloosterman's innovation \cite{Kloos} was to exploit cancellation
between the contributions corresponding to pairs of distinct rationals
$a_1/q, a_2/q$ with a common denominator $q$.

As suggested in \cite{HB96}, it would be desirable to apply a
Kloosterman refinement to two-dimensional problems, in which case one
would consider \[ \iint_{[0,1]^2}\sum_{m,n} r(m,n) e(\al_1 m + \al_2 n) d\al_1 d\al_2,\]
and hope to extract cancellation between portions of the integral
corresponding to neighbourhoods of distinct pairs of rationals $(a_1/q,
a_2/q)$ and $(b_1/q, b_2/q)$.  To do so, one would need to approximate
$\al_1, \al_2$ simultaneously by rationals with the same denominator,
in such a way that the intervals in the refinement all share the same
length. This can be accomplished by a 2-dimensional application of
Dirichlet's principle, which guarantees for any $S \geq 1$ the
existence of $1 \leq a_1,a_2 \leq q\leq S$, with $(a_1,a_2,q)=1$ such
that \beq\label{rsq} |\al_1 - a_1/q| \leq \frac{1}{q\sqrt{S}}, \qquad  |\al_2 - a_2/q|
\leq \frac{1}{q\sqrt{S}}.
\eeq
In the corresponding case in one dimension, Dirichlet's approximation
principle places $\al$ in an interval of length $(qS)^{-1}$, and one
can show that $\al$ can lie in at most two such intervals, if $q \leq
S$. In two dimensions, the intervals are longer, and a given pair
$(\al_1,\al_2)$ may lie in many such intervals simultaneously. This makes a true 2-dimensional Kloosterman refinement difficult to carry out.

The key feature in our application is that we apply the
2-dimensional Kloosterman refinement only to the minor arcs
contribution, for which we need simply an upper bound rather than an
asymptotic. Thus we include all approximations (\ref{rsq}) (accepting the
possible overlap and resulting loss in sharpness) and then carry out a
procedure to extract cancellation between the contributions of
distinct pairs of rationals. We note that while in a 1-dimensional
application of the Kloosterman refinement one typically encounters
exponential sums involving both $a$ and its inverse $\overline{a}$
modulo $q$, in our particular 2-dimensional application, such inverses
do not appear.  This is because in the usual 1-dimensional version one
has to take detailed account of the end points of the Farey arcs,
while in our 2-dimensional situation it suffices to use the simple
explicit squares (\ref{rsq}).

It is worth remarking that if one simply wanted to prove a result such
as Theorem \ref{thm_gen} for any $k > k_0$ sufficiently large one
could avoid the Kloosterman refinement. In particular, one could use the simpler methods of Birch \cite{Bir61} for the
treatment of the minor arcs for $k$ sufficiently large. But in order
to push the number of variables down to 5, we must use the more
technical argument presented in this paper.

\subsection{Notation}

We denote by $v_p(n)$ the $p$-adic order of $n$, and by $|n|_p =
p^{-v_p(n)}$ the $p$-adic valuation. We will use an underscore to
denote variables in 2-dimensional spaces, such as $\u{\al} \in \R^2$, $\u{a} \in \Z^2$. Similarly we will use boldface to denote variables in
$k$-dimensional spaces, such as $\xbf \in \Z^k$, and in
$2k$-dimensional spaces, such as $\lbf = (\lbf_1,\lbf_2) \in \Z^k
\times \Z^k$. We will use $Q(\xbf)$ to denote a quadratic form with
off-diagonal elements divisible by 2, as
well as $Q$ to denote the matrix associated with the quadratic form,
so that $Q(\xbf) = \xbf^t Q \xbf$; which meaning is intended will be
clear from context.  We will write $|\xbf|$ for the Euclidean length
of the vector $\xbf$, and $|\x|_p$ for the $p$-adic height.
We will use $\u{Q}$ to denote a pair of quadratic
forms $(Q_1,Q_2)$ and  $\u{a} \cdot \u{Q}$ to denote the linear
combination of forms $a_1Q_1 + a_2Q_2$. For a quadratic form $Q$ we
will let $||Q|| = \sup_{|\xbf|=1} |Q(\xbf)|$.  Throughout the paper
all constants, explicit or implied, will be allowed to depend on the forms $Q_1$ and $Q_2$, as well as on the
choice of the weight function $w$.

We say a weight function is smooth if it is $C^\infty$.   Denote by
$\partial_\xbf^\al$  the derivative with respect to $\xbf \in \R^k$ of order $\al$, where $\al$ is any multi-index
$\al=(\al_1,\ldots, \al_{k})$ of magnitude $|\al| = \al_1 + \cdots
+\al_{k}.$ We write
\[\frac{\partial^\al}{\partial \xbf^\al}\;\;\;\mbox{for}\;\;\;
\frac{\partial^{|\al|}}{\partial^{\al_1}x_1
 \cdots \partial^{\al_k}x_k}.\]
We use the notation $\iint_{\{\phi_1,\phi_2\}}$ to denote integration
over the region \[([-2\phi_1, -\phi_1]\cup [\phi_1, 2\phi_1]) \times ([-2\phi_2, -\phi_2]\cup [\phi_2, 2\phi_2]).\]

\section{Geometric conditions}
\subsection{Conditions related to Condition 1}\label{sec_geom_cond}
Recall Condition 1, which stated that the variety $V$ defined by
$Q_1(\xbf)=Q_2(\xbf)=0$ is nonsingular over $\overline{\Q}$. In
general, if we have any algebraically closed field $K$ of
characteristic zero or of
odd characteristic, we can consider the analogous
statement for $Q_1,Q_2$ over $K$.  We shall refer to this as
Condition 1 with respect to $K$.  Note that for any algebraically
closed field $K$ of characteristic zero,
the field of definition of the forms $Q_1,Q_2$ will
be taken to be $\Q$.
\label{note_Cond1}

We now fix such a field $K$ and define three more conditions. For convenience, we will refer to the relevant Jacobian matrix  as \[ J(\xbf) = \left( \begin{array}{c}
	 \nabla Q_1(\xbf) \\
	 \nabla Q_2(\xbf)
	 \end{array}
	 \right),\]
	 so that Condition 1 may be stated as $\rk(J(\xbf))=2$ for $\xbf \in V$ over $K$.
We also define the matrix
\beq\label{G_dfn}
\Delta(\xbf)=\Delta(\xbf;Q_1,Q_2)  = \left(\Delta_{ij}(\xbf)\right)_{i,j\le k} \eeq
where for each pair $i,j$ we write
\[ \Delta_{ij}(\xbf) = \left| \begin{array}{cc}
			\frac{\partial Q_1}{ \partial x_i} &
                       \frac{\partial Q_1}{ \partial x_j}\\
			\frac{\partial Q_2}{ \partial x_i} &
                       \frac{\partial Q_2}{ \partial x_j}
			\end{array} \right|\]
for the $i,j$-th minor of the Jacobian matrix $J(\xbf)$. Then $\rk(J(\xbf))=2$ precisely when $\Delta(\xbf) \neq 0$.

We begin by proving that Condition 1 is equivalent to two other
conditions.
Recall that we defined the determinant form
\[ F(x,y) = \det(xQ_1 + yQ_2).\]
This function will play a key role in the analysis throughout the
paper; note that $F$ is a binary form in $x$ and $y$ of degree $k$, and that its discriminant is a rational integer.
\begin{condition}\label{cond3}
The determinant form $F(x,y) = \det(xQ_1 + yQ_2)$ is not identically
zero and has distinct linear factors over $K$; that is to say the discriminant of $F(x,y)$ is nonzero in $K$.
\end{condition}

We say that $Q_1$ and $Q_2$ can by simultaneously diagonalized over
$K$ if there exists a basis $\{\ebf_1,\ldots, \ebf_k\}$ for
$K^k$ in which we may write \[Q_1(\sum x_i \ebf_i) = \sum_{i=1}^k
a_ix_i^2,\;\;\; Q_2(\sum x_i \ebf_i) = \sum_{i=1}^k b_ix_i^2,\]
with
coefficients  $a_i, b_i \in K$. \begin{condition}\label{cond4}
The forms $Q_1$ and $Q_2$ can be simultaneously diagonalized over
$K$, so that $Q_1 = \diag(a_i),$ $Q_2 = \diag(b_i)$ with $a_i, b_i
\in K$. Moreover, the ratios $a_i/b_i$ are well-defined as elements of
$K\cup\{\infty\}$ (that is to say $a_i$ and $b_i$ are not both zero), and are distinct, for $i=1,\ldots, k$. \end{condition}

Finally there is a fourth condition which is a consequence of the
three we have discussed, but not in general equivalent to them.

\begin{condition}\label{cond2}
For every coefficient pair $(\nu_1, \nu_2) \in K^2-\{0,0\}$, the
rank of the matrix associated to the quadratic form $\nu_1 Q_1 + \nu_2
Q_2$ satisfies \[ \rk (\nu_1Q_1 + \nu_2 Q_2) \geq k-1.\]
\end{condition}
Note that one cannot expect to have $\rk (\nu_1Q_1 + \nu_2 Q_2)=k$ for
every pair $(\nu_1,\nu_2)$.  For example, considering two
diagonal quadratic forms, it is clear that one may always choose a
linear combination that lowers the rank of the combination by one.  Thus
we cannot expect more than Condition 4 to hold, in general.

The main result of this section shows the relationships between these conditions; this equivalence has been
observed before, such as in Proposition 2.1 of \cite{Reid}; here we
will argue in more elementary terms. \begin{prop}\label{prop_cond}
Fix any algebraically closed field $K$ of characteristic zero or of
odd characteristic. Conditions 1, 2 and 3 with respect to $K$ are
equivalent, and they each imply Condition 4 with respect to $K$. Furthermore, Conditions 1, 2, and 3 with respect to $\Qbar$ are
equivalent to Conditions 1, 2, and 3 with respect to any algebraically
closed field $K$ of characteristic zero. \end{prop}
The last observation is visible when noting the equivalence of
Condition 2 with respect to $\Qbar$ and with respect to any $K$ of
characteristic zero, since the notion of $\disc(F) \neq 0$ is
equivalent in $\Q$ or in the prime field of $K$. In contrast, even if
Condition 1 holds with respect to $\Qbar$, Condition 1 may or may not
hold with respect to $\overline{\F_p}$ for each fixed odd prime; we
will distinguish between these possibilities in the next section when
we define good and bad primes.

We first prove that Condition 1 implies Condition 4, by showing the
contrapositive. Suppose that there exist $\nu_1,\nu_2 \in K$ for
which $\rk(\nu_1Q_1 + \nu_2 Q_2) \leq k-2$. Given such a pair, set
$Q_1' = \nu_1Q_1 + \nu_2Q_2$, and then choose $Q_2'$ to be any other
form such that the pair $\{Q_1',Q_2'\}$ generates the same
pencil as the original pair $\{Q_1,Q_2\}$.  Then by
assumption, \beq\label{rank'}
\rk(Q_1') \leq k-2.
\eeq
Via change of basis over $K$, we can diagonalize $Q_1' (\xbf) =
a_1x_1^2 + \cdots +a_{k-2}x_{k-2}^2,$ where we note that we may omit
the last two variables $x_{k-1},x_k$, due to the rank
assumption. After this change of basis, we also have some
representation for $Q_2'$, call it $Q_2'(x_1,\ldots, x_k).$ For any
$\xbf_0$ with $x_1=\cdots= x_{k-2}=0$, we obtain $Q_1'(\xbf_0)=0$ and
$\nabla Q_1'(\xbf_0)=0$, so that \beq\label{rk2}
\rk \left( \begin{array}{c}
	 \nabla Q_1'(\xbf_0) \\
	 \nabla Q_2'(\xbf_0)
	 \end{array}
	 \right)
	<2.
	 \eeq
So to show that Condition 1 fails, it is sufficient to find such
an $\xbf_0$ lying on the variety $Q_1'(\xbf)=Q_2'(\xbf)=0$. It is
automatic that $Q_1'(\xbf_0)=0$.  Moreover, when
$x_1=\cdots= x_{k-2}=0$ we see that $Q_2'(\xbf_0)$ is a quadratic form in the two
variables $x_{k-1},x_k$, so that it will have a non-trivial zero in $K$.  We
therefore obtain $\xbf_0 \neq \mathbf{0} \in K^k$ such that $Q_1'(\xbf_0) =
Q_2'(\xbf_0)=0$ and (\ref{rk2}) holds, contradicting Condition 1. We therefore conclude that Condition 1 implies Condition 4.

We next show that Condition 1 also implies Condition 2.  We begin by
observing that a change of basis via some invertible matrix $M$ converts the matrices $Q_1,Q_2$ into
$M^tQ_1M, M^tQ_2M$, and therefore multiplies
$F(x,y)$ by the non-zero constant $(\det M)^2$.  Moreover, replacing $Q_1,Q_2$ by another pair of forms generating the same
pencil has the effect of making an invertible linear substitution in
the variables $x,y$ occurring in $F(x,y)$.  The properties described in
Conditions 1 and 2 are clearly unchanged by these two types of
transformations.  Now, to prove our assertion we will again argue by
contradiction.  If Condition 2 fails we can make a linear change
between $Q_1$ and $Q_2$ so that $y^2|F(x,y)$. Thus
$\det(Q_1)=0$, since this is the coefficient of $x^k$ in $F(x,y)$.  However, since we now know that Condition 1 implies
Condition 4, we see that ${\rm rk}(Q_1)\ge k-1$ whence in fact ${\rm rk}(Q_1)=k-1$. After a change of basis we may then write
$Q_1(\xbf)$ as $Q_1'(x_1,\ldots,x_{k-1})$, and $Q_2(\xbf)$ as \[Q_2'(x_1,\ldots,x_{k-1})+(c_1x_1+\ldots+c_{k-1}x_{k-1})x_k+c_kx_k^2,\]
say. One then sees that the coefficient of $x^{k-1}y$ in
$\det(xQ_1 + yQ_2)$ must be $\det(Q_1')c_k$.  This coefficient must in fact
vanish in view of the condition $y^2|F(x,y)$. However $\det(Q_1')\not=0$
since ${\rm rk}(Q_1)=k-1$, and we therefore deduce that $c_k=0$.  The
point $\xbf_0=(0,\ldots,0,1)$ will therefore satisfy $Q_1'(\xbf_0)=Q_2(\xbf_0)=0$ and also $\nabla Q_1'(\xbf_0)=\mathbf{0}$.
We therefore have a singular point on the variety $V$, contradicting
Condition 1.  This suffices to establish our assertion that Condition
1 implies Condition 2.
\label{note_cond12}

\label{note_cond23}
We next assume that Condition 2 holds and deduce Condition 3. Since
the field $K$ is infinite Condition 2 implies that at least one linear
combination of $Q_1$ and $Q_2$ is nonsingular, and by a linear change
between the two forms we may assume that $\det(Q_2)\not=0$. Such a
linear change does not affect the validity or otherwise of Condition 3.
For each of the $k$ distinct roots $\lam_i$ of the equation $\det (Q_1 - \lam
Q_2)=0$, there exists an $\x_i \neq 0$, $\x_i \in K^k$,
such that \beq\label{eqnQ} Q_1 \x_i = \lam_i Q_2 \x_i.
\eeq
Note that as the $\x_i$ are eigenvectors for $Q_2^{-1}Q_1$,
corresponding to distinct eigenvalues, they are linearly independent,
and hence form a basis of $K^k$.
Taking the dot product of equation (\ref{eqnQ}) with any basis vector $\xbf_j$
we obtain
\beq\label{Mji}
\x_j^t Q_1 \x_i = \lam_i \x_j^t Q_2 \x_i.
\eeq
But the transpose of a real number is itself, and $Q_1, Q_2$ are
symmetric matrices, so \[ \x_j^t Q_1 \x_i = (\x_j^t Q_1 \x_i)^t = \x_i^t Q_1 \x_j = \lam_j
\x_i^t Q_2 \x_j = ( \lam_j \x_i^t Q_2 \x_j )^t = \lam_j \x_j^t Q_2
\x_i.\] Thus for any $i,j$,
\[  \lam_i ( \x_j^t Q_2 \x_i)=  \lam_j (\x_j^t Q_2 \x_i),\]
and so by the assumption that the $\lam_i$ are distinct we have
\[ \x_j^t Q_2 \x_i =0,\]
for any pair $i \neq j$. Hence by (\ref{Mji}), \[ \x_j^t Q_1 \x_i =0.\]
This means that the basis $\x_1, \ldots, \x_k$ diagonalizes both
matrices simultaneously, as desired. In the notation of Condition 3 we
will have $a_i=\x_i^t Q_1 \x_i$ and $b_i=\x_i^t Q_2 \x_i$ so that
$a_i/b_i=\lam_i$. We therefore see that these ratios are distinct as
required.  This completes the proof that Condition 2 implies Condition 3.
\label{note_23}

Finally, we show that Condition 3 implies Condition 1. We may assume
that $Q_1=\diag(a_i)$, $Q_2 = \diag(b_i)$ have been diagonalized. Let
$\mathbf{t} \neq \mathbf{0}$ be any point such that $Q_1(\mathbf{t}) =
Q_2(\mathbf{t})=0$. Suppose that $\rk(J(\mathbf{t}))<2$, so that there exists
$\lam \in K\cup\{\infty\}$ such that \beq\label{ab}
\nabla Q_1(\mathbf{t}) = \lam \nabla Q_2(\mathbf{t}).
\eeq We therefore see that $a_i t_i = \lam b_i t_i $
for all $i=1,\ldots, k$. Since $\mathbf{t} \neq \mathbf{0}$, there exists at least
one index $i$ for which $t_i \neq 0$ and hence we may deduce $\lam =
a_i/b_i$. If there are at least two indices $i \neq j$ with $t_i \neq
0$, $t_j \neq 0$, then $\lam = a_i/b_i$ and $\lam = a_j/b_j$, which
contradicts the condition that the ratios are distinct. Thus there can
only be one nonzero coordinate of $\mathbf{t}$, which we may assume is
$t_1$.  Then $Q_1(\mathbf{t})=a_1t_1^2$ and $Q_2(\mathbf{t})=b_1t_1^2$.
But we also assumed that
$Q_1(\mathbf{t}) = Q_2(\mathbf{t}) = 0$, and since $a_1$ and $b_1$
cannot both vanish in Condition 3, we obtain a contradiction. Hence Condition 3 implies Condition 1.
This completes the proof of Proposition \ref{prop_cond}.

\subsection{Definition of the good and bad primes}\label{sec_good_bad}
Recall from Proposition \ref{prop_cond} that under the assumption
of Condition 1, the determinant form $F$ has distinct linear factors over
$\overline{\Q}$. In particular, if we set
\[D_F:={\rm Disc}(F),\]
then $D_F$ is a non-zero rational integer.
\label{note_disc}\label{note_factorF}
Furthermore, if we write $K$ for the splitting field for $F$ over $\Q$, then we can factor $F$ as
\[F(x,y)=c^{-1}\prod_{i=1}^k(\lambda_i x-\mu_i y)\]
with $c\in\N$ and $\lambda_i,\mu_i\in\mathcal{O}_K$. \begin{dfn}
We shall say that a prime $p$ is ``bad'' if
\begin{equation}\label{rbad}
p\mid 2cD_F\prod_{\sigma}
\prod_{1\le i<j\le k}(\lambda_i\mu_j-\lambda_j\mu_i)^{\sigma},
\end{equation}
where $\sig$ runs over the Galois automorphisms of $K/\Q$. Otherwise we shall say that $p$ is ``good.''
\end{dfn}
\label{note_goodprimes}
It is an immediate consequence of this definition that  the set of bad
primes is finite, and is determined by the original forms $Q_1$ and $Q_2$.
Moreover, if $p$ is good, then the system $Q_1 = Q_2 =0$ is
nonsingular over $\overline{\F_p}$, in the sense that Condition 1
holds relative to $\overline{\F_p}$. For indeed, if this system is
singular over $\overline{\F_p}$ (with $p$ odd), then Condition 1 and
hence Condition 2 fails relative to $\overline{\F_p}$, so that $D_F=0$ in $\overline{\F_p}$, whence $p$ satisfies
(\ref{rbad}) and hence is bad.

Similarly, we will later call upon the following local version of
Condition \ref{cond2}: \begin{lemma}\label{lemma_Delta_p}
If $\rk(a_1Q_1 + a_2Q_2) \leq k-2$ over $\F_p$ for some $a_1,a_2 \in
\Z$ with $(a_1,a_2,p)=1$, then $p$ is bad. \end{lemma}

To establish this we note that if $p$ is odd,
the existence of such $a_1,a_2$ would
mean that Condition 4 fails for $K=\overline{\F_p}$.  By
Proposition \ref{prop_cond}, Condition 2 also fails for
$K=\overline{\F_p}$, and hence $D_F=0$ in $K=\overline{\F_p}$.  The
lemma then follows.
\label{note_Cond4local}

\subsection{Definition of Type I and Type II
 primes}\label{sec_local_smooth}

In our consideration of the singular series associated to a given pair
$(n_1,n_2)$, we will require an affine smoothness condition.
Given any pair $(n_1,n_2)\in \Z^2$, and any prime $p$, set \beq\label{Vp_dfn}
V_p(n_1,n_2) = \{ \xbf \in \overline{\F_p}^k : Q_1(\xbf) = n_1, Q_2(\xbf) = n_2 \; \text{in} \;\overline{\F_p} \}.
\eeq

\begin{dfn}
We will say that a prime $p$ is of ``Type I'' with respect to a fixed
pair of values $(n_1,n_2)$ if $p$ is good and $V_p(n_1,n_2)$ is
nonsingular over $\overline{\F_p}$, in the sense that $\Delta(\xbf) \neq 0$ in
$\overline{\F_p}$, for all $\xbf \in V_p(n_1,n_2)$. We will say that a prime $p$
is of ``Type II'' with respect to $(n_1,n_2)$ if $p$ is good but $p$
is not of Type I. \end{dfn}

It is worth remarking that while the distinction between good and bad
primes is completely independent of any values $(n_1,n_2)$, the
distinction between Type I and Type II primes is always with respect
to a fixed pair $(n_1,n_2)$. Note that the condition that $p$ is a Type I prime is a local version
of the requirement that the Jacobian matrix $J(\xbf)$ be full rank. On the other hand, if $p$ is a Type II prime, there exists $\xbf_0 \in
V_p(n_1,n_2)$ for which there is a pair $\al, \be$
with $(\al,\be)  \neq (0,0)$ in $\F_p$ that satisfies \beq\label{bad_p1}
\al \nabla Q_1(\xbf_0) + \be \nabla Q_2(\xbf_0) \con 0 \modd{p} .
\eeq
Such an $\xbf_0$ must be nonzero, since $p$ is good.

In order to prove the convergence of the singular series it will be crucial that for any pair $(n_1,n_2)$ we consider,
the number of Type II primes with respect to $(n_1,n_2)$ is finite.

\begin{lemma}\label{lemma_T_II}
Given $n_1,n_2$ (not both zero) such that $n_2Q_1(\xbf) -
n_1Q_2(\xbf)$ is globally nonsingular,  finitely many primes are of
Type II with respect to $n_1,n_2$. Indeed, any Type II prime must
satisfy $p|F(n_2,-n_1)$. \end{lemma}
\label{note_n1n2lemma}

Suppose $p$ is a Type II prime with respect to $n_1,n_2$ so that there
exists a non-zero point $\xbf_0  \in V_p(n_1,n_2)$ for which
(\ref{bad_p1}) holds. Then taking the dot product of $\xbf_0$ with
(\ref{bad_p1}) yields \beq\label{bad_p3}
2\al Q_1(\xbf_0) + 2\be Q_2(\xbf_0) \con 0 \modd{p}.
\eeq
But by assumption $\xbf_0$ lies on the variety $Q_1(\xbf_0)-n_1 \con
Q_2(\xbf_0)-n_2 \con 0$, so that (\ref{bad_p3}) implies
\[ 2\al n_1 + 2\be n_2 \con 0 \modd{p}.\]
Recall that $2$ is a bad prime, whereas Type II primes are good
primes, so $p \neq 2$. In the case that $p|\gcd(n_1,n_2)$, we have
$p|F(n_2,-n_1)$ and we are
finished. Otherwise, we may assume that $p \ndiv n_1$, say, and solve
for $\al \con -\be n_2 \bar{n}_1$. Hence in (\ref{bad_p1}),
\[-\be n_2 \bar{n}_1 \nabla  Q_1(\xbf_0) + \be \nabla Q_2(\xbf_0)
\con 0 \modd{p}.\]
If $p |\be$ then $\al \con -\be n_2 \bar{n}_1\con
0\modd{p}$, contradicting the fact that $(\al,\be)  \neq
(0,0)$. Hence $p\nmid\be$ and we conclude that
\beq\label{bad_p4}
 n_2 \nabla Q_1(\xbf_0) \con n_1 \nabla Q_2(\xbf_0) \modd{p}.
 \eeq
Regarding $Q_1,Q_2$ as matrices, (\ref{bad_p4}) is equivalent to the statement
\beq\label{matrix}
 (n_2Q_1 - n_1Q_2)\xbf_0 \con 0 \modd{p}.
 \eeq
However by assumption $\xbf_0 \not\con 0 \modd{p}$, whence
(\ref{matrix}) implies that $p$ divides $\det(n_2Q_1-n_1Q_2)$. Given a pair of values $n_1,n_2$, there are finitely many such primes
$p$, unless the determinant vanishes as an element of $\Z$.  However this
is precisely the
condition we have ruled out by assuming that $n_2Q_1(\xbf) -
n_1Q_2(\xbf)$ is globally nonsingular so that $\det(n_2Q_1 - n_1Q_2) \neq0$. This proves the lemma.

We remark that when $n_2Q_1 - n_1Q_2$ is globally singular there
may be infinitely many primes of Type II. For example, suppose that
$Q_1$ is itself singular, so that $\nabla Q_1(\xbf_0) = 0$ for some non-zero
integer vector $\xbf_0$.  Taking $n_1=0$ and $n_2=Q_2(\xbf_0)$ we then
see that $n_2Q_1(\xbf) - n_1Q_2(\xbf)$ is globally singular and
$\xbf_0$ is a singular point of $V_p(n_1,n_2)$ for every prime $p$.
Thus in our consideration of the major arcs, we restrict to those
pairs $(n_1,n_2)$ for which $\det(n_2Q_1 - n_1Q_2)\neq 0$,
which is the condition $F(n_2,-n_1)\not=0$ in Theorem \ref{thm_r1}.

In connection with Theorem \ref{thm_s2} we note that if $(n_1,n_2)=(0,0)$ there are no Type I primes, since $\xbf=\mathbf{0}$
is always a singular point on $V_p(0,0)$.

\subsection{Bounds for eigenvalues}\label{sec_eigen}
Given any pair $(\nu_1,\nu_2) \in \R^2$, Condition \ref{cond2}
allows that $\u{\nu} \cdot \u{Q}$ may be of rank $k-1$ and hence
$F(\nu_1,\nu_2)$ may vanish. Nevertheless, as we prove in the
following lemma, under Condition \ref{cond2}, at most one of the
eigenvalues of $\u{\nu} \cdot \u{Q}$ may be small.

\begin{lemma}\label{lemma_rho}
Let $\nu^* = \max( |\nu_1|,|\nu_2|)$, and let $\rho_1, \ldots,\rho_k$ denote the eigenvalues associated to  $\u{\nu} \cdot \u{Q}$, ordered
so that $|\rho_1|\leq\cdots\leq|\rho_k|$.  Then, under Condition \ref{cond2}, we have
\[ |\rho_2|\gg\nu^* \;\;\;\mbox{and}\;\;\; |\rho_k|\ll\nu^*.\]
\end{lemma}
We may assume that $|\nu_2| \leq 1 = |\nu_1|$ (by normalizing and
possibly interchanging the roles of $\nu_1,\nu_2$, $Q_1,Q_2$).  We
therefore study a linear combination of the type
$Q_1 + \lam Q_2$ with $|\lam| \leq 1$. Then \[|\rho_k|\le||Q_1+\lam Q_2||\ll_{Q_1,Q_2}1,\]
so that it remains to prove that
\[ |\rho_2| \gg_{Q_1,Q_2} 1.\]
We will argue by contradiction.  It will be convenient to write
$\rho_i(Q)$ to denote the $i$-th smallest eigenvalue of a real quadratic
form $Q$. We will then assume for a contradiction that for any positive integer
$n$ there is a $\lam_n\in[-1,1]$ for which $|\rho_2(Q_1+\lam_nQ_2)|\le 1/n$.
We can diagonalize $Q_1+\lam_nQ_2$ using an orthogonal matrix $M_n$,
say, so that $M_n^T(Q_1+\lam_nQ_2)M_n={\rm
 diag}(\rho_1^{(n)},\ldots,\rho_k^{(n)})$, say, with
\[|\rho_1^{(n)}|\le\ldots\le|\rho_k^{(n)}|.\]
By construction, both $\rho_1^{(n)}$ and $\rho_2^{(n)}$ tend to zero as $n$ goes to infinity.
The set of orthogonal matrices of order $k$ is compact, as is the
interval $[-1,1]$.  \label{note_compactarg}
Hence by choosing a suitable
subsequence we can suppose that $\lambda_n\rightarrow\lambda$ and
that $M_n\rightarrow M$, say.  We then see that $M^T(Q_1+\lam Q_2)M={\rm
 diag}(\rho_1,\ldots,\rho_k)$, say, with $\rho_1=\rho_2=0$.  This
contradicts Condition 4 and hence proves the lemma.

\section{Bounds for oscillatory integrals}\label{sec_osc_int}
In considering both the major and minor arcs, we will require
estimates for oscillatory integrals of the form \beq\label{J_dfn}
I(\Qcal;\lambf) = \int _{\R^{n}} e( \Qcal(\ubf) - \lambf \cdot \ubf) w(\ubf) d\ubf ,
\eeq
where for the moment $\lambf \in \R^n$, $\Qcal$ is any real quadratic form in
$n$ variables and $w$ is any smooth weight on $\R^n$ supported in
$[-1,1]^{n}$. \begin{lemma}\label{I_int}
Let $\rho_1, \ldots, \rho_n$ be the eigenvalues associated to the
quadratic form $\Qcal$. If $|\lambf| \geq 4||\Qcal||$, then
\beq\label{lam_bd1}
I(\Qcal;\lambf) \ll_{M,w} |\lambf|^{-M},
\eeq
for any $M \geq 1$. Moreover
\beq\label{lam_bd2}
|I(\Qcal;\lambf)| \ll_{w}\prod_{i=1}^{n}\min (1,\frac{1}{|\rho_i|^{1/2}}).
\eeq
\end{lemma}

We will estimate $I(\Qcal;\lambf)$ by the method of stationary phase.
Note first that if $|\lambf| \geq 4||\Qcal||$, then \[| \nabla_u \{\Qcal(\ubf) - \lambf \cdot \ubf\}| = |\nabla_u \Qcal(\ubf)
- \lambf| \geq|\lambf|/2\] on ${\rm supp}(w)$, since $|\nabla_u \Qcal(\ubf)|\le
2 ||\Qcal||\le|\lambf|/2$.  
On the other hand, for any multi-index $\al$ with $|\al|=2$ we have, \[ |\partial_u^\al \{\Qcal(u) -\lambf \cdot u\}| =  |\partial_u^\al
\Qcal(u) |\leq 2C ||\Qcal|| \leq \frac{C}{2} |\lambf|,\] for some constant $C$ depending on $n$. Moreover, when $|\al|\ge 3$
the left-hand side vanishes.
Thus an application of the first derivative test for an infinitely
differentiable function in high dimensions, as presented in Lemma 10
of \cite{HB96}, shows that $I(\Qcal;\lambf) \ll |\lambf|^{-M},$
for any $M \geq 1$, which proves (\ref{lam_bd1}).

To prove our second estimate we
will apply the second derivative test.
 Let $R$ be an orthogonal transformation such that $R^t\Qcal R = D$,
 where $D =\diag\{\rho_1,\ldots, \rho_n\}$. Under this change of variables,
\[ I(\Qcal;\lambf) = \int _{\R^{n}} e( \sum_{i=1}^n \rho_i u_i^2 - (R^t
\lambf) \cdot \ubf) w(R^{}\ubf) d\ubf.\]

Next, we apply the following lemma, in order to remove the presence of the weight $w$.
\begin{lemma}\label{lemma_weight}
Let $f,w$ be smooth functions of $\R^n$ and suppose that $w$ is
supported on $[-1,1]^n$. Then \[ | \int _{\R^{n}} f(\xbf) w(\xbf) d\xbf |\le \left\{\int_{\R^{n}}|\hat{w}(\ybf)|d\ybf\right\}
\sup_{\ybf\in\R^n}\left|\int_{[-1,1]^n}f(\xbf)
e(\xbf\cdot\ybf)d\xbf\right|,\]
where $\hat{w}$ is the Fourier transform of $w$.
\end{lemma}

To prove the lemma we express $w(\xbf)$ in terms of its Fourier
transform as
\[w(\xbf)=\int_{\R^n}\hat{w}(\ybf)e(\xbf\cdot\ybf)d\ybf,\]
where $\hat{w}$ is smooth and of rapid decay.
Using an interchange in the orders of integration we then find that
\begin{eqnarray*}
\int_{\R^{n}} f(\xbf)\chi_{[-1,1]^n}(\xbf)w(\xbf) d\xbf
&=&\int_{\R^{n}}f(\xbf)\chi_{[-1,1]^n}(\xbf)
\int_{\R^n}\hat{w}(\ybf)e(\xbf\cdot\ybf)d\ybf d\xbf\\
&=&\int_{\R^{n}}\hat{w}(\ybf)\int_{\R^{n}}f(\xbf)\chi_{[-1,1]^n}(\xbf)
e(\xbf\cdot\ybf)d\xbf d\ybf,
\end{eqnarray*}
and the lemma follows.

Continuing our treatment of Lemma \ref{I_int}, we now see that we can bound
$I(\Qcal;\lambf)$ using an $n$-fold product of one-dimensional integrals:
\[ |I(\Qcal;\lambf)| \ll_w \prod_{i=1}^{n} \sup_{y\in\R}
\left| \int_{-1}^1 e(\rho_i u^2 +yu) du \right|.\]
Applying the second derivative test, each one-dimensional integral is
bounded by $|\rho_i|^{-1/2}$.  Alternatively we may use the trivial bound for
the integral, and clearly (\ref{lam_bd2}) follows.

We will apply Lemma \ref{I_int} to the situation in which $\Qcal$ takes the
shape $\u{\nu} \cdot \u{F} = \nu_1F_1 + \nu_2 F_2$ for two fixed
quadratic forms $F_1,F_2$.  In this case we write \[ I(\u{\nu} \cdot \u{F};\lambf) = \int _{\R^{n}} e( \u{\nu} \cdot
\u{F}(\ubf) - \lambf \cdot \ubf) w(\ubf) d\ubf.\] We now employ Lemma \ref{I_int} to estimate the average of $|I(\u{\nu}
\cdot \u{F};\lambf)|$ over dyadic
ranges of $\nu_1, \nu_2$, using the notation
$\iint_{\{\phi_1,\phi_2\}}$ to denote integration over the region
\[([-2\phi_1, -\phi_1]\cup [\phi_1, 2\phi_1]) \times ([-2\phi_2, -\phi_2]\cup [\phi_2, 2\phi_2])\]
when $\phi_1,\phi_2>0$. \begin{lemma}\label{I_int_avg}
Let $\phi^* = \max(\phi_1,\phi_2)$. Take $F_1=Q_1$, $F_2=Q_2$ to be
the original quadratic forms acting on $\Z^k$. Then \[\iint_{\{ \phi_1,\phi_2\} } |I(\u{\nu} \cdot \u{F};\lambf) | d \u{\nu}
\ll \min((\phi^*)^2,(\phi^*)^{2-k/2} ).\] \end{lemma}

This will be applied to prove the convergence of the
singular integral on the major arcs. When we consider the oscillatory
integral on the minor arcs we will need a version related to a
different choice of quadratic forms acting on $\Z^{2k}$: \begin{lemma}\label{I_int_avg2}
Let $\phi^* = \max(\phi_1,\phi_2)$. Take
\[F_1(\xbf_1,\xbf_2)=Q_1(\xbf_1) - Q_1(\xbf_2),\;\;\;
F_2(\xbf_1,\xbf_2)=Q_2(\xbf_1) - Q_2(\xbf_2),\]
for $\xbf_1, \xbf_2\in \Z^k$. Then \[\iint_{\{ \phi_1,\phi_2\} } |I(\u{\nu} \cdot \u{F};\lambf) | d \u{\nu}
\ll \left\{\begin{array}{cc} (\phi^*)^2, & \phi^*\le 1,\\
(\phi^*)^{2-k}(1+\log \phi^*), & \phi^*\ge 1.\end{array}\right.\] \end{lemma}

\subsection{Proof of Lemmas \ref{I_int_avg} and \ref{I_int_avg2}}
For each fixed $\u{\nu}=(\nu_1,\nu_2)$, let $\nu^* = \max (|\nu_1|,
|\nu_2|)$ and let  $\rho_1,\ldots,\rho_k$ be the
eigenvalues associated to the quadratic form $\u{\nu} \cdot \u{Q}$,
ordered so that $|\rho_1|\leq\cdots\leq|\rho_k|$.
Recall from Lemma \ref{lemma_rho} that for each fixed $\u{\nu}$, at
most one eigenvalue of $\u{\nu} \cdot \u{Q}$ may be of smaller order than $\nu^*$, and that
for $i=2, \ldots, k$ we have $|\rho_i| \approx \nu^*$, where the implied constant
depends on the initial forms $Q_1, Q_2$.
Applying this in (\ref{lam_bd2}), it follows that
\beq\label{IF}
|I(\u{\nu} \cdot \u{F};\lambf) | \ll \min (1,
(\frac{1}{\nu^*})^{\frac{k-1}{2}}) \min(1,
\frac{1}{|\rho_1|^{1/2}}).
\eeq
Recalling the definition $F(\nu_1, \nu_2) = \det(\nu_1 Q_1 + \nu_2
Q_2),$ it follows that \[ \rho_1 \cdots \rho_k = F(\nu_1, \nu_2).\]
Since in the range of the
integral we have $|\u{\nu}| \ll \phi^*$, it follows in particular that $\rho_i
\ll \phi^*$ for $i=2, \ldots, k$, and hence \[ |\rho_1| \gg \frac{|F(\nu_1,\nu_2)|}{(\phi^*)^{k-1}}.\]
We therefore deduce that
\begin{eqnarray}\label{I_int_phi}
\lefteqn{\iint_{\{ \phi_1,\phi_2\}}  |I(\u{\nu} \cdot \u{F};\lambf) | d \u{\nu}}\hspace{1cm}\nonumber\\
&\ll&  \min(1, (\frac{1}{\phi^*})^{\frac{k-1}{2}}) \iint_{\{
 \phi_1,\phi_2\}} \min (1,
\left(\frac{(\phi^*)^{k-1}}{|F(\nu_1,\nu_2)|}\right)^{1/2} )d\u{\nu}. \end{eqnarray}
According to Condition 3 we can factor $F(\nu_1,\nu_2)$ over $\C$ as
\beq\label{dp}
F(\nu_1,\nu_2) = \prod_{i=1}^k (a_i \nu_1 - b_i \nu_2),
\eeq
with distinct ratios $a_i/b_i$ in $\C\cup\{\infty\}$.  We will fix an
admissible choice for the coefficients $a_i,b_i$ once and for all.
Write $\psi_i=a_i \nu_1 - b_i \nu_2$ and order the indices so that
\beq\label{Rorder}
|\psi_1|\le |\psi_2|\le\cdots. \eeq
Since
\[\nu_1=\frac{b_2\psi_1-b_1\psi_2}{a_1b_2-a_2b_1}\]
and similarly for $\nu_2$, we conclude that $\nu^*\ll |\psi_2|$, so
that $|\psi_2|\gg\phi^*$ and hence
\beq\label{Flb} |F(\nu_1, \nu_2)| \gg (\phi^*)^{k-1}|\psi_1|.
\eeq
Certainly $|\psi_1| \geq \min_i |a_i\nu_1-b_i \nu_2|$. Thus, \[ \left(\frac{(\phi^*)^{k-1}}{|F(\nu_1,\nu_2)|}\right)^{1/2}  \leq
\frac{1}{\min_i |a_i \nu_1 - b_i \nu_2|^{1/2}}  \leq \sum_{1 \leq i
 \leq k} \frac{1}{|a_i \nu_1 - b_i \nu_2|^{1/2}} . \] Integrating over the appropriate region gives us an upper bound for
the integral in (\ref{I_int_phi}): \begin{eqnarray*}
\lefteqn{ \iint_{\{\phi_1,\phi_2\}} \min (1,
\left(\frac{(\phi^*)^{k-1}}{|F(\nu_1,\nu_2)|}\right)^{1/2} )d\u{\nu}}
\hspace{2cm}\\
&\ll &\sum_{1 \leq i \leq k} \iint_{\{\phi_1,\phi_2\}} \min (1,
\frac{1}{|a_i \nu_1 - b_i \nu_2|^{1/2}} )d\u{\nu}.
\end{eqnarray*}
For each fixed index $i$, we may exchange the roles of $Q_1$ and $Q_2$
if necessary so that $a_i \neq 0$, and set $c_i=b_i/a_i$.
Then the contribution of the integral corresponding to the index $i$ is majorized by \[ \iint_{\{\phi_1,\phi_2\}} \min(1, \frac{1}{|\nu_1-c_i \nu_2|^{1/2}}) d\u{\nu},
\]
which after the change of variables $v_1 = \nu_1-c_i\nu_2,$
$v_2 = \nu_2$, is the sum of four integrals of the type
\begin{eqnarray*}
\int_{\phi_2}^{2\phi_2} \int_{\phi_1-c_iv_2}^{2\phi_1-c_iv_2} \min(1, \frac{1}{|v_1|^{1/2}}) dv_1 dv_2  &\le&
2\phi^*\int_0^{(2+|c_i|)\phi^*}\min(1, \frac{1}{v^{1/2}}) dv\\
&\ll&\phi^*\min(\phi^*\,,\,(\phi^*)^{1/2}).
\end{eqnarray*}
Incorporating this upper bound in (\ref{I_int_phi}) completes the
proof of Lemma \ref{I_int_avg}.

Lemma \ref{I_int_avg2} is proved in the same manner as Lemma
\ref{I_int_avg}, with one important difference. Again let $\rho_1,
\rho_2, \ldots , \rho_k$ be the eigenvalues of $\u{\nu} \cdot \u{Q}$
acting on $\Z^k$, where $\u{Q} = (Q_1,Q_2)$. Then the particular
choice of $\u{F}$ in Lemma \ref{I_int_avg2} means that the eigenvalues
of $\u{\nu}\cdot \u{F}$ occur in pairs, as $\pm \rho_1, \pm
\rho_2,\ldots, \pm \rho_k$. In particular, \emph{two} of the
eigenvalues can now be small. However, this is easy to handle, as we
simply replace (\ref{IF}) with \[  |I(\u{\nu} \cdot \u{F},\lambf) | \ll \min (1,
(\frac{1}{\nu^*})^{k-1} ) \min(1, \frac{1}{|\rho_1|}),\] and the argument proceeds along the same lines as in the previous case.

\subsection{The singular integral}
We end this section by considering the singular integral, defined by
\beq\label{sing_int_dfn}
\Jcal_w(\u{\mu}) = \iint_{\R^2} \int_{\R^k} e( \u{\theta}
\cdot (\u{Q}(\xbf) -\u{\mu}))w(\xbf) d\xbf d\theta_1 d\theta_2. \eeq
We also define the truncated singular integral as
\beq\label{sing_int_dfn+}
\Jcal_w(\u{\mu};R) =\int_{-R}^R \int_{-R}^R \int_{\R^k} e( \u{\theta} \cdot
( \u{Q}(\xbf)-\u{\mu}) )w(\xbf) d\xbf d\theta_1 d\theta_2,
\eeq
so that $\Jcal_w (\u{\mu})= \lim_{R \maps \infty} \Jcal_w(\u{\mu};R)$, if the limit exists. We now prove the following proposition.
\begin{prop}\label{prop_sing_int_conv}
For $k\ge 5$, the singular integral $\Jcal_w(\u{\mu})$ is absolutely
convergent in the sense that
\[\iint_{\R^2} \left|\int_{\R^k} e( \u{\theta}
\cdot (\u{Q}(\xbf) -\u{\mu}))w(\xbf)d\xbf \right| d\theta_1 d\theta_2<\infty.\] Moreover it is bounded uniformly in $\u{\mu}$. The
rate of convergence may be quantified for $R\ge 2$ as \[ |\Jcal_w(\u{\mu}) - \Jcal_w(\u{\mu};R)| \ll R^{\frac{4-k}{2}}\log R.\]
\end{prop}

We first verify that $\Jcal_w(\u{\mu})$ converges. Letting $\theta^* = \max (|\theta_1|, |\theta_2| )$, we may write
\[ |\Jcal_w(\u{\mu}) - \Jcal_w(\u{\mu};R) | \ll \iint_{\theta^* > R}
|I(\u{\theta})|d \u{\theta},\] where \beq\label{dfn_I}
I( \u{\theta})  = \int_{\R^k} e(\u{\theta} \cdot \u{Q}(\xbf))  w(\xbf) d\xbf.
\eeq
Lemma \ref{I_int_avg} implies that \beq\label{R+}
\iint_{\{ \phi_1,\phi_2\}} |I(\u{\theta}) |  d \u{\theta}  \ll
(\phi^*)^{\frac{4-k}{2}},
\eeq
for any dyadic range $\{ \phi_1,\phi_2\}$ with $\phi^* =
\max(\phi_1,\phi_2) \geq 1$.  Since $I(\u{\theta})\ll 1$ we also have
\[\iint_{\{ \phi_1,\phi_2\}} |I(\u{\theta}) |  d \u{\theta}  \ll
\phi_1\phi_2.\]
We now sum over dyadic ranges for $\phi_1$ and $\phi_2$, using this
latter bound when $\min(\phi_1,\phi_2)\le (\phi^*)^{-k}$, and (\ref{R+})
otherwise.  This shows that
\begin{eqnarray}\label{plane}
\iint_{\theta^* > R}|I(\u{\theta})|d \u{\theta}& \ll & \sum_{2^n>R}\left\{\sum_{\substack{m\in\Z \\ m\le -kn}}2^m.2^n+
\sum_{\substack{m\in\Z \\ -kn<m\le n}}(2^n)^{(4-k)/2}\right\}\nonumber\\
&\ll & \sum_{2^n>R}\left\{2^{n(1-k)}+n(2^n)^{(4-k)/2}\right\}\nonumber\\
&\ll & R^{\frac{4-k}{2}}(1+ \log R).
\end{eqnarray}
Thus $\Jcal_w(\u{\mu};R)$ converges absolutely for $k\ge 5$. A
similar argument shows that $\Jcal_w(\u{\mu})$ is uniformly bounded with
respect to $\u{\mu}$.
\label{note_singint}

Our second major result on the singular integral gives an interpretation
in terms of the density of real solutions of $\u{Q}(\xbf)=\u{\mu}$.
\label{note_Kep}

\begin{prop}\label{prop_sing_int_val}
Let $k\ge 5$.  Then \[\ep^{-2}\int_{\max|Q_i(\xbf)-\mu_i|\le\ep}w(\xbf)
\left(1-\frac{|Q_1(\xbf)-\mu_1|}{\ep}\right)
\left(1-\frac{|Q_2(\xbf)-\mu_2|}{\ep}\right)d\xbf\]
tends to $\Jcal_w(\u{\mu})$
as $\ep\rightarrow 0$.
\end{prop}
We remark for later reference that Propositions \ref{prop_sing_int_conv} and \ref{prop_sing_int_val} hold even when $\u{\mu}=\u{0}$.

For the proof we define
\[K_{\ep}(\u{\theta})=
\left(\frac{\sin\pi\ep\theta_1}{\pi\ep\theta_1}\right)^2
\left(\frac{\sin\pi\ep\theta_2}{\pi\ep\theta_2}\right)^2.\]
Then $K_{\ep}(\u{\theta})=1+O(\ep^{3/2})$ if
$\max(|\theta_1|,|\theta_2|)<\ep^{-1/4}$ and $K_{\ep}(\u{\theta})\ll
1$ in general.  It follows that
\begin{eqnarray*}
\lefteqn{\iint_{\max(|\theta_1|,|\theta_2|)<\ep^{-1/4}}|K_{\ep}(\u{\theta})-1|.
|I(\u{\theta})|d\u{\theta}}\hspace{2cm}\\
&\ll&\ep^{3/2}\iint_{\max(|\theta_1|,|\theta_2|)<\ep^{-1/4}}
|I(\u{\theta})|d\u{\theta}\\
&\ll &\ep,
\end{eqnarray*}
since $I(\u{\theta})\ll 1$. Similarly we have
\begin{eqnarray*}
\iint_{\max(|\theta_1|,|\theta_2|)\ge\ep^{-1/4}}  |K_{\ep}(\u{\theta})-1|.
|I(\u{\theta})|d\u{\theta}
&\ll&
\iint_{\max(|\theta_1|,|\theta_2|)\ge\ep^{-1/4}}|I(\u{\theta})|d\u{\theta} \\
&\ll& \ep^{\frac{k-4}{8}}(1+ \log (\ep^{-1/4})),
\end{eqnarray*}
by (\ref{plane}). We therefore deduce that
\[\iint_{\R^2}|K_{\ep}(\u{\theta})-1|.|I(\u{\theta})|d\theta\rightarrow
0\]
as $\ep\rightarrow 0$,  provided that $k\ge 5$. It then follows that
\[\lim_{\ep\rightarrow 0}\iint_{\R^2}K_{\ep}(\u{\theta})I(\u{\theta})
e(-\u{\theta}\cdot\u{\mu})d\u{\theta}=\Jcal_w(\u{\mu})\]
when $k\ge 5$.

However
\[\iint_{\R^2}K_{\ep}(\u{\theta})I(\u{\theta})
e(-\u{\theta}\cdot\u{\mu})d\u{\theta}
=\int_{\R^k}w(\xbf)L(\u{Q}(\xbf)-\u{\mu})d\xbf,\]
with
\[L(\u{\lambda})=\iint_{\R^2}K_{\ep}(\u{\theta})
e(\u{\theta}\cdot\u{\lambda})d\u{\theta}.\]
We can evaluate this last integral as a product of two one-dimensional
integrals of the form
\[\int_{-\infty}^{\infty}
\left(\frac{\sin\pi\ep t}{\pi\ep t}\right)^2e(t\lambda)dt=
\left\{\begin{array}{cc} \ep^{-1}(1-\ep^{-1}|\lambda|), & |\lambda|\le\ep,\\ 0, & |\lambda|>\ep,\end{array}\right.\]
and the proposition follows.

\section{Exponential sums: the major arcs}\label{sec_exp_maj}

Let \[ S_q(\u{a}) = \sum_{\xbf \modd{q}} e_q(\u{a} \cdot \u{Q}(\xbf)),\]
\[ S_q(\u{a};\u{n}) = S_q(\u{a}) e_q(-\u{a} \cdot \u{n}),\]
and
\beq\label{T_dfn}
T(\u{n};q) = \sum_{\substack{1 \leq a_1,a_2 \leq q\\ (a_1,a_2,q)=1}}
S_q(\u{a};\u{n}).
\eeq
We will require estimates for $T(\u{n};q)$ in order to show that the
singular series converges, as well as to give a precise rate of
convergence and a lower bound for the singular series. The work on
$T(\u{n};q)$ follows  standard arguments and is less technical than
the methods we will use for exponential sums encountered in the minor
arcs, although it still requires us to make a distinction between
Type I and Type II primes with respect to each fixed value
$(n_1,n_2)$, as defined in Section \ref{sec_local_smooth}. Note that $T(\u{n};q)$ is a multiplicative function with respect to
$q$: namely if $q=q_1q_2$ with $(q_1,q_2)=1$, then \[ T(\u{n};q) = T(\u{n};q_1)T(\u{n};q_2).\]
Thus it is sufficient to consider the case where
$q=p^e$.

We first note that the following simple relation holds, for all primes $p$:
\begin{lemma}\label{TN_lemma}
For any prime $p$ and any $e\ge 1$ we have
\beq\label{T_rec}
T(\u{n};p^e) = p^{2e} N(\u{n};p^e) - p^{k+2(e-1)}N(\u{n};p^{e-1}),
\eeq
where
\[ N(\u{n};q) = \# V_{q}(n_1,n_2)= \#\{ \xbf \modd{q} : Q_1(\xbf) \con
n_1, Q_2(\xbf) \con n_2 \modd{q} \}.\] \end{lemma}
We can simply write \beq\label{Tp}
T(\u{n};p^e) = \sum_{1 \leq a_1,a_2 \leq p^{e}} S_{p^e}(\u{a};\u{n})
- \sum_{1 \leq a_1,a_2 \leq p^{e-1}} S_{p^e}(p\u{a};\u{n})
\eeq
and observe that in the second term,
\beq\label{Sp}
S_{p^e}(p\u{a};\u{n}) = \sum_{\xbf \modd{p^e}} e_{p^{e-1}}(\u{a}\cdot
\u{Q}(\xbf) - \u{a} \cdot \u{n})
	= p^kS_{p^{e-1}}(\u{a};\u{n}).
	\eeq
\label{note_cong}
Next we note that \begin{eqnarray}
\label{Np}
\sum_{1 \leq a_1,a_2 \leq p^{e}} S_{p^e}(\u{a};\u{n})
	&=& \sum_{\xbf \modd{p^e}} \sum_{1 \leq a_1,a_2 \leq p^e}
       e_{p^e}(\u{a} \cdot (\u{Q}(\xbf)-\u{n}))\nonumber \\
	 &=& p^{2e}N(\u{n};p^e).
\end{eqnarray}
Inserting (\ref{Sp}) in (\ref{Tp}) and applying the representation
(\ref{Np}) to both resulting terms proves the lemma.

Moreover, we note that the sum $\sum_e p^{-ek}T(\u{n};p^e)$ telescopes:
\begin{lemma}\label{telescope}
For any prime $p$,
\[1+\sum_{e=1}^{E} p^{-ek}T(\u{n};p^e) =p^{-E(k-2)}N(\u{n};p^E).\]
\end{lemma}
For in fact,
\begin{eqnarray*}
\sum_{e=1}^{E} p^{-ek}T(\u{n};p^e)
	&=&		\sum_{e=1}^E  p^{-e(k-2)} N(\u{n};p^e) -
       \sum_{e=1}^E p^{-(k-2)(e-1)}N(\u{n};p^{e-1})	\\
	&=& p^{-E(k-2)} N(\u{n};p^E) + \sum_{e=1}^{E-1}  p^{-e(k-2)} N(\u{n};p^e) \\
		&& \qquad  - \; \sum_{e=1}^{E-1}
               p^{-e(k-2)}N(\u{n};p^e) - N(\u{n};1) \\
	& = & p^{-E(k-2)}N(\u{n};p^E) - 1.
	\end{eqnarray*}

Thus the key to understanding $T(\u{n};p^e)$ is bounding $N(\u{n};p^e)$. We would expect that $N(\u{n};q)$ is of order $q^{k-2}$, up to a
smaller error term, since there are $q^k$ choices of $\xbf$ modulo
$q$, and the probability that a certain value $(n_1,n_2)$ is taken by
$(Q_1(\xbf), Q_2(\xbf))$ is $q^{-2}$. In the case of Type I primes,
$V_p(n_1,n_2)$ is smooth and we can apply Deligne's estimates to get a
good error term for $N(\u{n};p)$.  We can then control $N(\u{n};p^e)$
by lifting solutions modulo $p$ to solutions modulo $p^e$.
In the case of Type II primes, we obtain a slightly
worse error term for prime moduli, and we avoid $N(\u{n};p^e)$ for
prime power moduli.

\begin{prop}\label{lemma_N_bd}
For all good primes (and hence primes of Type I or Type II) we have
\beq\label{N2p}
N(\u{n};p) = p^{k-2} + O(p^{\frac{k-1}{2}}).
\eeq
For Type I primes $p$, we have
\beq\label{N1p}
N(\u{n};p) = p^{k-2} +O(p^{\frac{k-2}{2}})
\eeq
and
\beq\label{N1pe}
N(\u{n};p^e) =p^{e(k-2)} + O(p^{e(k-2)}p^{-\frac{(k-2)}{2}}).
\eeq
The implied constants may depend on $Q_1$ and $Q_2$ but are
independent of $p,e,n_1$ and $n_2$. \end{prop}

For $T(\u{n};q)$ we obtain the following bounds:
\begin{prop}\label{T_I_lemma}
For all good primes (and hence primes of Type I or Type II) we have
\beq\label{T2p}
T(\u{n};p) = O(p^{\frac{k+3}{2}}) \eeq
and
\beq\label{T2pe}
T(\u{n};p^e)  \ll  p^{e(\frac{k+4}{2})}.
\eeq
For bad primes $p$, there exists a constant $c_p$ such that for all
$e \geq 1$,  we have
\beq\label{T3pe}
| T(\u{n};p^e)|  \leq c_p p^{e(\frac{k+4}{2})}.
\eeq
For Type I primes $p$ we have
\beq\label{T1p}
T(\u{n};p) =  O(p^{\frac{k+2}{2}}) \eeq
and
\beq\label{T1pe}
T(\u{n};p^e) = 0\;\;\; (e\ge 2).
\eeq
The implied constants may depend on $Q_1$ and $Q_2$, but are independent
of $p,e,n_1$ and $n_2$. \end{prop}
Since the finite set of bad primes is determined by the original choice of
$Q_1$ and $Q_2$, and since our $\ll$ constants are allowed to depend
on this choice, we may replace (\ref{T3pe}) by the bound
$T(\u{n};p^e)\ll p^{e(\frac{k+4}{2})}$.

\subsection{Upper bounds for $N$}\label{sec_N_bounds}
We first prove Proposition \ref{lemma_N_bd}; we will prove Proposition
\ref{T_I_lemma} in Section \ref{sec_T_bounds} . We projectivize the counting problem, writing
\beq\label{N_bds}
N(\u{n};p) = (p-1)^{-1}(N^{(1)}(\u{n};p)-N^{(2)}(\u{n};p))
\eeq
with
\[N^{(1)}(\u{n};p)=\#\{ (x_0,\xbf) \in \F_p^{k+1} :
 \u{Q}(\xbf) \con \u{n}x_0^2 \modd{p} \}\]
and
\[N^{(2)}(\u{n};p)=\#\{ \xbf \in \F_p^k:
 \u{Q}(\xbf) \con \u{0} \modd{p} \},\]
say.
The last term is independent of $\u{n}$, and $p$ is a good prime
(whether it is Type I or Type II), so that the variety
defined by $\u{Q}(\xbf) =\u{0}$ is smooth over $\F_p$. Thus we may
apply Deligne's bound to obtain \beq\label{Np0}
N^{(2)}(\u{n};p)= \#\{ \xbf \in \F_p^k: \u{Q}(\xbf) \con \u{0} \modd{p}
\}  =p^{k-2} + O(p^{\frac{k-1}{2}}).
\eeq

In order to bound the first term on the right hand side of
(\ref{N_bds}), which may involve a singular variety in the case where
$p$ is of Type II, we recall the following theorem of Hooley \cite{Hoo}.
\begin{prop}\label{prop_Sal}
If $V$ is a projective complete intersection of dimension $n$ defined over the finite field $\F_p$, with singular locus of dimension $s$, then the number of $\F_p$-rational points is $(p^{n+1}-1)/(p-1)+O(p^{(n+s+1)/2})$.
\end{prop} We apply this to the variety \[ Z = \{ (x_0,\xbf) \in \F_p^{k+1} : \u{Q}(\xbf) - \u{n} x_0^2 = 0\} \]
which has projective dimension $k-2$ and singular locus of projective dimension
$s$, say, giving
\[ N^{(1)}(\u{n};p)
	= p^{k-1} + O(p^{\frac{k+1+s}{2}}).\]
\label{note_projratl}

We first verify that if $p$ is a Type I prime, then $Z$ is nonsingular
over $\F_p$. Otherwise, there is a nontrivial pair $(\xbf,x_0)$ for
which \[ \rk\left( \begin{array}{cc} 2x_0 n_1 & 	 \nabla Q_1(\xbf) \\
					2x_0 n_2 & 	\nabla Q_2(\xbf)
					\end{array} \right) <2,\]
so that certainly $\nabla Q_1(\xbf)$ and $\nabla Q_2(\xbf)$ are
 proportional over
$\F_p$. If $x_0 \neq 0$, this would imply that
$\nabla Q_1(\overline{x_0}\xbf)$ and $\nabla
Q_2(\overline{x_0}\xbf)$ are proportional,
while also $\u{Q}(\overline{x_0}\xbf) =\u{n}$
in $\F_p$.  This would contradict the fact that $p$ is Type I. If $x_0=0$,
then  $\nabla Q_1(\xbf)$ would be proportional to $\nabla Q_2(\xbf)$ while also
$\u{Q}(\xbf) =\u{0}$, which contradicts the fact that $p$ is good (and
hence satisfies Condition 1 over $\overline{\F_p}$).

Thus if $p$ is a Type I prime then $Z$ is nonsingular over
$\F_p$, so that $s=-1$ and \[ N^{(1)}(\u{n};p)
	= p^{k-1} + O(p^{\frac{k}{2}}).\]
Combining this with (\ref{Np0}) in (\ref{N_bds}) proves (\ref{N1p}).

If $p$ is a Type II prime the variety $Z$ may be singular and we must estimate the dimension of the singular locus. Let $H$ be the hyperplane $\{ (x_0,\xbf): x_0
=0\}$. Then $Z \intersect H$ is nonsingular, since by Condition 1 the
variety $V = \{ \xbf : \u{Q} (\xbf) = \u{0} \}$ is
smooth. Thus the projective dimension of the singular locus of $Z
\intersect H$ is $-1$, whence the singular locus of $Z$ itself can
have dimension at most $0$. Thus for Type II primes we have $s\leq 0$
and hence \[ N^{(1)}(\u{n};p)	= p^{k-1} + O(p^{\frac{k+1}{2}}).\]
Now combining this with (\ref{Np0}) in (\ref{N_bds}) proves (\ref{N2p}).

For prime power moduli with $p$ of Type I, we have the recursion \beq\label{N_rec}
N(\u{n};p^e) = N(\u{n};p^{e-1}) p^{k-2}.
\eeq
The proof of this follows a standard route, lifting solutions from $\F_p$ via Hensel's Lemma.  This is always possible
since $V_p(n_1,n_2)$ is smooth.  We leave the details to the reader.
Given (\ref{N_rec}) the second result (\ref{N1pe}) of
Proposition \ref{lemma_N_bd} follows immediately from (\ref{N1p}).

\subsection{Upper bounds for $T$}\label{sec_T_bounds}
We now prove Proposition \ref{T_I_lemma}.
The result (\ref{T1p}) for Type I primes follows
directly from the bound (\ref{N1p}) for $N(\u{n};p)$,
via the relation (\ref{T_rec}). Similarly, for Type II
primes, the bound (\ref{T2p}) for $T(\u{n};p)$ follows directly from
(\ref{N2p}). To prove (\ref{T1pe}) for Type I primes we again use (\ref{T_rec}), coupled now with (\ref{N_rec}).

For Type II primes, it is not effective to bound $T(\u{n};p^e)$ via
$N(\u{n};p^e)$ when $e \geq 2$, since we cannot lift
solutions. Instead, we recall that $T(\u{n};p^e)$ is a sum of twists
of \[ S_{p^e}(\u{a})  =  \sum_{\xbf \modd{p^e}} e_{p^e}(\u{a} \cdot \u{Q}(\xbf)) ,\]
and we bound $S_{p^e}(\u{a})$
directly, ignoring the role of $\u{n}$. This is also how we will
obtain the bound (\ref{T3pe}) for bad primes; this makes inherent
sense, since the property of being bad is independent of
$\u{n}$.

Thus let $p$ be any prime, good or bad. Writing $q=p^e$ for convenience, and setting $\xbf=\ybf+\zbf$, we find that
\begin{eqnarray*}
|S_q(\u{a})|^2 & = & \sum_{\xbf, \ybf \modd{q}} e_q(\u{a} \cdot
\u{Q}(\xbf) - \u{a} \cdot \u{Q}(\ybf) )\\
			& = & \sum_{\ybf, \zbf \modd{q}} e_q(\u{a}
                       \cdot \u{Q}(\zbf) +2\zbf^t (\u{a} \cdot
                       \u{Q})\ybf )\\
			& = & q^k \sum_{\substack{\zbf \modd{q} \\ q |
                           2\zbf^t (\u{a} \cdot \u{Q})}} e_q(\u{a}
                       \cdot \u{Q}(\zbf) ),
			\end{eqnarray*}
		so that
\beq\label{SZ}
|S_q(\u{a})|^2 \leq q^k Z(\u{a};q),
\eeq
where $Z(\u{a};q)$ is defined by
\beq\label{Z_dfn}
Z(\u{a};q):=\#\{\b{z}\modd{q}:\,q\mid 2\b{z}^T(\u{a} \cdot \u{Q})\}.
\eeq
The analysis of $Z(\u{a};q)$ lies more naturally within the realm of
the dichotomy of good and bad primes, which are the focus of Section
\ref{sec_min_sum}; thus for the moment, we merely state the following
claim: \begin{lemma}\label{lemma_Z_good}
For all good primes (and hence for all primes either of Type I or of Type II),
we have
\beq\label{N_bd_gb0}
Z(\u{a};p^e) \leq \gcd (F(a_1,a_2),p^e).
\eeq
For each bad prime $p$, there exists a constant $C_p$ such that for
all $e \geq 1$,  we have
\beq\label{N_bd_bg1}
Z(\u{a};p^e) \leq C_p \gcd (F(a_1,a_2),p^e).
\eeq
\end{lemma}
We will prove this result in Section \ref{sec_first_bound}; see
(\ref{N_bd_gb}). For now, we will assume the bounds (\ref{N_bd_gb0}) and
(\ref{N_bd_bg1}) and apply them to (\ref{SZ}). In order to simplify
notation, we temporarily adopt the convention that $C_p^*=1$ if $p$ is
a good prime, and $C_p^* = \sqrt{C_p}$ if $p$ is a bad prime.

We estimate $|T(\u{n};p^e)|$ as
\[
|T(\u{n};p^e)|  \leq
C_p^* p^{ek/2}
\twosum{\u{a}\modd{p^e}}{(\u{a},p)=1}\gcd(F(a_1,a_2),p^e)^{1/2}.
\]
\label{note_Faa}
The sum on the right is
\begin{eqnarray}
\lefteqn{\sum_{f=0}^e\,p^{f/2}\#\{\u{a}\modd{p^e}:\,(\u{a},p)=1,\,
\gcd(F(a_1,a_2),p^e)=p^f\}}\nonumber \\ &\le& \sum_{f=0}^e\,p^{f/2}\#\{\u{a}\modd{p^e}:\,(\u{a},p)=1,\,\,p^f\mid
F(a_1,a_2)\} \nonumber \\ &\le&  p^{2e}+
\sum_{f=1}^e\,p^{f/2} p^{2(e-f)}\#\{\u{a}\modd{p^f}:\,(\u{a},p)=1,\,\,p^f | F(a_1,a_2)\}. \label{F_sum}
\end{eqnarray}
For a given value of $f$, the number of allowable $\u{a}$ with $p\nmid
a_2$ is
\beq\label{Huxley}
\phi(p^f)\#\{u\modd{p^f}:\,p^f\mid F(u,1)\}.
\eeq
According to Huxley \cite{Hux79}  the polynomial
congruence $F(u,1)\equiv 0\modd{p^f}$ has at most $k|D_F|^{1/2}$ roots, where
$D_F$ is the discriminant of $F(x_1,x_2)$.  (It is important to note
here that Huxley's result requires $D_F\not=0$, which is certainly true in our case).  When $p\mid a_2$ we have the same estimate, by
reversing the roles of $a_1$ and $a_2$.  We therefore conclude that
\[\twosum{\u{a}\modd{p^e}}{(\u{a},p)=1}\gcd(F(a_1,a_2),p^e)^{1/2}\ll p^{2e}+
\sum_{f=1}^ep^{2e-3f/2}\phi(p^f)\ll p^{2e}.\]
\label{note_p2e}
This implies that
\[ T(\u{n};p^e) \ll p^{e\left( \frac{k+4}{2} \right) }\]
in the case of good primes, and
\[ T(\u{n};p^e) \ll C_p^* p^{e\left( \frac{k+4}{2} \right) }\]
in the case of bad primes.
This suffices for the proof of Proposition \ref{T_I_lemma}.

\subsection{The congruence problem: lower bounds for $N$}
In order to prove that the singular series is nonvanishing---and in
particular to give an effective lower bound for it---we will require a lower bound for $N(\u{n};p^e)$ with an explicit dependence
on $p$, for any prime $p$ (whether good, bad, Type I or Type II). This
effectively means we must show that the
local congruence problem $Q_1(\xbf) \con n_1$, $Q_2(\xbf) \con n_1$
modulo $p^e$ has sufficiently many solutions. Define for each $p$ the local density
\[ \sig_p(\u{n})  = \sum_{e=0}^\infty p^{-ek}T(\u{n};p^e) = \sum_{e=0}^\infty
p^{-ek}\sum_{\substack{1 \leq a_1,a_2 \leq p^e\\ (a_1,a_2,p^e)=1}}
S_{p^e}(\u{a};\u{n}) .\] We will give a lower bound for $\sig_p(\u{n})$ for appropriate $\u{n}$. \begin{prop}\label{prop_lower_bd}
For each prime $p$ there is a constant $\varpi_p>0$ such that if the system $\u{Q}(\xbf)= \u{n}$ has a solution $\xbf_0$ over $\Z_p$, then
\[\sig_p(\u{n}) \geq p^{-(k-2)}\max_{i,j}|\Delta_{ij}(\xbf_0)|_p^{2(k-2)}
\geq\varpi_p|F(n_2,-n_1)|_p^{k-2}.\]
\end{prop}
\label{note_alt44}

When $k\ge 5$, Lemma \ref{telescope} coupled with (\ref{T2pe}) and (\ref{T3pe})
shows that for all primes we have
\beq\label{sig_N}
\sig_p(\u{n}) = \lim_{e \maps \infty} \frac{N(\u{n};p^e)}{p^{e(k-2)}}.
\eeq
\label{note_sig1}

We will need to understand how
the size of $|F(n_2,-n_1)|_p$ controls the
$p$-adic valuation of $\Delta(\xbf)$.  The necessary information is
given by our next result.

\begin{lemma}\label{4rcontrol}
For each prime $p$ there is a constant $c_p>0$ as follows.  Let $\xbf_0\in\Z_p^k$ and suppose that $\u{Q}(\xbf_0)= \u{n}$.  Then
\[|F(n_2,-n_1)|_p\le c_p\max_{i,j}|\Delta_{ij}(\xbf_0)|_p^2.\]
\end{lemma}
\label{note_lem44}

For the proof we use Condition 3 to diagonalize $Q_1,Q_2$
simultaneously over $\overline{\mathbb{Q}_p}$, writing $Q_1=M^TAM$ and $Q_2=M^TBM$
with a nonsingular matrix $M$ and diagonal matrices $A={\rm
 diag}(a_i)$ and $B={\rm diag}(b_i)$.  Then since $\u{Q}(\xbf_0) = \u{n}$,
\[F(n_2,-n_1)=\det\left(Q_2(\xbf_0)Q_1-Q_1(\xbf_0)Q_2\right)=
\det(M)^2\det\left(B(\ybf)A-A(\ybf)B\right),\]
where $\ybf=M\xbf_0$.  If we set $\Delta_{ij}=a_ib_j-a_jb_i$ then the matrix
$\Delta(\ybf;A,B)$, as defined in (\ref{G_dfn}), has entries
$\Delta_{ij}(\ybf)=4\Delta_{ij}y_iy_j$; note that for $i \neq j$,
$\Del_{ij}$ is nonzero by Condition 3.  Moreover \[\det\left(B(\ybf)A-A(\ybf)B\right)=\prod_{i=1}^k\left(\sum_{\substack{1\le
       j\le k\\ j\not=i}}\Delta_{ij}y_j^2\right),\]
which is a form of degree $2k$ in $\ybf$. On expanding out the product we see that every resulting term contains
a factor $y_i^2y_j^2$ for some pair $i \neq j$, and this factor is a
constant multiple of $\Del_{ij}(\ybf)^2$.  It is therefore possible
to choose certain forms $G_{ij}$, depending on $A$
and $B$, for which there is an identity in $\ybf$ of the shape
\[\det\left(B(\ybf)A-A(\ybf)B\right)=\sum_{i,j}\Delta_{ij}(\ybf)^2G_{ij}(\ybf).\]
Since by assumption $|\xbf_0|_p \leq 1$, we have $|\ybf|_p\ll_p 1$, and hence
$|G_{ij}(\ybf)|_p \ll_p 1$. We then deduce that
\[|F(n_2,-n_1)|_p\ll_p\max_{i,j}|\Delta_{ij}(\ybf)^2|_p.\]
To complete the proof of the lemma it remains to observe that
\[\frac{\partial A(\ybf)}{\partial y_i}=
\sum_j(M^{-1})_{ji}\frac{\partial Q_1(\xbf)}{\partial x_j}\]
and similarly for $B$ and $Q_2$. This may be verified by a tedious
piece of linear algebra. \label{note_tedious}
It follows that the determinants
$\Delta_{ij}(\ybf)$ are linear combinations of the entries
$\Delta_{ij}(\xbf_0)$ in $\Del(\xbf_0; Q_1,Q_2)$, so that
\[\max_{i,j}|\Delta_{ij}(\ybf)^2|_p\ll_p\max_{i,j}|\Delta_{ij}(\xbf_0)^2|_p.\]
This establishes the lemma.

We can now complete the proof of Proposition \ref{prop_lower_bd}. After permuting the indices as necessary, Lemma \ref{4rcontrol} allows us to assume that that \beq\label{y4}
|\Delta_{12}(\xbf_0)|_p^2=\max_{i,j}|\Delta_{ij}(\xbf_0)|_p^2\ge c_p^{-1}|F(n_2,-n_1)|_p.
\eeq
Since Proposition \ref{prop_lower_bd} is trivial unless
$\Delta_{12}(\xbf_0)$ is non-zero we may set
\beq\label{y5}
|\Delta_{12}(\xbf_0)|_p=p^{-v},
\eeq
say. Suppose now that
$e\ge 2v+1$ and take any integers $y_3,\ldots,y_k\modd{p^e}$ all
divisible by $p^{2v+1}$.  There are $p^{(k-2)(e-2v-1)}$ such $(k-2)$-tuples.
We then claim that the congruences
\[\u{Q}(\xbf_0+(w_1,w_2,y_3,y_4,\ldots,y_k))\equiv\u{n}\modd{p^e}\]
will have an integer solution $(w_1,w_2)$.  Thus
$N(\u{n};p^e)\ge p^{(k-2)(e-2v-1)}$, whence (\ref{sig_N}) yields
\[\sigma_p(\u{n})\ge p^{-(k-2)(2v+1)} = p^{-(k-2)}
|\Del_{12}(\xbf_0)|_p ^{2(k-2)} \geq (c_pp)^{-(k-2)}
|F(n_2,-n_1)|_p^{k-2}\] in view of (\ref{y4}) and (\ref{y5}). This will suffice for
Proposition \ref{prop_lower_bd}, with $\varpi_p = (c_pp)^{-(k-2)}.$

To prove the claim we  set
\[\ybf=\xbf_0+(0,0,y_3,y_4,\ldots,y_k).\]
Then
$\ybf\equiv\xbf_0\modd{p^{2v+1}}$, whence
 $\u{Q}(\ybf)\equiv\u{n}\modd{p^{2v+1}}$ and \[\Delta_{12}(\ybf)\equiv\Delta_{12}(\xbf_0)\modd{p^{2v+1}}.\]
In particular we will have $|\Delta_{12}(\ybf)|_p=p^{-v}$, in view of (\ref{y5}).
\label{note_Del12}
We therefore
require an integer solution $(w_1,w_2)$ to the simultaneous congruences
\[q_1(w_1,w_2)\equiv q_2(w_1,w_2)\equiv 0 \modd{p^e},\]
where
\[q_i(w_1,w_2)=Q_i(\xbf_0+(w_1,w_2,y_3,y_4,\ldots,y_k)) - Q_i(\xbf_0),\;\;\;(i=1,2).\]
The existence
of suitable $w_1,w_2$ now follows from a standard application of
Hensel's Lemma, since $p^{2v+1}|q_1(0,0), q_2(0,0)$ and
\[\det\left(\begin{array}{cc}\frac{\partial q_1(0,0)}{\partial w_1} &
\frac{\partial q_1(0,0)}{\partial w_2}\\
\frac{\partial q_2(0,0)}{\partial w_1}&
\frac{\partial q_2(0,0)}{\partial w_2}
\end{array}\right)\not\equiv 0\modd{p^{v+1}}.\]
\label{note_q12}
This completes the proof of Proposition \ref{prop_lower_bd}.

Theorem \ref{thm_s2} relates to the case $(n_1,n_2) = (0,0)$ for
 which we replace Proposition \ref{prop_lower_bd} with the following
 result.
  \begin{prop}\label{prop_lower_bd0}
For each prime $p$ such that the system $\u{Q}(\xbf)= \u{0}$ has a nonzero local solution in $\Z_p$ we have
$\sig_p(\u{0})>0$.
\end{prop}
This will follow immediately from the following analogue of Lemma \ref{4rcontrol}.
\begin{lemma}\label{4rcontrol0}
For each prime $p$ such that the system $\u{Q}(\xbf)= \u{0}$ has a nonzero local solution $\xbf_0$ over $\Z_p$ we have
\[\max_{i,j}|\Del_{i,j}(\xbf_0)|_p>0.\]
\end{lemma}
For the proof we merely observe that if $\Del_{i,j}(\xbf_0)=0$ for
every pair $i,j$ then $\rk (J(\xbf_0))<2$, which contradicts Condition \ref{cond1} over $\overline{\Q_p}$, and hence equivalently
over $\overline{\Q}$.

\section{Exponential sums: minor arcs}\label{sec_min_sum}
In this section we consider exponential sums of the form
\[ S(\lbf;q) = \sum_{\substack{a_1,a_2 \\ (a_1,a_2,q)=1}}
\sum_{\substack{\rbf_1 \modd{q}\\ \rbf_2 \modd{q}}} e_q(\u{a} \cdot
\u{Q}(\rbf_1)- \u{a} \cdot \u{Q}(\rbf_2))e_q(\rbf_1\cdot  \lbf_1 +
\rbf_2 \cdot \lbf_2)\] where $\lbf = (\lbf_1, \lbf_2) \in \Z^k \times \Z^k$. Such sums arise
naturally in the analysis of the minor arcs, for which the original
generating function is squared.

We first observe that $S(\lbf;q)$ satisfies a standard multiplicative property: if $q=q_1q_2$ with $(q_1,q_2)=1$, then one may easily check that
\begin{equation}\label{rmult} S(\lbf;q) = S(\lbf;q_1)S(\lbf;q_2).
\end{equation}
We may thus reduce our consideration to the case in which $q=p^e$ for some prime $p$ and $e\geq 1$.

The sum $S(\lbf;q)$ is in $2k+2$ variables, so a square-root
cancellation bound would take the form $S(\lbf;q) \ll q^{k+1}$. This
is out of reach; for general primes $p$, the best we can prove is that
$S(\lbf; p^e) \ll_{p,e} p^{e(k+2)}$. However, for most primes we can
do a bit better, and this is the heart of the work in showing that the
minor arcs make a small enough contribution.

\subsection{Preliminaries}\label{rss1}

Recall the determinant form
\[F(x,y):=\det(xQ_1+yQ_2)\]
which we factored over its splitting field $K$ as
\beq\label{Fxy}
F(x,y)=c^{-1}\prod_{i=1}^k(\lambda_i x-\mu_i y)
\eeq
with $c\in\N$ and $\lambda_i,\mu_i\in\mathcal{O}_K$. According to Condition 4,
for each index $i$ the matrix $(\mu_i,\lambda_i).\u{Q}$ will have rank exactly $k-1$. \label{note_k1}
In particular the null space of $(\mu_i,\lambda_i).\u{Q}$ is one-dimensional,
and we may choose a generator $\b{e}_i\in\mathcal{O}_K^k$ for it. We
then define
\begin{equation}\label{rGprod}
G(\b{x}):=\prod_{i=1}^k\prod_{\sigma}\left(\b{x}^T\b{e}_i^{\sigma}\right),
\end{equation}
where $\sigma$ runs over ${\rm Gal}(K/\Q)$.  This
will be a non-zero form in $k$ variables of a certain fixed degree, with coefficients in
$\Z$.

A further collection of forms will also appear in our analysis. To define these we consider the $(k+1) \times (k+1)$ matrix
\begin{equation}\label{rmat}
M(\u{\alpha};\b{x},\b{y}):=
\left(\begin{array}{c|c} \rule[-3mm]{0mm}{1mm}\u{\alpha} \cdot \u{Q} & \b{y}\\ \hline
\rule[3mm]{0mm}{1mm}\b{x}^T & 0\end{array}\right),
\end{equation}
where $\b{x}$ and $\b{y}$ are vectors of length $k$, and write
$H_{rs}(\u{\alpha};\b{x},\b{y})$ for its $(r,s)$-minor, this being a
determinant of size $k\times k$.  We then define
\begin{equation}\label{rHprod}
H_{rsi}(\b{x},\b{y}):=\prod_{\sigma}
H_{rs}\big((\mu_i^{\sigma},\lambda_i^{\sigma});\b{x},\b{y}\big)
\end{equation}for $r,s\le k+1$ and $i\le k$, where as usual $\sig$ runs over $\text{Gal}(K/\Q)$. Finally we set
\[H(\b{x},\b{y}):=\prod_{i=1}^k
\left\{\sum_{r,s\le k+1}H_{rsi}(\b{x},\b{y})^2\right\}.\]
Note that the coefficients of $H_{rsi}$ and $H$ are rational integers.

\subsection{Results}

Recall the definition of ``bad'' primes as given in (\ref{rbad}). The
bound we obtain for the exponential sum $S(\lbf;p^e)$ will depend on whether $p$ is good or bad.

\begin{lemma}\label{lemma_bad_p}
For each bad prime $p$, there exists a constant $c_p$ such that
\[ |S(\lbf;p^e)| \leq c_p (e+1)p^{e(k+2)}.\]
\end{lemma}
If there were square-root cancellation we would have a bound
$O(p^{e(k+1)})$, so our bound is worse by a factor of $p^e$. Moreover we have an implied constant which depends on $p$.
However this is acceptable since there are only finitely many bad primes $p$.

For good primes, we will prove: \begin{lemma}\label{prop_min_sum_bd0}
For good primes $p$ we have
\[ S(\lbf; p^e) \ll (e+1)p^{e(k+2)} \]
in all cases.  Moreover if we set $\lbf_3:=\lbf_1+\lbf_2$ and $\lbf_4:=\lbf_1-\lbf_2$ we will have
\begin{equation}\label{rpropest}
S(\lbf; p^e) \ll (e+1)p^{ek}\left(p^{2e-1}+
\#\{\u{b}\modd{p^e}:\, p\mid\det M(\u{b};\lbf_3,\lbf_4)\}\right)
\end{equation}
unless $p$ divides each of
\[G(\lbf_3),\;\; G(\lbf_4)\;\;\mbox{and}\;\;
H(\lbf_3,\lbf_4).\]
\end{lemma}

This is a truly unpleasant result, and its proof and later application
are the most awkward parts of our entire argument.

\subsection{A first bound}\label{sec_first_bound}

In this section we shall prove the following estimate.
\begin{lemma}\label{r1}
For good primes $p$ we have
\[|S(\lbf;p^e)|\ll (e+1)p^{e(k+2)},\]
while for bad primes we have
\[|S(\lbf;p^e)|\ll_p (e+1)p^{e(k+2)}\]
uniformly in $e$.
\end{lemma}
Evidently Lemma \ref{r1} implies
Lemma \ref{lemma_bad_p} and the first part of Lemma \ref{prop_min_sum_bd0}.

For the proof we write $q=p^e$ for convenience, and recall that
\begin{equation}\label{rSdef}
S(\lbf;q) = \twosum{\u{a}\modd{q}}{(\u{a},q)=1}\;
\twosum{\b{x}\modd{q}}{\b{y}\modd{q}}
e_q (\u{a}\cdot \u{Q}(\b{x})-\u{a}\cdot
\u{Q}(\b{y})+\lbf_1^T\b{x}+\lbf_2^T\b{y}). \end{equation}
Making the change of variables $\b{x}=\b{z}+\b{y}$ we now obtain \begin{eqnarray*}
S(\lbf;q)& = &\twosum{\u{a}\modd{q}}{(\u{a},q)=1}\;
\sum_{\b{y},\b{z}\modd{q}}
e_q(\u{a}\cdot \u{Q}(\b{z})+2\b{z}^T(\u{a}\cdot \u{Q})\b{y}+\lbf_1^T\b{z}
+\lbf_3^T\b{y})\\
& = &q^k\twosum{\u{a}\modd{q}}{(\u{a},q)=1}\;
\twosum{\b{z}\modd{q}}{q\mid 2\b{z}^T(\u{a}\cdot \u{Q})+\lbf_3^T}
e_{q}(\u{a}\cdot \u{Q}(\b{z})+\lbf_1^T\b{z}).
\end{eqnarray*}
Hence \[ |S(\lbf;q)| \leq q^k\twosum{\u{a}\modd{q}}{(\u{a},q)=1}
\# \Scal(\lbf_3,\u{a};q),\]
where \begin{equation}\label{N_a}
\Scal(\lbf_3,\u{a};q)=\{\b{z}\modd{q}:\,q\mid 2\b{z}^T(\u{a}\cdot
\u{Q})+\lbf_3^T\}.
\end{equation}
Note that if a set $\Scal(\lbf_3,\u{a};q)$ is non-empty then it must
be a coset  of $\Scal(\b{0},\u{a};q)$ in $(\Z/q\Z)^k$. For indeed, if
$\Scal(\lbf_3,\u{a};q)$ is non-empty then there exists a solution $\zbf_1$ of
$2\zbf_1^T (\u{a} \cdot \u{Q})+\lbf_3^T\con\mathbf{0} \modd{q}$, so that every
solution $\zbf$ lying in $\Scal(\lbf_3,\u{a};q)$ can be written as
$\zbf = \zbf_1  + \zbf_0$, where $\zbf_0 \in \Scal(\mathbf{0},
\u{a};q)$.

As a consequence,
\[ |S(\lbf;q)| \leq
q^k\twosum{\u{a}\modd{q}}{(\u{a},q)=1}Z(\u{a};q)\]
where as in (\ref{Z_dfn}) we have defined
\begin{equation}\label{rN}
Z(\u{a};q):=\#\{\b{z}\modd{q}:\,q\mid 2\b{z}^T(\u{a}\cdot \u{Q})\}.
\end{equation}

Our task is now to analyse $Z(\u{a};q)$. We will first prove Lemma
\ref{lemma_Z_good}, which gives an upper bound for $Z(\u{a};q)$, and then use it to prove the
following result, from which Lemma \ref{r1}
is now an immediate consequence.
\begin{lemma}\label{r2}
For good primes $p$ we have
\[\twosum{\u{a}\modd{p^e}}{(\u{a},p)=1}Z(\u{a};p^e)\ll (e+1)p^{2e}\]
while for bad primes we have
\[\twosum{\u{a}\modd{p^e}}{(\u{a},p)=1}Z(\u{a};p^e)\ll_p (e+1)p^{2e}.\]
\end{lemma}

We proceed to prove Lemma \ref{lemma_Z_good}.  We begin by putting the matrix $\u{a}\cdot \u{Q}$ into Smith normal form as \begin{equation}\label{snf}
 S(\u{a}\cdot \u{Q})T = \diag(d_i),
\end{equation}
say, using unimodular integer matrices $S,T$.
\label{note_Smith}
We may assume here that the integer diagonal entries $d_i$ are ordered so that their $p$-adic valuations satisfy
\begin{equation}\label{rord}
|d_1|_p \leq |d_2|_p \leq \ldots\leq |d_k|_p\le 1.
\end{equation}
If $p\mid d_i$ for some
$i\ge 2$, then
$S(\u{a}\cdot \u{Q})T$ has rank at most $k-2$ over $\F_p$. Since $S,T$
have determinants $\pm 1$ over $\Z$, they also have determinants $\pm
1$ over $\F_p$, and are hence invertible over $\F_p$. Thus if $p \mid
d_i$ for some $i \geq 2$, it would then follow that $\u{a}\cdot \u{Q}$ also has rank at most $k-2$ over
$\F_p$.  Thus, by the local version of Condition \ref{cond2} given
in Lemma \ref{lemma_Delta_p}, we can only have $p|d_i$ with $i\ge 2$ in the case in which $p$ is bad.

\label{note_Z}
We now see that
\[Z(\u{a};p^e) = \#\{\b{w}\modd{p^e}:\,p^e\mid 2\b{w}^T\diag(d_i)\}.\]
For each $i$, the congruence $2w_id_i\equiv 0\modd{p^e}$
has exactly $(2d_i,p^e)$ solutions modulo $p^e$, so that \[Z(\u{a};p^e)=\prod_{i=1}^k(2d_i,p^e).\]
We now call on the following result, which we shall prove in a moment.
\begin{lemma}\label{bad_p_unif_v}
For each bad $p$ there is a constant $c_p$ such that
\[ \nu_p(2d_2) \leq c_p \]
for every pair $\u{a}$ with $(a_1,a_2,p)=1$.
\end{lemma}

We have already observed that $(d_i,p)=1$ for $i\ge 2$ and good primes
$p$, and we now conclude that
\[Z(\u{a};p^
e)\le\left\{\begin{array}{cc}
(d_1,p^e),& p\,\mbox{ good},\\
\rule{0cm}{5mm}p^{(k-1)c_p}(2d_1,p^e),& p\,\mbox{ bad}.\end{array}\right.\]
In order to relate $d_1$ to the vector $\u{a}$ we note from
(\ref{snf}) that
\[(d_1,p^e)\mid\det\left(S(\u{a}\cdot \u{Q})T\right), \]
and since $S$ and $T$ are unimodular we see that \[(d_1,p^e)\mid \left(\det(\u{a}\cdot \u{Q}), p^e\right)=(F(a_1,a_2),p^e).\] It then follows that
\beq\label{N_bd_gb}
Z(\u{a};p^e)\le\left\{\begin{array}{cc}
(F(a_1,a_2),p^e),& p\,\mbox{ good},\\
\rule{0cm}{5mm}2p^{(k-1)c_p}(F(a_1,a_2),p^e),& p\,\mbox{ bad}.\end{array}\right.
\eeq
This establishes Lemma \ref{lemma_Z_good}, subject to the proof of
Lemma \ref{bad_p_unif_v}.
\label{note_Z2}

\label{note_Z3}
To handle Lemma \ref{r2} we now sum over $\u{a} \modd{p^e}$ with
$(\u{a},p^e)=1$, via a similar procedure as in (\ref{F_sum}), with the
only modification that a factor $(e+1)$ appears. This yields: \begin{eqnarray*}
\lefteqn{\twosum{\u{a}\modd{p^e}}{(\u{a},p)=1}(F(a_1,a_2),p^e)}
\hspace{2cm}\\
&\ll&
\sum_{f=0}^e\,p^f\#\{\u{a}\modd{p^e}:\,(\u{a},p)=1,\,(F(a_1,a_2),p^e)=p^f\}\\
&\ll& (e+1)p^{2e},
\end{eqnarray*} and Lemma \ref{r2} follows.

\label{note_compact}
We conclude the present subsection by establishing Lemma
\ref{bad_p_unif_v}.  We remark that the lemma may be seen as a
$p$-adic analogue of Lemma \ref{lemma_rho}, and our proof will follow
similar lines.  We argue by contradiction, using a compactness
argument over the $p$-adics $\Z_p$. Our assumption, contrary to the conclusion of Lemma \ref{bad_p_unif_v}, is that there exists a sequence of unimodular matrices $S_i, T_i$, and of pairs $\u{a}^{(i)}$,  all of whose
entries are in $\Z_p$, and such that the values
$|d_2^{(i)}|_p$ of the corresponding Smith normal forms tend to
zero. By compactness of $\Z_p$, there will be a subsequence of triples $(S_i,T_i,\u{a}^{(i)})$ converging to
$S^{(*)},T^{(*)},\u{a}^{(*)}$, say. It follows that $S^{(*)}$ and
$T^{(*)}$ are unimodular matrices, and also that at least one of $a_1^{(*)}$ and $a_2^{(*)}$ is a $p$-adic unit.  We then see from (\ref{snf}) that
\[  S^{(*)}(\u{a}^{(*)}\cdot \u{Q})T^{(*)} = \diag(d_i^{(*)})\]
with $d_2^{(*)}=0$.  However (\ref{rord}) shows that if
$|d_2^{(i)}|_p$ tends to zero then so also must $|d_1^{(i)}|_p$. We
deduce that $d_1^{(*)}=d_2^{(*)}=0$, whence $\diag(d_i^{(*)})$ has
rank at most $k-2$. Since the matrices $S^{(*)}$ and $T^{(*)}$ are invertible over $\Z_p$ we
see that $\u{a}^{(*)}\cdot \u{Q}$ also has rank at most $k-2$ over
$\Z_p$.  This finally contradicts Condition 4 with respect to $\overline{\Q_p}$ since the
pair $\u{a}^{(*)}$ cannot vanish in $\Z_p$. (Recall that Condition 4
over $\overline{\Q_p}$ is implied by Condition 1 over
$\overline{\Q_p}$, which in turn is equivalent to Condition 1 over
$\Qbar$.) This completes the proof of Lemma \ref{bad_p_unif_v}.

\subsection{Good primes}

\label{note_hom}
We now turn to the good primes $p$, for which it remains to prove the
second part of Lemma \ref{prop_min_sum_bd0}. We begin with a
slight variant of the previous analysis in which we make a symmetric change of variables, writing $\b{x}=\b{u}+\b{v}$ and
$\b{y}=\b{u}-\b{v}$.  Since $p=2$ is a bad prime, by convention, this
is permissible.  We now find for $q=p^e$ that
\begin{eqnarray*}
S(\lbf;q) & = & \twosum{\u{a}\modd{q}}{(\u{a},q)=1}\;
\sum_{\b{u},\b{v}\modd{q}}
e_q(4\b{u}^T(\u{a}\cdot \u{Q})\b{v}+\lbf_3^T\b{u}+\lbf_4^T\b{v})\\
& = &q^k\twosum{\u{a}\modd{q}}{(\u{a},q)=1}\;
\twosum{\b{u}\modd{q}}{q|4\b{u}^T(\u{a}\cdot \u{Q})+\lbf_4^T}
e_{q}(\lbf_3^T\b{u}).
\end{eqnarray*}
Since $p$ is odd, the factor 4 can be absorbed into $\u{a}$.  We may
simplify the exponential sum further by using homogeneity, averaging over an auxiliary variable $r$ as follows.
We have
\begin{eqnarray*}
\phi(q) S(\lbf;q) & = &q^k\sum_{(r,q)=1}\;\twosum{\u{a}\modd{q}}{(\u{a},q)=1}
\;\twosum{\b{u}\modd{q}}{q|r^{-1}\b{u}^T(r\u{a}\cdot \u{Q})+\lbf_4^T}
e_{q}(\lbf_3^T\b{u})\\
& = &q^k\sum_{(r,q)=1}\;\twosum{\u{b}\modd{q}}{(\u{b},q)=1}\;
\twosum{\b{w}\modd{q}}{q|\b{w}^T(\u{b}\cdot \u{Q})+\lbf_4^T}
e_{q}(r\lbf_3^T\b{w})
\end{eqnarray*}
on replacing $r\u{a}$ by $\u{b}$ and $r^{-1}\b{u}$ by $\b{w}$.  We can
now perform the summation over $r$, on recalling that $q=p^e$.  This
produces
\beq\label{SN1N2}
\phi(p^e) S(\lbf;p^e)=p^{ek}\left(-p^{e-1}N_1(p^e)+p^eN_2(p^e)\right),
\eeq
where
\[N_1(p^e):=\#\{\u{b},\b{w}\modd{p^e}:\,(\u{b},p)=1,\,
p^e\mid\b{w}^T(\u{b}\cdot \u{Q})+\lbf_4^T,\, p^{e-1}\mid \lbf_3^T\b{w}\},\]
and
\[N_2(p^e):=\#\{\u{b},\b{w}\modd{p^e}:\,(\u{b},p)=1,\,
p^e\mid\b{w}^T(\u{b}\cdot \u{Q})+\lbf_4^T,\, p^e\mid \lbf_3^T\b{w}\}.\]

\label{note_N1}
We may easily dispose of $N_1(p^e)$ by ignoring the condition $p^{e-1}\mid
\lbf_3^T\b{w}$. The set
\[\{\u{b},\b{w}\modd{p^e}:\,(\u{b},p)=1,\,p^e\mid\b{w}^T(\u{b}\cdot
\u{Q})+\lbf_4^T\}\] is either empty, or is a coset (in $\Z/p^e\Z$) of
\beq\label{bz_set}
\{\u{b},\b{w}\modd{p^e}:\,(\u{b},p)=1,\,p^e\mid\b{w}^T(\u{b}\cdot \u{Q})\},
\eeq
by the same reasoning applied to (\ref{N_a}). Certainly the set
(\ref{bz_set}) is contained in \[ \{\u{b},\b{w}\modd{p^e}:\,(\u{b},p)=1,\,
p^e\mid 2\b{w}^T(\u{b}\cdot \u{Q})\},\]
so that in the notation of (\ref{rN}), \[N_1(p^e)\le\twosum{\u{b}\modd{p^e}}{(\u{b},p)=1}\;Z(\u{b};p^e).\]
Then Lemma \ref{r2} yields
\[N_1(p^e)\ll (e+1)p^{2e}\]
for good primes $p$. It follows that the corresponding contribution to $S(\lbf;p^e)$, after dividing (\ref{SN1N2}) through by $\phi(p^e)$,
is therefore $O((e+1)p^{e(k+2)-1})$, which is satisfactory for (\ref{rpropest}).

We turn now to $N_2(p^e)$, which will be somewhat harder. We decompose $N_2(p^e)$ into two parts
$N_3(p^e)+N_4(p^e)$, the first of which corresponds to pairs $\u{b}$ for which $\u{b}\cdot \u{Q}$ is invertible modulo $p$ and the second to those
for which it is not.  For $N_3(p^e)$ the condition $p^e\mid
\b{w}^T(\u{b}\cdot \u{Q})+\lbf_4^T$ determines $\b{w}\modd{p^e}$
uniquely. \label{note_w_unq} Moreover, for a fixed $\u{b} \modd{p^e}$, the system of linear equations
\begin{equation}\label{rsys}
\b{w}^T(\u{b}\cdot \u{Q})+\lbf_4^T=\b{0}^T,\; \lbf_3^T\b{w}=0
\end{equation}
will have a solution over $\F_p$ if and only if the matrix
$M(\u{b};\lbf_3,\lbf_4)$ given by (\ref{rmat})
has rank $k$ over $\F_p$.  It follows that
\[N_3(p^e)\le\#\{\u{b}\modd{p^e}:\, p\mid\det M(\u{b};\lbf_3,\lbf_4)\},\]
which leads to a satisfactory contribution in (\ref{rpropest}).

\label{note_N4}
We handle $N_4(p^e)$ by splitting it as $N_5(p^e)+N_6(p^e)$, the first
of which counts those solutions for which the matrix
\[M_0(\u{b};\lbf_3):=\left(\begin{array}{c} \rule[-3mm]{0mm}{1mm}\u{b}\cdot \u{Q} \\ \hline  \rule[3mm]{0mm}{1mm}\lbf_3^T \end{array}\right)\]
has rank $k$ over $\F_p$, and the second of which counts solutions for
which the rank is $k-1$. Here we recall that the rank of
$M_0(\u{b};\lbf_3)$ cannot be lower than $k-1$ over $\F_p$, since $p$ is a good prime.

\label{note_N5}
In the case of $N_5(p^e)$, we have
\[\b{w}^T(\u{b}\cdot \u{Q})\equiv\b{0}^T,\,\lbf_3^T\b{w}\equiv 0\modd{p}\]
only when $p\mid\b{w}$.  Then the only solution modulo $p^e$
will be $\b{w}\equiv\b{0}\modd{p^e}$, so that (again by the coset
argument) the congruences \[\b{w}^T(\u{b}\cdot
\u{Q})+\lbf_4^T\equiv\b{0}^T,\,\lbf_3^T\b{w}\equiv 0\modd{p^e}\] have at most one solution $\wbf \modd{p^e}$. Moreover, as in our
treatment of $N_3(p^e)$, in order for (\ref{rsys}) to have a solution over $\F_p$ the matrix $M(\u{b};\lbf_3,\lbf_4)$ must have rank $k$ over $\F_p$ and hence
\[N_5(p^e)\le\#\{\u{b}\modd{p^e}:\, p\mid\det M(\u{b};\lbf_3,\lbf_4)\},\]
which is again satisfactory.

\label{note_N6}
To handle $N_6(p^e)$ we must first understand better the condition that $\u{b}\cdot \u{Q}$ is singular modulo
$p$.  We recall the notation defined in \S \ref{rss1} and we take $P$
to be a prime ideal of $\mathcal{O}_K$ lying over $p$.  Since
$p$ is a good prime neither $p$ nor $P$ can divide the constant $c$ in
(\ref{Fxy}).  Moreover, since $p$ is good, none of the determinants $\lambda_i\mu_j-\lambda_j\mu_i$ can be divisible
by $P$ and so in particular, for each $i$, the ideal $P$ can divide at
 most one of $\lam_i$ or $\mu_i$.  For
 $N_6(p^e)$ we have $p \mid \det (\u{b} \cdot \u{Q})$ and hence
 $P\mid\det(\u{b}\cdot \u{Q})$.  It follows that
\[P\mid\prod_{i=1}^k(\lambda_i b_1-\mu_i b_2),\]
so that $P\mid \lambda_i b_1-\mu_i b_2$ for some index $i$. Since $(b_1,b_2,p)=1$ and $P$ cannot divide both $\lam_i$ and
 $\mu_i$, we deduce that
$\u{b}\equiv\alpha(\mu_i,\lambda_i)\modd{P}$ for some $\alpha\in\mathcal{O}_K$ not divisible by $P$. As in our treatment of $N_3(p^e)$ we obtain a codition for
(\ref{rsys}) to be solvable modulo $p$.  Here we require the matrices
$M_0(\u{b};\lbf_3)$ and $M(\u{b};\lbf_3,\lbf_4)$ to have the same rank over $\F_p$.  Since the first of these has rank $k-1$ for
$N_6(p^e)$ we see that all the $k\times k$ minors of
$M(\u{b};\lbf_3,\lbf_4)$ must vanish modulo $p$.  It follows that
$p\mid H_{rs}(\u{b};\lbf_3,\lbf_4)$ for each $r,s\le k+1$, and hence
that $P\mid H_{rs}(\alpha(\mu_i,\lambda_i);\lbf_3,\lbf_4)$.  However the
polynomial $H_{rs}(\u{x};\b{y},\b{z})$ is homogeneous in $\u{x}$, and we may
therefore conclude that there is an index $i$ for which each of
$H_{rs}\left((\mu_i,\lambda_i);\lbf_3,\lbf_4\right)$, for $r,s\le
k+1$, is divisible by $P$.  We then deduce that $H_{rsi}(\lbf_3,\lbf_4)$ is
divisible by $P$ for each $r,s \leq k+1$, since all the other factors in (\ref{rHprod}) are
algebraic integers.  This then implies firstly that $P$ divides
$H(\lbf_3,\lbf_4)$ and then that $p\mid
H(\lbf_3,\lbf_4)$, because the coefficients of $H$ are rational
integers. It follows that $N_6(p^e)$ vanishes unless $p\mid H(\lbf_3,\lbf_4)$.

Finally, recall the requirement for $N_6(p^e)$: if $\u{b}$ and $\b{w}$ are counted by $N_6(p^e)$, then $M_0(\u{b};\lbf_3)$ has rank $k-1$ over $\F_p$.  Thus
$\lbf_3=(\u{b}\cdot \u{Q})\b{x}$ for some $\b{x}\in\F_p^k$ and then, with some
abuse of notation,
\[(\u{b}\cdot \u{Q})\b{x}\equiv \lbf_3\modd{P}.\]
By the argument above one then has
\[\left(\alpha(\mu_i,\lambda_i)\cdot \u{Q}\right)\b{x}\equiv\lbf_3\modd{P}\]
for some some index $i$ and $\alpha\in\mathcal{O}_K$ such that $P\ndiv \al$. Pre-multiplying by the null-vector $\b{e}_i$ described in \S \ref{rss1} we deduce that
$P\mid \b{e}_i^T\lbf_3$.  It then follows that $G(\lbf_3)$ is
divisible by $P$, since all the other factors in (\ref{rGprod}) are
algebraic integers. Just as in the case of $H(\lbf_3,\lbf_4)$ we now
deduce that $N_6(p^e)$ vanishes unless $p\mid G(\lbf_3)$.
In exactly the same way, if $p$
divides $\b{w}^T(\u{b}\cdot \u{Q})+\lbf_4^T$ we may conclude that
$p\mid G(\lbf_4)$. Thus $N_6(p^e)$ vanishes unless $p$ divides each of
$H(\lbf_3,\lbf_4)$, $G(\lbf_3)$ and $G(\lbf_4)$.  This now
suffices for Lemma \ref{prop_min_sum_bd0}.

\subsection{Averages of $S(\lbf;q)$}

We will now use Lemmas \ref{lemma_bad_p} and \ref{prop_min_sum_bd0} to
establish the following estimate.
\begin{lemma}\label{ravlem}
Let $\ep>0$ be given and suppose that $k\ge 4$.  Then for any $Q,L\in \N$ we have
\[\sum_{q\le Q}\;\sum_{|\lbf|\le L}|S(\lbf;q)|\ll_{\ep}
Q^{k+2+\ep}L^{2k+\ep}+Q^{k+3+\ep}L^{2k-3+\ep}.\] \end{lemma}

We shall write $\Sigma$ for the double sum to be estimated, and we
begin by considering those terms for which $G(\lbf_3)$, $G(\lbf_4)$
and $H(\lbf_3,\lbf_4)$ all vanish.  We write $\Sigma_1$ for the
corresponding contribution to $\Sigma$.  In this case Lemmas
\ref{lemma_bad_p} and \ref{prop_min_sum_bd0} show that there is a
constant $C$ such that $|S(\lbf;p^e)|\le C(e+1)p^{e(k+2)}$ for all
prime powers $p^e$.  One merely takes $C$ to be the maximum of the
implied constant from Lemma \ref{prop_min_sum_bd0} and the various constants
$c_p$ in Lemma \ref{lemma_bad_p} for the finite number of bad primes.
In view of the multiplicative property (\ref{rmult}) we then see that
\[|S(\lbf;q)|\le C^{\omega(q)}d(q)q^{k+2}\ll_{\ep}q^{k+2+\ep}, \]
where $\omega(q)$ is the number of distinct prime factors of $q$, and
$d(q)$ is the divisor function. In particular, we record for later reference that in the case $\lbf={\bf 0}$,
\beq\label{S0}
|S({\mathbf{0}};q)|\ll_{\ep}q^{k+2+\ep}.
\eeq
Continuing our estimation, it follows that
\[\Sigma_1\ll_{\ep}Q^{k+3+\ep}\#\{\lbf:\,|\lbf|\le L,\,
G(\lbf_3)=G(\lbf_4)=H(\lbf_3,\lbf_4)=0\}.\]
Examining the definition of $H(\b{x},\b{y})$ we see that there must be
some index $i\le k$ such that
\[\Sigma_1\ll_{\ep}Q^{k+3+\ep}\#\{\lbf:\,|\lbf|\le L,\,
G(\lbf_3)=G(\lbf_4)=H_{rsi}(\lbf_3,\lbf_4)=0,\, (r,s\le k+1)\}.\]
Moreover, if $H_{rsi}(\lbf_3,\lbf_4)=0$ then
$H_{rs}\big((\mu_i^\sig,\lambda_i^\sig);\lbf_3,\lbf_4\big)=0$, for
one of the automorphisms $\sig$, and hence
$H_{rs}\big((\mu_i,\lambda_i);\lbf_3,\lbf_4\big)=0$. It follows that
\begin{eqnarray*}
\Sigma_1&\ll_{\ep}&Q^{k+3+\ep}\#\{\lbf:\,|\lbf|\le L,\,
G(\lbf_3)=G(\lbf_4)=H_{rs}\big((\mu_i,\lambda_i);\lbf_3,\lbf_4\big)
=0\\ &&\hspace{7cm}(r,s\le k+1)\}.
\end{eqnarray*}
According to Theorem 1 of Browning and Heath-Brown \cite{BHB05}, the
number of admissible vectors $\lbf$ is $O(L^{m+1})$, where $m$ is the
dimension of the projective variety defined over the field $K$ by
\begin{equation}\label{rvar}
\{[\xbf,\ybf]\in\Proj^{2k-1}:G(\b{x})=G(\b{y})=
H_{rs}\big((\mu_i,\lambda_i);\b{x},\b{y}\big)=0\;\;\;(r,s\le k+1)\}.
\end{equation}
We will then have
\begin{equation}\label{rres}
\Sigma_1\ll_{\ep}Q^{k+3+\ep}L^{m+1}.
\end{equation}
\label{note_thmTHB}

It remains to provide an upper bound for the dimension $m$. Both of the varieties $G(\b{x})=0$ and $G(\b{y})=0$ are unions of hyperplanes, $L_j(\b{x})=0$ or $L_j(\b{y})=0$ say, while the variety \[H_{rs}\big((\mu_i,\lambda_i);\b{x},\b{y}\big)=0\;\;\;(r,s\le k+1)\]
describes the vanishing of all the $k\times k$ minors of a matrix
\[B:=\left(\begin{array}{c|c} A & \b{y}\\ \hline
\rule[3mm]{0mm}{1mm}\b{x}^T & 0\end{array}\right).\]
Here $A=(\mu_i,\lambda_i)\cdot \u{Q}$ is a $k\times k$ symmetric
matrix, which is known to have rank
$k-1$ by construction. To examine this further we diagonalize $A$ over $\Qbar$ as $R^TDR$, say,
where $R$ is nonsingular, and $D={\rm diag}(0,d_2,\ldots, d_k)$ with $d_2\ldots d_k\not=0$. If $B$ has rank at most $k-1$ the matrix
\[\left(\begin{array}{c|c}
\rule[-3mm]{0mm}{1mm}R^{-T}&\b{0}\\ \hline
\rule[3mm]{0mm}{1mm}\b{0}^T & 1\end{array}\right)
\left(\begin{array}{c|c}\rule[-3mm]{0mm}{1mm}A &\b{y}\\ \hline
\rule[3mm]{0mm}{1mm}\b{x}^T & 0\end{array}\right)
\left(\begin{array}{c|c}\rule[-3mm]{0mm}{1mm}R&\b{0}\\ \hline
\rule[3mm]{0mm}{1mm}\b{0}^T & 1\end{array}\right)=
\left(\begin{array}{c|c}
\rule[-3mm]{0mm}{1mm}D&R^{-T}\b{y}\\ \hline \rule[3mm]{0mm}{1mm}\b{x}^TR & 0\end{array}\right)
\]
will also have rank at most $k-1$.
However if we write $\b{x}^TR=\b{u}$ and $R^{-T}\b{y}=\b{v}$ then the
$(1,1)$-minor of the above matrix is
\begin{equation}\label{rexp}
-u_2v_2\widehat{d_2}-\ldots-u_kv_k\widehat{d_k},
\end{equation}
where
\[\widehat{d_j}=\prod_{{\substack{2\le h\le k\\ h\not=j}}}d_h\not=0.\]
In terms of the coordinates $\b{u}$ and $\b{v}$ the intersection
$G(\b{x})=G(\b{y})=0$ is a union of linear spaces of codimension 2,
of the type $L(\b{u})=L'(\b{v})=0$.  However it is clear that
the expression (\ref{rexp}) cannot vanish on such a linear space,
providing that
$k\ge 4$. (One may verify that (\ref{rexp}) is a quadratic form in
$(\ubf,\vbf)$ of rank $2k-2$, and in general a quadratic form that
vanishes on a codimension $r$ linear space has rank $\leq 2r$.) It
therefore follows that the variety (\ref{rvar}) has codimension at least 3 in $\mathbb{P}^{2k-1}$, whence $m\le 2k-4$.
The bound (\ref{rres}) then yields
\[\Sigma_1\ll_{\ep}Q^{k+3+\ep}L^{2k-3}\]
which is satisfactory for Lemma \ref{ravlem}.
\label{note_linearspace}

We turn now to the case in which $F(\lbf)\not=0$, where we have set
\[F(\lbf):=G(\lbf_3)^2+G(\lbf_4)^2+H(\lbf_3,\lbf_4)^2\]
for convenience. Here we will write $\Sigma_2$ for the sum to be
estimated. We begin by claiming that Lemmas \ref{lemma_bad_p} and
\ref{prop_min_sum_bd0} may be combined to say that there is a constant
$C$ for which
\begin{equation}\label{rcl}
|S(\lbf;p^e)|\le
C(e+1)p^{e(k+2)}\{p^{-1}(p,F(\lbf))+p^{-2}f(\lbf;p)\}
\end{equation}
for all primes $p$, where we define
\[f(\lbf;q):=\#\{\u{b}\modd{q}:\, q\mid\det
M(\u{b};\lbf_3,\lbf_4)\}\]
for any square-free $q$.   Since $(p,F(\lbf))\ge 1$, the claim is plainly true for the bad primes, on
taking $C\ge p\,c_p$ for all such primes.  The claim is also straightforward in
the case in which $p$ is a good prime dividing each of $G(\lbf_3)$,
$G(\lbf_4)$ and $H(\lbf_3,\lbf_4)$, since then $(p,F(\lbf))=p$.
In the remaining case all that is needed is to observe that
\[\#\{\u{b}\modd{p^e}:\, p\mid\det M(\u{b};\lbf_3,\lbf_4)\}
=p^{2e-2}f(\lbf;p).\]
This establishes our claim.

We will find it convenient to introduce the notation
\[\kappa(q):=\prod_{p\mid q}p.\]
Since the function $f(\lbf;q)$ is multiplicative on square-free
integers $q$ we can now extend the estimate (\ref{rcl}) to general
moduli $q$ by multiplicativity, to give
\[|S(\lbf;q)|\le C^{\omega(q)}d(q)q^{k+2}\twosum{q_1q_2=q}{(q_1,q_2)=1}
\frac{(\kappa(q_1),F(\lbf))}{\kappa(q_1)}
\frac{f(\lbf;\kappa(q_2))}{\kappa(q_2)^2}.\]
We then deduce that
\[\Sigma_2\ll Q^{k+2+\ep}\twosum{|\lbf|\le L}{F(\lbf)\not=0}
\Sigma_3(\lbf)\Sigma_4(\lbf),\]
with
\[\Sigma_3(\lbf):=\sum_{q\le Q}\frac{(\kappa(q),F(\lbf))}{\kappa(q)}\]
and
\[\Sigma_4(\lbf):=\sum_{q\le Q}\frac{f(\lbf;\kappa(q))}{\kappa(q)^2}.\]
We estimate $\Sigma_3(\lbf)$ using Rankin's trick.  We have
\begin{eqnarray*}
\Sigma_3(\lbf)&\le&
Q^{\ep}\sum_{q=1}^{\infty}\frac{(\kappa(q),F(\lbf))}{q^{\ep}\kappa(q)}\\
&=&Q^{\ep}\prod_p\left\{
\sum_{e=0}^{\infty}\frac{(\kappa(p^e),F(\lbf))}{p^{e\ep}\kappa(p^e)}\right\}\\
&=&Q^{\ep}\prod_{p\mid F(\lbf)}\left\{1+p^{-\ep}+p^{-2\ep}+\cdots\right\}
\prod_{p\nmid F(\lbf)}\left\{1+p^{-1-\ep}+p^{-1-2\ep}+\cdots\right\}.
\end{eqnarray*}
Hence if \begin{equation}\label{rcdef}
c=c(\ep):=1+2^{-\ep}+4^{-\ep}+\cdots
\end{equation}
we will have
\[\Sigma_3(\lbf)\le Q^{\ep}c^{\omega(F(\lbf))}\prod_p\{1+cp^{-1-\ep}\}
\le Q^{\ep}c^{\omega(F(\lbf))}\zeta(1+\ep)^c\ll_{\ep}(QL)^{\ep}.\]
\label{note_Sig3}

\label{note_Sig4}
To handle $\Sigma_4(\lbf)$ we observe that for $q \leq Q$,
\[\left[\frac{Q}{\kappa(q)}\right]^2f(\lbf;\kappa(q))\le
\#\{b_1,b_2\le Q:\, \kappa(q)\mid\det M(\u{b};\lbf_3,\lbf_4)\}.\]
Since $[\theta]\ge \theta/2$ for any real $\theta\ge 1$ we deduce that
\[\Sigma_4(\lbf)\ll Q^{-2}\sum_{b_1,b_2\le Q}
\#\{q\le Q:\, \kappa(q)\mid\det M(\u{b};\lbf_3,\lbf_4)\}.\]
Let $\Sigma_5(\lbf)$ be the contribution to $\Sigma_4(\lbf)$ from terms in which
$\Delta :=\det M(\u{b};\lbf_3,\lbf_4)$ is non-zero; and write
$\Sigma_6(\lbf)$ for the contribution in which $\det M(\u{b};\lbf_3,\lbf_4)=0$.
For $\Sigma_5(\lbf)$ we use Rankin's
trick again, which shows that
\begin{eqnarray*}
\#\{q\le Q:\, \kappa(q)\mid\det M(\u{b};\lbf_3,\lbf_4)\}&\le&
Q^{\ep}\sum_{\substack{q=1\\ \kappa(q)\mid\Delta}}^{\infty}q^{-\ep}\\
&=& Q^{\ep}\prod_{p\mid\Delta}\left\{1+p^{-\ep}+p^{-2\ep}+\cdots\right\}\\
&\le & Q^{\ep}c^{\omega(|\Delta|)}
\end{eqnarray*}
with $c$ as in (\ref{rcdef}).  Since $\Delta$ is bounded by a suitable
power of $QL$ we deduce that
\[\#\{q\le Q:\, \kappa(q)\mid\det M(\u{b};\lbf_3,\lbf_4)\}\ll_{\ep}
Q^{2\ep}L^{\ep}\]
so that \[\Sigma_5(\lbf)\ll_{\ep}Q^{2\ep}L^{\ep}.\]

It remains to deal with $\Sigma_6(\lbf)$, for which we clearly have
\[\Sigma_6(\lbf)\ll Q^{-1}\#\{b_1,b_2\le Q:\,
\det M(\u{b};\lbf_3,\lbf_4)=0\}.\]
Combining our various estimates we now see that
\begin{eqnarray}\label{rs1}
\Sigma&\ll_{\ep}&Q^{k+3+\ep}L^{2k-3}+Q^{k+2+4\ep}L^{2k+2\ep}\nonumber\\
&&\hspace{3mm}\mbox{}+
Q^{k+1+2\ep}L^{\ep}\sum_{|\lbf|\le L}\#\{b_1,b_2\le Q:\,
\det M(\u{b};\lbf_3,\lbf_4)=0\}.
\end{eqnarray}
For each fixed $\lbf$ the expression $\det M(\u{b};\lbf_3,\lbf_4)$
is a binary form in $b_1,b_2$ of degree $k-1$.  It follows that there
are at most $O(Q)$ pairs $b_1,b_2\le Q$ for which $\det M(\u{b};\lbf_3,\lbf_4)=0$, unless the form vanishes identically.
We therefore have
\begin{equation}\label{rs2}
\sum_{|\lbf|\le L}\#\{b_1,b_2\le Q:\,\det
M(\u{b};\lbf_3,\lbf_4)=0\}\ll QL^{2k}+Q^2n(L),
\end{equation}
where
\[n(L):=\#\{\lbf:|\lbf|\le L,\,\det M(\u{x};\lbf_3,\lbf_4)\equiv 0\},\]
so that $n(L)$ counts vectors $\lbf$ for which $\det
M(\u{x};\lbf_3,\lbf_4)$ vanishes with respect to $\u{x}$.
\label{note_b}

To analyse $n(L)$ we consider the choice $\u{x}=(\mu_i,\lambda_i)$,
using the notation from \S \ref{rss1}.  This produces a matrix
$(\mu_i,\lambda_i)\cdot \u{Q}$ of rank exactly $k-1$. Suppose that
\[M:=\left(\begin{array}{c|c} \rule[-3mm]{0mm}{1mm}
(\mu_i,\lambda_i)\cdot \u{Q} & \lbf_4\\ \hline
\rule[3mm]{0mm}{1mm}\lbf_3^T & 0\end{array}\right)\]
is singular. If $\lbf_4$ is not in the column space for
$(\mu_i,\lambda_i)\cdot \u{Q}$ then
\[M_0:=\left(\begin{array}{c|c} (\mu_i,\lambda_i)\cdot \u{Q} & \lbf_4\end{array}\right)\]
will have rank $k$, and hence has linearly independent rows. Since $M$ is singular it would follow that $(\lbf_3^T\mid 0)$ is in the row space for $M_0$, whence in
particular $\lbf_3^T$ would be in the row space for
$(\mu_i,\lambda_i)\cdot \u{Q}$. We therefore conclude that either $\lbf_3$ or $\lbf_4$ must lie in
the column space of $(\mu_i,\lambda_i)\cdot \u{Q}$.  In \S \ref{rss1} we chose
a null vector $\b{e}_i$ for $(\mu_i,\lambda_i)\cdot \u{Q}$, and we now see that
either $\b{e}_i^T\lbf_3=0$ or $\b{e}_i^T\lbf_4=0$.  This conclusion
is valid for any index $i\le k$ so that for each vector $\lbf$
counted by $n(L)$ there will be a subset
$S\subseteq\{1,\ldots,k\}$ for which
\[\b{e}_i^T\lbf_3=0\;\;\; (i\in S)\;\;\;\mbox{and}\;\;\; \b{e}_i^T\lbf_4=0\;\;\; (i\not\in S).\]
These conditions restrict $\lbf$ to a $k$-dimensional subspace of $\Z^{2k}$.
Any such subspace can contain at most
$O(L^k)$ integral points with $|\lbf|\le L$, whence $n(L)\ll L^k$.
It then follows from (\ref{rs1}) and (\ref{rs2}) that
\[\Sigma\ll_{\ep}Q^{k+3+\ep}L^{2k-3}+Q^{k+2+4\ep}L^{2k+2\ep}+
Q^{k+1+2\ep}L^{\ep}(QL^{2k}+Q^2L^k),\]
which suffices for Lemma \ref{ravlem}, on redefining $\ep$.

\section{Exponential sums: The singular series}\label{sec_ss}

We will now use our results about exponential sums to establish some
key results about the singular series, given by
\beq\label{sing_ser_dfn}
 \Sf (\u{n})
	=\sum_{ q =1}^\infty  \frac{1}{q^k} \sum_{\substack{1 \leq
           a_1,a_1 \leq q\\ (a_1,a_2,q)=1}} S_q(\u{a};\u{n})
=\sum_{ q =1}^\infty  \frac{1}{q^k}T(\u{n};q),
	\eeq
with
\[S_q(\underline{a};\underline{n}) = S_q(\u{a})e_q(-\underline{a}
\cdot \underline{n}) .\]

We will prove the following propositions.
\begin{prop}\label{prop_sing_ser_conv}
The singular series $\Sf(\u{n})$ is absolutely convergent for every
$\u{n}$ for which
$F(n_2,-n_1)\not=0$, providing that $k \geq 5$. Indeed
\beq\label{Rcon}
\sum_{q\ge R}\frac{1}{q^k} |T(\u{n};q)|\ll |\u{n}|^\ep R^{-1/3},
\eeq
for such $\u{n}$, where the implied constant depends only on $Q_1,Q_2$ and $\ep$. \end{prop}
\begin{prop}\label{prop_sing_ser_conv2}
There is a constant $p_0$ depending only on $Q_1$ and $Q_2$ with the
following property.  Let $k\geq 5$ and $F(n_2,-n_1)\not=0$. Suppose further that the
system $\u{Q}(\xbf)=\u{n}$ is locally solvable in $\Z_p$ for every
prime $p$. Then for any $\ep>0$ we will have
\[\Sf(\u{n})\gg |\u{n}|^{-\ep}\prod_{p\le p_0}|F(n_2,-n_1)|_p^{k-2}.\] \end{prop}

Of course (\ref{Rcon}) yields not only the statement about absolute
convergence, by taking $R=1$, but also the estimate
\[\Sf(\u{n})\ll|\u{n}|^{\ep}\]
for $F(n_2,-n_1)\not=0$.

\label{note_sumR}
We begin by establishing (\ref{Rcon}). For any $R \geq 1$ we have
\[ \sum_{q\ge R}q^{-k} |T(\u{n};q)|\le R^{-1/3}
\sum_{q=1}^{\infty}q^{1/3-k} |T(\u{n};q)|.\]
Let
\[\psi_p:=\sum_{e=0}^{\infty}p^{-e(k-1/3)}|T(\u{n};p^e)|,\]
so that
\[ \sum_{q\ge R}q^{-k} |T(\u{n};q)|\le R^{-1/3}\prod_p \psi_p\]
by multiplicativity. By Proposition \ref{T_I_lemma} we have
\[\psi_p=1+O(p^{4/3-k/2})=1+O(p^{-7/6})\]
for Type I primes, and
\begin{eqnarray}\label{13x}
\psi_p&=&1+O(p^{11/6-k/2})+O(\sum_{e=2}^{\infty}p^{(7/3-k/2)e})\nonumber\\
&=& 1+O(p^{11/6-k/2})+O(p^{14/3-k})\\
&=&O(1)\nonumber
\end{eqnarray}
for Type II primes. Finally, for bad primes, we will have
\[\psi_p=1+O_p(\sum_{e=1}^{\infty}p^{(7/3-k/2)e})\ll_p1.\]
The product of $\psi_p$ for Type I primes is thus $O(1)$, and similarly
for bad primes, since the collection of bad primes finite and is determined
purely by $Q_1$ and $Q_2$.  Finally, for the Type II primes, if
$|\psi_p|\le C$ say for such primes $p$, then the corresponding
product is at most
\[C^{\omega(F(n_2,-n_1))}\ll|\u{n}|^{\ep}.\]
The bound (\ref{Rcon}) then follows.

Turning to Proposition \ref{prop_sing_ser_conv2}, we begin by
observing that
\beq\label{Rpd}
\Sf(\u{n})=\prod_p\sigma_p
\eeq
with \[\sigma_p=1+\sum_{e=1}^{\infty}p^{-ek}T(\u{n};p^e).\]
By another application of Proposition \ref{T_I_lemma} we have
\beq\label{x1}
\sigma_p=1+p^{-k}T(\u{n};p)=1+O(p^{1-k/2})=1+O(p^{-3/2})
\eeq
for Type I primes and
\begin{eqnarray}\label{x2}
\sigma_p&=&1+p^{-k}T(\u{n};p)+\sum_{e=2}^{\infty}p^{-ek}T(\u{n};p^e)\nonumber\\
&=&1+O(p^{(3-k)/2})+O(\sum_{e=2}^{\infty}p^{(2-k/2)e})\nonumber\\
&=&1+O(p^{(3-k)/2})+O(p^{4-k})\nonumber\\
&=&1+O(p^{-1})
\end{eqnarray}
for Type II primes.  Suppose that $\sigma_p\ge 1-Ap^{-3/2}$ for some
explicit constant $A\ge 1$, for Type I primes.  Then $\sigma_p\ge (1-p^{-3/2})^{2A}$
for Type I primes $p\ge 2A$, since we have $1-At\ge(1-t)^{2A}$ for any
positive real $t\le (2A)^{-1}$.  Similarly, if $\sigma_p\ge 1-Ap^{-1}$
for Type II primes we will have $\sigma_p\ge (1-p^{-1})^{2A}$ for
$p\ge 2A$.  The contribution of such primes to the product (\ref{Rpd})
is therefore
\[\ge \prod_{p}(1-p^{-3/2})^{2A}\prod_{p|F(n_2,-n_1)}2^{-2A}
\gg d(|F(n_2,-n_1)|)^{-2A}\gg |\u{n}|^{-\ep},\]
where $d(*)$ is the usual divisor function.
\label{note_A}

It remains to consider bad primes, along with the remaining primes $p\le 2A$. Let $p_0$ be the largest of all these primes, so that $p_0$
depends only on the original forms $Q_1$ and $Q_2$.  Then according to
Proposition \ref{prop_lower_bd} we have
\[\sig_p(\u{n}) \geq \varpi_p |F(n_2,-n_1)|_p^{k-2}\]
for $p\le p_0$. The required lower bound then follows since
$\varpi_p \gg 1$ for $p\le p_0$.

In the special case $(n_1,n_2)=(0,0)$ similar results continue to hold.
\begin{prop}\label{prop_sing_ser_conv0}
The singular series $\Sf(\u{0})$ is absolutely convergent providing that $k \geq 6$. Indeed we then have
\beq\label{Rcon0}
\sum_{q\ge R}\frac{1}{q^k} |T(\u{0};q)|\ll R^{-1/3},
\eeq
where the implied constant depends only on $Q_1$ and $Q_2$. \end{prop}
\begin{prop}\label{prop_sing_ser_conv20}
Let $k\geq 6$ and suppose further that for each prime $p$ the system
$\u{Q}(\xbf)=\u{0}$ has a nonzero solution $\xbf_p \in \Z_p^k$.  Then
$\Sf(\u{0})>0$.
\end{prop}
Proposition \ref{prop_sing_ser_conv20} is an immediate consequence of
Propositions \ref{prop_lower_bd0} and \ref{prop_sing_ser_conv0}.

To prove Proposition \ref{prop_sing_ser_conv0} we note that there are
no Type I primes, and that
\[\psi_p=1+O(p^{11/6-k/2}) + O(p^{14/3-k})\]
for Type II primes, as in (\ref{13x}).  Thus
$\psi_p=1+O(p^{-7/6})$ for $k\ge 6$ and all good primes $p$.  Moreover
$\psi_p=O_p(1)$ for bad primes $p$ as before.  One may then establish
(\ref{Rcon0}) by exactly the same argument used for (\ref{Rcon}).
Finally we note that the absolute convergence of $\Sf(\u{0})$ is a
direct consequence of (\ref{Rcon0}).

\section{Application of the circle method}\label{sec_major}
We now begin our detailed application of the circle method.

\subsection{Division into major and minor arcs}

Let $Q=B^\Del$ for some fixed $\Del >0$ to be chosen later. We define the box \[ I(a_1,a_2;q)  = \left[\frac{a_1}{q} - \frac{Q}{B^2}, \frac{a_1}{q}
 + \frac{Q}{B^2}\right] \times  \left[\frac{a_2}{q} - \frac{Q}{B^2},
 \frac{a_2}{q} + \frac{Q}{B^2} \right]\]
for any integers $1 \leq a_1, a_2 \leq q \leq Q$.
Such boxes are disjoint provided that
\[ Q=B^{\Del}, \qquad \Del< 2/3,\]
since if $(a_1,a_2;q) \neq (a_1',a_2';q')$, then for $j=1$ or 2 we have
\[ \left| \frac{a_j}{q} - \frac{a_j'}{q'} \right| \geq \frac{1}{qq'}
\geq \frac{1}{Q^2} > 2\frac{Q}{B^2},\] for large enough $B$.  We now define the major arcs to be \[ \Mf(\Delta) = \Union_{1 \leq q \leq Q} \Union_{\substack{1 \leq a_1,
   a_2 \leq q\\ (a_1,a_2,q)=1}} I(a_1,a_2;q)\] and we take the minor arcs to be the complement of the major arcs in $[0,1]^2$:
\[ \mf(\Delta) = [0,1]^2 \setminus \Mf (\Delta).\]

\subsection{The major arcs: the singular integral and singular series}
Recall the definition
\[S(\al_1,\al_2) = \sum_{\x \in \Z^k} e(\al_1 Q_1(\xbf) + \al_2 Q_2(\xbf))w_B(x),\]
and the representation function
\[ R_B(n_1,n_2) = \iint_{[0,1]^2} S(\al_1,\al_2) e(-\u{\al} \cdot
\u{n}) d\al_1 d\al_2 .\] The contribution of the major arcs to $R_B(n_1,n_2)$ may now be represented as
\begin{eqnarray*}\label{major_int_cont}
\lefteqn{\iint_{\Mf(\Delta)} S(\al_1,\al_2) e(-\u{\al} \cdot \u{n}) d\al_1
d\al_2}\hspace{1cm} \nonumber\\
&= &\sum_{1 \leq q \leq Q} \sum_{\substack{1 \leq a_1,a_1 \leq
           q\\ (a_1,a_2,q)=1}} \iint_{I(a_1,a_2;q)} S(\al_1,\al_2) e(-\u{\al} \cdot \u{n}) d\al_1 d\al_2.
\end{eqnarray*}
Our goal is to approximate this by a main term of size $B^{k-4}$,
times a singular integral and a singular series. \label{note_Del4}

\begin{prop}\label{prop_major_int}
Suppose either that $k\ge 5$ and $|\u{n}|\ll B^2$ with $F(n_2,-n_1)\not=0$, or that $k\ge 6$ and $(n_1,n_2)=(0,0)$. Then for any fixed positive $\Delta\le 1/8$
we have
\[ \iint_{\Mf(\Delta)} S(\u{\al}) e(-\u{\al} \cdot \u{n})d\al_1
d\al_2= \Sf(\u{n}) \Jcal_w(B^{-2}\u{n}) B^{k-4} + E(\u{n}),\] with
\beq\label{eqn_E}
E(\u{n}) \ll B^{k-4 - \Del/4},
\eeq
where $\Sf(\u{n})$ and $\Jcal_w(\u{\mu})$ are given by
(\ref{sing_ser_dfn}) and (\ref{sing_int_dfn}) respectively.
\end{prop}

We will first prove by a standard argument that:

\begin{lemma}\label{xxx}
If $\u{\al}$ belongs to one of the major arcs $I(a_1,a_2;q)$, and
$\u{\theta} = \u{\al} - \u{a}/q$, then \beq\label{S_2terms}
S(\al_1,\al_2) = q^{-k} B^k S_q(a_1,a_2) I(B^2 \u{\theta})+O(B^{k-1+2\Delta}), \eeq
with $I(\u{\phi})$ given by (\ref{dfn_I}).
\end{lemma}

Write $\theta_j = \al_j - a_j/q$ for $j=1,2$.  Then
\begin{eqnarray}
S(\al_1,\al_2) &=& \sum_{\x \in \Z^k} e(\u{\al} \cdot
\u{Q}(\xbf))w_B(\xbf)\nonumber \\	&=& \sum_{\ybf \modd{q}}
\sum_{\zbf \in \Z^k} e(\u{\al} \cdot \u{Q}(q\zbf + \ybf)) w_B(q\zbf+
\ybf) \nonumber \\
	& = &  \sum_{\ybf \modd{q}} e_q(\u{a} \cdot \u{Q}(\ybf)
       )\sum_{\zbf \in \Z^k} f(\zbf), \label{S_sum}
	\end{eqnarray}
where
\[ f(\zbf) = 	 e(\u{\theta} \cdot \u{Q}(q\zbf+ \ybf))  w_B(q\zbf + \ybf) .\] Note that $f(\zbf)$ is supported in a $k$-dimensional cube $K$ centred
at the origin, with
side-length of order $B/q+1 \ll B/q$. We will now replace the summation over $\zbf$ by integration over a
continuous variable, incurring a small error in the process. Note
first that for any $\wbf  \in [0,1]^k$ we have
\[ |f(\zbf + \wbf) - f(\zbf)| \leq k\max_{\ubf \in [0,1]^k} |\nabla
f(\zbf + \ubf)| \] by the mean-value theorem.  Thus \begin{eqnarray*}
	| \int_{\R^k} f(\zbf) d\zbf - \sum_{\zbf \in \Z^k} f(\zbf)|
		&\ll& (B/q)^k \max_{\zbf \in K} |\nabla f(\zbf)| \\
		&\ll & (B/q)^k \max_{\zbf \in K} (q/B + q|\u{\theta}|.
               |q\zbf+\ybf|) \\
		&\ll &  q^{1-k}B^{k-1}+ |\u{\theta}| q^{1-k}B^{k+1} .
		\end{eqnarray*}
Consequently, \[ \sum_{\zbf \in \Z^k} f(\zbf) = \int_{\R^k} e(\u{\theta} \cdot
\u{Q}(q\zbf+ \ybf) ) w_B(q\zbf+ \ybf)  d \zbf + O( |\u{\theta}|
q^{1-k}B^{k+1} + q^{1-k}B^{k-1}),\] which upon setting $B \xbf = \zbf q + \ybf$ becomes \[ \sum_{\zbf \in \Z^k} f(\zbf) = \frac{B^k}{q^k} \int_{\R^k} e(B^2
\u{\theta} \cdot \u{Q}(\xbf) )w(\xbf)  d \xbf + O( |\u{\theta}|
q^{1-k}B^{k+1} + q^{1-k}B^{k-1}).\] Applying this to the innermost sum in (\ref{S_sum}), we then see that
\beq\label{S_2terms'}
S(\al_1,\al_2) = q^{-k} B^k S_q(a_1,a_2) I(B^2\u{\theta}) + O
(B^{k-1}q(|\u{\theta}| B^2 +1)). \eeq
This proves the lemma, upon noting that $|\u{\theta}| \leq
B^{-2+\Delta}$ and $q \leq B^{\Delta}$ in the major arcs, so that the
error term is no more than $O(B^{k-1 + 2\Delta}).$

Our goal is now to integrate $S(\al_1, \al_2)$ over the full
collection of major arcs. Note that the
measure of the total collection of major arcs is \[\ll B^{\Delta} \cdot B^{2\Delta} \cdot (B^{-2 + \Delta})^2
\ll B^{-4+ 5\Delta}.\] Thus Lemma \ref{xxx} immediately implies: \begin{lemma}\label{lemma_major_R}
We have
\begin{eqnarray*}
\iint_{\Mf(\Delta)} S(\al_1,\al_2) e(-\u{\al} \cdot \u{n}) d\al_1
d\al_2&=& B^{k-4}\Jcal_w(B^{-2}\u{n};B^\Del)\sum_{q\le B^\Delta}q^{-k}T(\u{n};q)\\
&&\hspace{2cm}  + O (B^{k-5 +7\Delta})
\end{eqnarray*}
with $\Jcal_w(\u{n};R)$ given by (\ref{sing_int_dfn+}).
\end{lemma}

Finally, we apply the results of Propositions \ref{prop_sing_int_conv}
and \ref{prop_sing_ser_conv} to the truncated singular integral and
singular series in order to pass to the limit on the right hand
side. We obtain: 
\begin{eqnarray*}
\iint_{\Mf(\Delta)} S(\al_1,\al_2) e(-\u{\al} \cdot \u{n})d\al_1
d\al_2&=& \Sf(\u{n}) \Jcal_w(B^{-2}\u{n}) B^{k-4}\\ &&\hspace{5mm}\mbox{}+O(B^{k-4 - \Del/3+\varepsilon})+O(B^{k-5 + 7\Del}),
\end{eqnarray*}
for $F(n_2,-n_1)\not=0$,
as long as $k \geq 5$. Similarly, if $k\ge 6$ and $(n_1,n_2)=(0,0)$
we may apply Propositions \ref{prop_sing_int_conv} and
\ref{prop_sing_ser_conv0}. Proposition \ref{prop_major_int} then follows, upon restricting $\Del \leq 1/8$ and replacing
$B^{-\Del/3+\ep}$ by $B^{-\Del/4}$.

\section{Proof of Theorem \ref{thm_r1}}\label{sec_thm_proofs} 

\subsection{The mean square argument}
We are now ready to make precise the mean square argument sketched in
\S \ref{sec_method}.  Proposition \ref{prop_major_int} establishes
(\ref{maj_ME}) with \[M(n_1,n_2)=\Sf(\u{n}) \Jcal_w(B^{-2}\u{n}) B^{k-4}\]
and
$E(n_1,n_2)=O(B^{k-4 - \Del/4})$, whence (\ref{min_E}) yields
\begin{eqnarray*}
\lefteqn{\sum_{\substack{\max(|n_1|,|n_2|)\le N\\ F(n_2,-n_1)\not=0}} |R_B(n_1,n_2) -\Sf(\u{n})\Jcal_w(B^{-2}\u{n})B^{k-4}|^2}\hspace{2cm}\\
&\ll&  \iint_{\mf(\Delta)} |S(\al_1,\al_2)|^2 d\al_1 d\al_2 +
B^{2k-4-\Delta/2}, \end{eqnarray*}
provided that $N\ll B^2$.

The crucial result is the following upper bound for the minor arcs integral.
\begin{prop}\label{prop_minor_arcs}
For any $k  \geq 5$, any $\ep>0$, and any $\Del\in(0,1/6)$, we have
\beq\label{min_int}
\iint_{\mf(\Delta)} |S(\al_1,\al_2)|^2d\al_1 d\al_2  \ll
B^{2k-4-2\Del + \ep}.
\eeq
\end{prop}
The choice $\Delta=1/8$ then establishes Theorem \ref{thm_r1}.
\label{note_16}

The proof of Proposition \ref{prop_minor_arcs} is the most delicate part of the
paper.  We begin by employing
a 2-dimensional Dirichlet approximation with a parameter $S \geq 1$. Thus for every pair
$\al_1, \al_2$ in $[0,1]$ there exist $1 \leq q \leq S$ and $1 \leq
a_1, a_2 \leq q$ with $(a_1,a_2,q)=1$ such that \[ \left| \al_1 - \frac{a_1}{q} \right| \leq \frac{1}{ q\sqrt{S}}, \quad \text{and} \quad \left| \al_2 - \frac{a_2}{q} \right| \leq \frac{1}{ q\sqrt{S}}.\]
Given $\al_1,\al_2$ and approximations $\al_1 = a_1/q + \theta_1,
\al_2 = a_2/q+ \theta_2$ of the above type, then if
$\al_1,\al_2 \in \mf(\Delta)$ at least one of the inequalities \beq\label{maj_cond}
q \leq B^\Delta, \qquad |\theta_1| \leq B^{-2 + \Delta}, \qquad
|\theta_2| \leq B^{-2 + \Delta},
\eeq
must fail to hold. For our application we shall choose \[S=B^{4/3},\]
which is essentially optimal.

We bound the integral (\ref{min_int}) from above, using a collection of dyadic sums
\beq\label{Sig_dfn}
\Sig (R, \phi_1,\phi_2) = \sum_{R \leq q < 2R}
\sum_{\substack{1 \leq a_1,a_2 \leq q\\ (a_1,a_2,q)=1}}
\iint_{\{\phi_1,\phi_2\}} |S(a_1/q+\theta_1,a_2/q+\theta_2)|^2
d\u{\theta}.
\eeq
The reader should recall that$\iint_{\{\phi_1,\phi_2\}}$ denotes an
integral over the range
\[([-2\phi_1, -\phi_1]\cup [\phi_1, 2\phi_1]) \times ([-2\phi_2, -\phi_2]\cup [\phi_2, 2\phi_2]).\] We will prove:
\begin{prop}\label{prop_sig_bd}
For any $k \geq 5$, any $\ep>0$, and any $\Del\in(0,1/6)$, we have
\beq\label{Sig_R}
 \Sig (R, \phi_1,\phi_2)  \ll B^{2k-4-2\Delta+\ep}
 \eeq
for $R\ll B^{4/3}$ and $\phi_1,\phi_2\ll R^{-1}B^{-2/3}$,
unless all three conditions \beq\label{maj_cond'}
R \leq \frac{1}{2}B^\Delta, \qquad \phi_1 \leq \frac{1}{2}B^{-2 + \Delta}, \qquad \phi_2
\leq \frac{1}{2}B^{-2 + \Delta}, \eeq
hold. \end{prop}

Before proving Proposition \ref{prop_sig_bd} we show how it implies Proposition \ref{prop_minor_arcs}. When $q\ge B^{\Delta}$ we handle
the squares \[\{ (\al_1,\al_2) = (a_1/q + \theta_1, a_2/q +
\theta_2): \max(|\theta_1|,|\theta_2|) \leq B^{-k} \}\]
by a trivial estimate, producing an overall contribution $\ll
S^3B^{-2k}.B^{2k}$ to (\ref{min_int}), since we trivially have $S(\al_1,\al_2)\ll B^k$. When $q\le B^{\Delta}$ we know that
\[\max(|\theta_1|,|\theta_2|)\ge B^{-2+\Delta}\]
since at least one of
the conditions (\ref{maj_cond}) is known to fail. The remaining cases
may then be covered by $O((\log B)^3)$ dyadic
intervals for $q$, $\theta_1$ and $\theta_2$.

Proposition \ref{prop_sig_bd} now yields
\[ \iint_{\mathfrak{m}(\Del)} |S(\al_1,\al_2)|^2 d\al_1 d\al_2 \ll
S^3+(\log B)^3 \sup \Sig(R,\phi_1,\phi_2),\] where the supremum is taken over all dyadic parameters with $0<R< S$
and $B^{-k} \leq \phi_1,\phi_2 \leq (R\sqrt{S})^{-1}$ such that not
all three conditions (\ref{maj_cond'}) hold. This clearly suffices for
Proposition \ref{prop_minor_arcs}, with the choice $S=B^{4/3}$.

\subsection{Proof of Proposition \ref{prop_sig_bd}}

Recalling the definition (\ref{S_dfn1}) of $S(\al_1,\al_2)$, we may
expand the integrand in (\ref{Sig_dfn}) and write $\xbf_j = \lbf_j q +
\rbf_j$, where $\lbf_j \in \Z^{k}$ and $\rbf_j \in (\Z/q\Z)^{k}$ for
$j=1,2$. This produces
\[|S(a_1/q+\theta_1,a_2/q+\theta_2)|^2=
\sum_{\substack{\rbf_1 \modd{q}\\ \rbf_2\modd{q}}} e_q(\u{a}\cdot\u{Q}(\rbf_1)-\u{a}\cdot\u{Q}(\rbf_2))
\Sigma(\u{\theta},\rbf_1,\rbf_2,q),\]
where we have temporarily set
\[\Sigma(\u{\theta},\rbf_1,\rbf_2,q)=
\sum_{\lbf_1, \lbf_2 \in \Z^k} e(\u{\theta} \cdot
\u{Q} (\lbf_1q + \rbf_1) - \u{\theta} \cdot \u{Q} (\lbf_2q +
\rbf_2)) w_B(\lbf_1q+\rbf_1)w_B(\lbf_2q+\rbf_2).\]
Applying Poisson summation to the sum over $\lbf_1,\lbf_2$,  we may
rewrite this as
\[\Sigma(\u{\theta},\rbf_1,\rbf_2,q)=(\frac{B}{q})^{2k}
\sum_{\lbf \in \Z^{2k}}e_q(\rbf \cdot \lbf)J(B^2\u{\theta},B\lbf/q),\]
where for any $\lambf\in \R^{2k}$ we have defined \[ J(\u{\nu},\lambf)  = \iint_{\R^{2k}} e( \unu \cdot \u{Q}(\ubf_1)
- \unu \cdot \u{Q}(\ubf_2))  w(\ubf_1) w(\ubf_2) e(-\ubf \cdot \lambf)
d\ubf.\] It follows that
\[|S(a_1/q+\theta_1,a_2/q+\theta_2)|^2=(\frac{B}{q})^{2k}
\sum_{\lbf \in \Z^{2k}}S(\lbf;q)J(B^2\u{\theta},B\lbf/q)\]
with
\[S(\lbf;q) = \sum_{\substack{1 \leq a_1,a_2 \leq q\\ (a_1,a_2,q)=1}}
\sum_{\substack{\rbf_1 \modd{q}\\ \rbf_2 \modd{q}}} e_q( \u{a} \cdot
\u{Q} (\rbf_1)- \u{a} \cdot \u{Q}(\rbf_2))e_q(\rbf \cdot \lbf).\]
We then conclude that
\beq\label{Sig_SI}
\Sig (R, \phi_1,\phi_2) = B^{2k}\sum_{R \leq q < 2R}
\frac{1}{q^{2k}}	\sum_{\lbf \in \Z^{2k}}S(\lbf;q)
\Ical_{\{\phi_1,\phi_2\}} (\lbf;q),
\eeq
where
\[\Ical_{\{\phi_1,\phi_2\}} (\lbf;q) = \iint_{\{\phi_1,\phi_2\}} J(B^2\u{\theta},B\lbf/q)d\u{\theta}= B^{-4}\iint_{B^2\{\phi_1,\phi_2\}} J(\u{\nu},B\lbf/q)d\u{\nu}.\]
We now observe that $J(\u{\nu},\lambf)$ is essentially the integral $I(\u{\nu}\cdot\u{F};\lambf)$
occurring in Lemma \ref{I_int_avg2}.  Indeed Lemma \ref{I_int} shows
that
\[J(\u{\nu},\lambf)\ll_M |\lambf|^{-M}\]
for any fixed $M>0$, when $|\lambf|\gg |\u{\nu}|$.  We may therefore
deduce the following bound, via Lemma \ref{I_int_avg2}.
\begin{lemma}\label{J_lemma}
Let $\phi^* = \max( \phi_1,\phi_2)$. Then
\beq\label{I_bd_cases1}
\iint_{B^2\{ \phi_1,\phi_2\}} J(\u{\nu},B \lbf /q)d \u{\nu} \ll
(B^2\phi^*)^2 \min(1, (B^2\phi^*)^{-k})\log B.
\eeq
Moreover, for any $M>0$ we have \beq\label{I_bd_cases12}
\iint_{B^2\{ \phi_1,\phi_2\}} J(\u{\nu},B \lbf /q)d \u{\nu} \ll_{M}
	 (B^2\phi^*)^2 (B|\lbf|/R)^{-M},
\eeq
if $|\lbf| \gg  RB\phi^*$.
\end{lemma}

We trivially have $S(\lbf;q)\ll q^{2k+2}$, so on writing \[L=RB^{-1+\ep}(1+ B^2 \phi^*)\]
with a small $\ep>0$, and assuming that $M>2k+1$, we see that
\begin{eqnarray*}
\lefteqn{B^{2k}\sum_{R \leq q < 2R} \frac{1}{q^{2k}}	\sum_{|\lbf|\ge L}|S(\lbf;q)
\Ical_{\{\phi_1,\phi_2\}} (\lbf;q)|}\hspace{2cm}\\
& \ll_{\ep,M} &B^{2k-4} R^3(B^2\phi^*)^2(R/B)^M\sum_{|\lbf|\ge L}|\lbf|^{-M}\\
&\ll_{\ep,M} & B^{2k} R^3(R/B)^M L^{2k+1-M}\\
&\ll_{\ep,M} &B^{2k} R^3(R/B)^M (RB^{-1+\ep})^{2k+1-M}\\
&\ll_{\ep,M} &B^{2k}R^3(RB^{-1+\ep})^{2k+1}B^{-\ep M}.
\end{eqnarray*}
Thus, taking $M$ as a suitably large multiple of $\ep^{-1}$, we see
that terms with $|\lbf|\ge L$ contribute $O(1)$ to $\Sig(R,
\phi_1,\phi_2)$. This is satisfactory for Proposition
\ref{prop_sig_bd}.

We dispose next of the term $\lbf=\mathbf{0}$.  As noted in (\ref{S0}),
\[S(\mathbf{0};q)\ll q^{k+2+\ep}\]
for any $\ep>0$.  Then, applying (\ref{I_bd_cases1}) with $\lbf=\mathbf{0}$
we conclude that the contribution to $\Sig(R,\phi_1,\phi_2)$ is
\begin{eqnarray*}
&\ll& B^{2k-4}(B^2\phi^*)^{2}\min(1,(B^2\phi^*)^{-k})(\log B)\sum_{R\le
 q<2R}q^{-2k}q^{k+2+\ep}\\
&\ll&
B^{2k-4}\min((B^2\phi^*)^{2},(B^2\phi^*)^{2-k})(\log B)R^{-k+3+\ep}.
\end{eqnarray*}
When $R\le \frac{1}{2} B^{\Delta}$ we have $\phi^*\ge
\frac{1}{2}B^{-2+\Delta}$ and the above is
\[\ll B^{2k-4}(B^2\phi^*)^{2-k}(\log B)R^{-k+3+\ep}
\ll B^{2k-4}B^{\Delta(2-k)}(\log B).\]
This too is satisfactory for Proposition \ref{prop_sig_bd}.
Similarly when $R\ge \frac{1}{2} B^{\Delta}$ we see that our bound becomes
\[\ll B^{2k-4}(\log B)R^{-k+3+\ep}\ll
B^{2k-4}(\log B)B^{\Delta(-k+3+\ep)}\]
which again is satisfactory.
\label{note_BR}

It remains to handle the range $1\le|\lbf|\le L$. Using the bound (\ref{I_bd_cases1}) we see that the
contribution to $\Sig(R,\phi_1,\phi_2)$  will be
\beq\label{REST}
\ll B^{2k-4}\min((B^2\phi^*)^{2},(B^2\phi^*)^{2-k})(\log B)
\sum_{R\le q<2R}q^{-2k}\sum_{1\le|\lbf|\le L}|S(\lbf;q)|.
\eeq
We may now apply Lemma \ref{ravlem}, which shows the above to be
\begin{eqnarray}\label{REST1}
\lefteqn{\ll B^{2k-4}\min((B^2\phi^*)^{2},(B^2\phi^*)^{2-k})(\log B)
R^{-2k}}\hspace{3cm}\nonumber\\
&& \times (R^{k+2+\ep}L^{2k+\ep}+R^{k+3+\ep}L^{2k-3+\ep}).
\end{eqnarray}
When $B^2\phi^*\ge 1$ we have $L\ll RB^{1+\ep}\phi^*$ and this becomes
\begin{eqnarray}
&\ll &B^{2k-4+3k\ep}(B^2\phi^*)^{2-k}
R^{-2k}(R^{k+2}(RB\phi^*)^{2k}+R^{k+3}(RB\phi^*)^{2k-3})\nonumber \\
&=&B^{2k-4+3k\ep}(R^{k+2}B^{4}{\phi^*}^{k+2}+R^{k}B{\phi^*}^{k-1}). \label{RBR}
\end{eqnarray}
Here we use the crude bound $RL\ll B^2$ to show that $(RL)^{\ep}\ll
B^{2\ep}$. On using first the assumption that $\phi^*\ll
R^{-1}B^{-2/3}$, and then that $R\ll B^{4/3}$, we see that
(\ref{RBR}) is
\beq\label{minor_term1}
\ll B^{2k-4+3k\ep}(B^{-(2k-8)/3}+RB^{-(2k-5)/3})\ll B^{2k-4+3k\ep}.B^{-(2k-9)/3}.
\eeq
Since $k\ge 5$ this is satisfactory for Proposition \ref{prop_sig_bd}, after re-defining $\ep$.

Finally, if $B^2\phi^*\le 1$, then $L\ll RB^{-1+\ep}$, so that
(\ref{REST1}) becomes
\[\ll B^{2k-4}(\log B)
R^{-2k}(R^{k+2+\ep}(RB^{-1+\ep})^{2k+\ep}+R^{k+3+\ep}(RB^{-1+\ep})^{2k-3+\ep}).\]
We may simplify this if $\ep$ is small enough to get
\beq\label{minor_term2}
\ll B^{2k-4+3k\ep}(R^{k+2}B^{-2k}+R^{k}B^{3-2k}) \ll B^{2k-4+3k\ep}.B^{-(2k-9)/3},
\eeq
since $R\ll B^{4/3}$. For $k\ge 5$ this is also satisfactory for Proposition \ref{prop_sig_bd}, after re-defining $\ep$.

We remark here that it is now visible that the most significant terms
(\ref{minor_term1}) and (\ref{minor_term2}) contributing to the minor
arcs involve a saving $O(B^{-(2k-9)/3})$, so that the result is
non-trivial for $k>9/2$.  Thus to handle $k=4$ one would have to do
more than shave off a small power of $B$.

\section{Proof of Theorems \ref{thm_r2} and \ref{thm_s2}}
\label{sec_sing_int}

We first prove Theorem \ref{thm_r2}. The statements about
$\Sf(\u{n})$ follow from Propositions \ref{prop_sing_ser_conv} and \ref{prop_sing_ser_conv2}, while the
uniform boundedness of $\Jcal_w(\u{\mu})$ is part of
Proposition \ref{prop_sing_int_conv}.  Thus it remains to consider
lower bounds for the singular integral $\Jcal_w(\u{\mu})$.

We begin by establishing the following result.
\begin{lemma}\label{lemma_small_x}
Let $A_1>A_2>0$ be given.  Then there exists $\Lambda>0$, dependent
only on $A_1,A_2$ and $\u{Q}$, such
that if $A_2\le |\u{\nu}|\le A_1$ and if $\u{Q}(\xbf)=\u{\nu}$ has a solution
 $\xbf\in\R^k$,  then in fact there exists a real solution of $\u{Q}(\xbf) = \u{\nu}$ satisfying $|\xbf| \leq \Lambda$.
\end{lemma}

We may clearly reduce to the case $|\u{\nu}|=1$ by rescaling.
To find a suitable value for $\Lambda$ we consider two cases. First suppose
that $\u{Q}(\xbf)\not=\u{0}$ for all $\xbf$ with $|\xbf|=1$. Then by continuity and compactness, we deduce that $\inf_{|\xbf| =1} |\u{Q}(\xbf)|>0$.  Thus for any $\xbf$ such that $\uQ(\xbf) = \u{\nu}$ we have
\[ 1 = |\u{\nu}| = |Q(\xbf)| \geq |\xbf|^2\inf_{|\xbf| =1} |\u{Q}(\xbf)|.\]
It follows that if we define
\[\Lambda=\left(\inf_{|\xbf| =1} |\u{Q}(\xbf)|\right)^{-1/2}\]
then $|\xbf|\leq\Lambda$.

In the alternative case there exists $\b{a}$ with $|\b{a}|=1$ such
that $\u{Q}(\b{a}) = \u{0}$. Then by Condition 1 we have $\rk( J(\b{a})) = 2$.  We shall suppose that
\[\left| \begin{array}{cc}
 \frac{\partial Q_1(\b{a})}{ \partial x_1} &
                       \frac{\partial Q_1(\b{a})}{ \partial x_2}\\
			\frac{\partial Q_2(\b{a})}{ \partial x_1} &
                       \frac{\partial Q_2(\b{a})}{ \partial x_2}
			\end{array} \right|\not=0,\]
as we may, without any loss of generality.
We can therefore apply the Implicit Function Theorem to the mapping $F:\R^4\rightarrow\R^2$ given by
\[F(t_1,t_2,\mu_1,\mu_2)=\u{Q}(t_1+a_1,t_2+a_2,a_3,\ldots,a_n)-\u{\mu}.\]
Since $F(\b{0},\b{0})=\u{0}$ we deduce that there is a $\delta>0$ and
a continuous function
\[G:\{\u{\mu}: |\u{\mu}|\le\delta\}\rightarrow\R^2\]
such that $G(\b{0})=\b{0}$ and $F(G(\u{\mu}),\u{\mu})=\b{0}$.
Since $G$ is continuous it is bounded for $|\u{\mu}|\le\delta$,
by $\kappa$, say.  Thus if $|\u{\mu}|\le\delta$ there will be a solution \[\xbf=(t_1+a_1,t_2+a_2,a_3,\ldots,a_n)\]
of $\u{Q}(\xbf)=\u{\mu}$
satisfying $|\xbf|\le 1+\kappa$. Hence, given $|\u{\nu}| = 1$, we take $\del>0$ as found above and set $\u{\mu} = \del \u{\nu}$ so that there exists $\xbf$ with $|\xbf| \leq 1 + \kappa$ such that
$\u{Q}(\xbf) = \u{\mu} = \del \u{\nu}$. Then by setting $\ybf = \del^{-1/2} \xbf$ we provide a solution to $\u{Q}(\ybf) = \u{\nu}$ with $|\ybf| \leq \del^{-1/2} (1+\kappa)$.
This establishes the lemma in the second case, with
$\Lambda=\delta^{-1/2}(1+\kappa)$.

Our next result is a real analogue of Lemma \ref{4rcontrol}.
\begin{lemma}\label{ana}
Let $A_1>A_2>0$ be given. Let $\xbf\in\R^k$ with $A_2\le|\xbf|\le
A_1$, and suppose that
$\u{Q}(\xbf)=\u{\nu}$. Then \[F(\nu_2,-\nu_1)\ll \max_{i,j}|\Delta_{ij}(\xbf)|^2.\]
\end{lemma}
Again we may reduce to the case $|\u{\nu}|=1$ by rescaling.
Since the proof of the lemma
is completely analogous to that of Lemma \ref{4rcontrol} we leave the details to the reader.

To complete the proof of Theorem \ref{thm_r2} we require one further
lemma.
\begin{lemma}\label{last}
Let $A_1>A_2>0$ be given.  Then there is a constant $\kappa_0=\kappa_0(A_1,A_2,Q_1,Q_2)$ with $0<\kappa_0
 \leq 1$ such that for any $\kappa\leq \kappa_0$ the following
 holds. Let $\b{a}\in\R^k$ with $A_2\le |\b{a}|\le A_1$.  Write \[M=\max_{i,j}|\Delta_{ij}(\b{a})|\]
and assume that $M>0$. Re-order the indices $i$ so that $|\Delta_{12}(\b{a})|=M$.  Then for any $x_3,\ldots,x_k$ with \beq\label{t23}
\max_{3\le i\le k}|x_i|\le (\kappa M)^2
\eeq
there exist $x_1$ and $x_2$ in the square $S$ given by
\beq\label{t24}
S=\{(x_1,x_2)\in\R^2:\,\max(|x_1|,|x_2|)\le \kappa M\}
\eeq
such that
\[\u{Q}(\b{a}+\xbf)=\u{Q}(\b{a}).\]
\end{lemma}

It will suffice to show that the conclusion of the lemma holds
whenever $\kappa>0$ is sufficiently small in terms of $A_1, A_2, Q_1$ and $Q_2$. For the proof it will be convenient to write
\[\Mcal(\xbf)=\left(\begin{array}{cc}
 \frac{\partial Q_1(\b{x})}{ \partial x_1} &
                       \frac{\partial Q_1(\b{x})}{ \partial x_2}\\
			\frac{\partial Q_2(\b{x})}{ \partial x_1} &
                       \frac{\partial Q_2(\b{x})}{ \partial x_2}
			\end{array} \right).\]
Then $||\Mcal(\xbf)||\ll 1$ if $|\xbf|\ll 1$, and \beq\label{t21}
|\det(\Mcal(\b{x}))|\ge M/2
\eeq
if $|\xbf-\b{a}|\ll M$ with a small enough implied constant.  It follows that
\beq\label{t21'}
|| \Mcal(\xbf)^{-1}||\ll M^{-1}
\eeq
when $|\xbf-\b{a}|\ll M$.

We now write each $\xbf\in\R^k$ in the shape $(\b{u},\b{v})$ where $\b{u}$ corresponds to the first two
variables $x_1,x_2$ and $\b{v}$ to the remaining variables
$x_3,\ldots,x_k$. It will also be convenient to write
$\b{a}=(\b{b},\b{c})$ accordingly.
For each vector $\b{v}\in\R^{k-2}$ we now consider the function
$F_{\b{v}}:\R^2\rightarrow\R^2$ given by
\[F_{\b{v}}(\b{u})=
\b{u}-\Mcal(\b{b},\b{v})^{-1}\{\u{Q}(\b{b}+\b{u},\b{v})-\u{Q}(\b{b},\b{c})\}.\]
It follows from the definition that
\beq\label{t22}
F_{\b{v}}(\b{u})=-\Mcal(\b{b},\b{v})^{-1}
\{\u{Q}(\b{b},\b{v})+\u{Q}(\b{u},\b{0})-\u{Q}(\b{b},\b{c})\},
\eeq
whence (\ref{t21'}) yields
\[F_{\b{v}}(\b{u})\ll M^{-1}(|\b{u}|^2+|\b{v}-\b{c}|)\]
if $|\b{v}-\b{c}|\ll M$ with a sufficiently small implied constant.
Hence choosing $\kappa \in(0,1)$ sufficiently small, if $\b{u}=(x_1,x_2)$ lies in
$S$ and if $\b{v}-\b{c}=(x_3,\ldots,x_k)$ satisfies
(\ref{t23}), then
$F_{\b{v}}$ maps the square (\ref{t24}) to itself.
\label{note_F_square}

Moreover if we have
two vectors $\b{u}^{(1)},\b{u}^{(2)}$ in $S$, then
\[F_{\b{v}}(\b{u}^{(1)})-F_{\b{v}}(\b{u}^{(2)})=
\Mcal(\b{b},\b{v})^{-1}\{\u{Q}(\b{u}^{(2)},\b{0})-\u{Q}(\b{u}^{(1)},\b{0})\},\]
and
\[Q_i(\b{u}^{(2)},\b{0})-Q_i(\b{u}^{(1)},\b{0})\ll |\b{u}^{(1)}-\b{u}^{(2)}|\max(|\b{u}^{(1)}|,|\b{u}^{(2)}|)\ll \kappa M
|\b{u}^{(1)}-\b{u}^{(2)}|.\]
It follows that (\ref{t21'}) yields
\[|F_{\b{v}}(\b{u}^{(1)})-F_{\b{v}}(\b{u}^{(2)})|
\ll \kappa |\b{u}^{(1)}-\b{u}^{(2)}|\]
if $|\b{v}-\b{c}|\ll M$.  We therefore conclude
that if $\kappa$ is small enough then the function
$F_{\b{v}}$ is a contraction mapping on $S$
whenever (\ref{t23}) holds. It follows that $F_{\b{v}}$ has a fixed point $\b{u}\in S$,
which means that
\[\u{Q}(\b{b}+\b{u},\b{v})=\u{Q}(\b{b},\b{c})\]
by construction of $F_{\b{v}}$.  The lemma now follows.

\label{note_box_text}
We are finally in a position to complete the proof of Theorem \ref{thm_r2}.
Suppose that $1/2\le \max(|\mu_1|,|\mu_2|)\le 1$, whence
$1/\sqrt{2}\le |\u{\mu}|\le\sqrt{2}$.  Suppose further that
$Q_1(\xbf)=\mu_1$, $Q_2(\xbf)=\mu_2$ has a solution $\xbf=\b{a}\in\R^k$.
According to Lemma \ref{lemma_small_x} we may assume that
$|\b{a}|\le\Lambda$ for some $\Lambda$ depending only on $Q_1$ and
$Q_2$. Moreover, since $\max(|Q_1(\b{a})|,|Q_2(\b{a})|)\ge 1/2$ we
deduce that $|\b{a}|\gg 1$.
We shall take $C=2+\Lambda$ in Theorem \ref{thm_r2}, with $\Lambda$ as above.  Then if $w(\xbf)>0$ for all $\xbf$ with $\max |x_i|\le C$ we may use compactness to show that there is a
constant $c_0>0$ such that $w(\xbf)\ge c_0$ for all such $\xbf$.

We now write $M=\max |\Del_{ij} (\abf) |$ as in Lemma \ref{last}, whence
Lemma \ref{ana} shows that $|F(\mu_2,-\mu_1)|\le c_0 M^2$ for some
$c_0$ depending only on $Q_1$ and $Q_2$.  Now take
$\kappa_0=\kappa_0(\sqrt{2},1/\sqrt{2},Q_1,Q_2)$ as in Lemma
\ref{last}.  Then if $0<\kappa\le\kappa_0$ and
\beq\label{t25}
|(x_3,\ldots,x_k)|\le \kappa^2c_0^{-1}|F(\mu_2,-\mu_1)|
\eeq
we will have $|(x_3,\ldots,x_k)|\leq (\kappa M)^2$, whence Lemma
\ref{last} will produce
values of $x_1,x_2$  with $\max(|x_1|,|x_2|) \leq \kappa M$ such that $\u{Q}(\b{a}+\xbf)=\u{\mu}$.  By taking
$\kappa$ sufficiently small we may ensure that $|\xbf|\le 1$.

We now use the same notation $\b{a}=(\b{b},\b{c})$ as before, and set
$\b{u}_0=(x_1,x_2)$ and $\b{v}=(x_3,\ldots,x_k)$, so that for any
$\b{v}$ satisfying (\ref{t25}) there is a corresponding $\b{u}_0$ such
that $\u{Q}(\b{u}_0,\b{v})=\u{\mu}$.  We proceed to consider values of
$\u{Q}$ near to $\u{\mu}$.  If $c_1>0$ is small enough, then for any $\b{v}=(x_3,\ldots,x_k)$ satisfying (\ref{t25}), and any
$\ep\in(0,1)$, there is
a square \[ \{ \ubf = (u_1,u_2) : |u_1-x_1|\le c_1 \ep,\;\;\;|u_2-x_2|\le c_1 \ep \} \]
on which $|\u{Q}(\b{b}+\b{u},\b{c}+\b{v})-\u{\mu}|\le \ep/2$. Moreover,
if $c_1$ is small enough then we will have $|(\bbf+ \ubf, \cbf+ \vbf)|
\leq \Lambda + 2$ for all such $(\ubf,\vbf)$.  We
therefore see that
\begin{eqnarray*}
\lefteqn{\int_{\substack{|\xbf|\le \Lambda+2\\\max|Q_i(\xbf)-\mu_i|\le\ep}}
\left(1-\frac{|Q_1(\xbf)-\mu_1|}{\ep}\right)
\left(1-\frac{|Q_2(\xbf)-\mu_2|}{\ep}\right)d\xbf}\hspace{1cm}\\
&\gg&  \ep^2|F(\mu_2,-\mu_1)|^{k-2},\hspace{6cm}
\end{eqnarray*}
whence
\begin{eqnarray*}
\lefteqn{\int_{\max|Q_i(\xbf)-\mu_i|\le\ep}w(\xbf)
\left(1-\frac{|Q_1(\xbf)-\mu_1|}{\ep}\right)
\left(1-\frac{|Q_2(\xbf)-\mu_2|}{\ep}\right)d\xbf}\hspace{1cm}\\
&\gg& \ep^2|F(\mu_2,-\mu_1)|^{k-2}.\hspace{6cm}.
\end{eqnarray*}
The claimed lower bound for $\Jcal_w(\u{\mu})$ then follows from
Proposition \ref{prop_sing_int_val}.

We turn now to the proof of Theorem \ref{thm_s2}, which follows similar lines.  The positivity of the
singular series follows immediately from Propositions \ref{prop_sing_ser_conv0} and \ref{prop_sing_ser_conv20}, so that it
remains to show that $\Jcal_w(\u{0})>0$. We first do this under the
assumption that $w(\x)$ is supported on the hypercube $\max|x_i|\le 3$.
We are supposing also that
$\u{Q}(\b{a})=\u{0}$ for some non-zero real vector $\b{a}$, and by
homogeneity we may take $|\b{a}|=1$. Such a solution has
$\rk(J(\b{a}))=2$ by Condition 1, whence
$\max_{i,j}|\Del_{i,j}(\b{a})|>0$. We may therefore complete the proof that
$\Jcal_w(\u{0})>0$ using Lemma \ref{last}, just as we did for
Theorem \ref{thm_r2}, but with each occurence of the numbers $\Lambda$
and  $F(\mu_2,-\mu_1)$ replaced by the value 1.  Finally we note that
we can rescale the weight $w(\x)$ and the parameter $B$ without
affecting the conclusion that $\Jcal_w(\u{0})>0$.  Since our theorem
assumes that the support of $w(\x)$ includes a small hypercube around
the origin this allows us to suppose that in fact the support includes
the set $\max|x_i|\le 3$. This observation completes the proof.

\section{Proof of Theorems \ref{thm_gen} and \ref{thm_prime}}
\label{sec_final}

It is easy to see, using a dyadic subdivision, that in proving Theorem
\ref{thm_gen} it will suffice to handle pairs of integers $n_1,n_2$
with $N/2\le\max(|n_1|,|n_2|)\le N$. We choose a weight $w$ such that
$w(\xbf)>0$ for $|\xbf|\le C$, with $C$ as in Theorem \ref{thm_r2}, and
we take $B=N^{1/2}$. We classify
pairs $\u{n}$ contributing to $\mathcal{E}(N)$ into three cases,
which may overlap. Case I will consist of pairs for which
$\Jcal_w(B^{-2}\u{n})\ge B^{-(k-2)/(16k^2)}$ and $\Sf(\u{n})\ge
B^{-(k-2)/(32k)}$.  Case II will be that in which $\Jcal_w(B^{-2}\u{n})\le
B^{-(k-2)/(16k^2)}$ or $F(n_2,-n_1)=0$, while Case III will have
$\Sf(\u{n})\le B^{-(k-2)/(32k)}$ and $F(n_2,-n_1)\not=0$.

The number of pairs $\u{n}$ for which $F(n_2,-n_1)=0$ is clearly
$O(N)$. For the remaining pairs, if there is no integer
solution $\xbf$ with $\u{Q}(\xbf)=\u{n}$ then $R_B(\u{n})=0$. For
pairs $\u{n}$ in Case I the
corresponding summand in (\ref{tt}) is then \[\gg (B^{-(k-2)/(16k^2)}.B^{-(k-2)/(32k)}.B^{k-4})^2.\]
It follows from Theorem \ref{thm_r1} that there are \[\ll B^{4-1/(4k^2)}\ll N^{2-1/(2k^2)}\]
pairs $(n_1,n_2)$ for which this first
case holds.  This is satisfactory for Theorem \ref{thm_gen}.

We have already treated those $\u{n}$ for which $F(n_2,-n_1)=0$, so we
turn to Case II with the assumption that $\Jcal_w(B^{-2}\u{n}) \leq
B^{-(k-2)/(16k^2)}$ but $F(n_2,-n_1)\neq0$. Since we are assuming that
the system of equations
$Q_1(\xbf)=n_1$, $Q_2(\xbf)=n_2$ has a real solution $\xbf$, we may
deduce from Theorem \ref{thm_r2} with $\u{\mu}=B^{-2}\u{n}$ that
$|F(\mu_2,-\mu_1)|\ll B^{-1/(16k^2)}$, whence
$|F(n_2,-n_1)|\ll N^kB^{-1/(16k^2)}$. We factor $F(x_1,x_2)$
over $\C$ as in (\ref{dp}), and write $\psi=an_2+bn_1$ for the
smallest factor $\psi_i$ of $F(n_2,-n_1)$. Then, as in the proof of
(\ref{Flb}), we see that $|F(n_2,-n_1)|\gg N^{k-1}|\psi|$, since
$\max(|n_1|,|n_2|)\ge N/2$.  It follows that
$|\psi|\ll NB^{-1/(16k^2)}$ for some factor $\psi$ on the
right of (\ref{dp}).  Assuming that the coefficient $a$, say, is
non-zero we deduce that for each $n_1$ the value for $n_2$ is
restricted to an interval of length $O(NB^{-1/(16k^2)})$.  We
then deduce that the number of possible pairs $\u{n}$ in Case II is
$O(N^2B^{-1/(16k^2)})$, which is satisfactory for Theorem
\ref{thm_gen}.

We turn finally to Case III.  Since we are assuming that
the system of equations
$Q_1(\xbf)=n_1$, $Q_2(\xbf)=n_2$ is solvable in every $p$-adic ring
$\Z_p$, we may deduce from Theorem \ref{thm_r2} that if $(n_1,n_2)$ belongs to Case III then
\[\prod_{p\le p_0}|F(n_2,-n_1)|_p\ll_{\ep}B^{\ep-1/(32k)}\]
for any fixed $\ep>0$. Let $P=\prod_{p\le p_0}p$.  We then see that
there is a divisor $q\gg_{\ep}B^{-\ep+1/(32k)}$ of $F(n_2,-n_1)$ such
that $q|P^{\infty}$. Since $F(n_2,-n_1)$ is non-zero and $\ll B^{2k}$, we conclude that
$q\ll B^{2k}$.  Thus the exponent to which any given prime divides $q$
will be $O(\log B)$.  It follows that $q|P^{E}$ for some $E\ll\log B$,
so that the number of possibilities for $q$ is $O((\log B)^{p_0})$.
Since $p_0$ depends only on $Q_1$ and $Q_2$ we deduce that the number
of possibilities for $q$ is $O_{\ep}(B^{\ep})$, following our
convention that the implied constant may depend also on the forms
$Q_1$ and $Q_2$.  We therefore conclude that the number of pairs $\u{n}$ in
Case III is controlled by up to $B^\ep$ factors $q$ of size
$q\gg_{\ep}B^{-\ep+1/(32k)}$, where for each such $q$, the
corresponding number of $\u{n}$ belonging to Case III is \[\ll_{\ep}\#\{\u{n}\in\Z^2:\, |n_1|,|n_2|\le N,\,
q|F(n_2,-n_1)\}.\]
We shall estimate this for each $q$ very crudely.
For any fixed $n_1\not=0$ the polynomial $F(x,n_1)$ in $x$ has discriminant
$n_1^{k(k-1)}D_F$, where $D_F\not=0$ is the discriminant of
$F(x_1,x_2)$. According to Huxley \cite{Hux79} the congruence $F(X,n_1)\equiv 0\modd{q}$ then has at most $k^{\omega(q)}|D_F|^{1/2}a^{k(k-1)/2}$ roots modulo $q$, where
$a=a(n_1,q)$ is the largest factor of $n_1$ which divides $q^{\infty}$. If we now fix $A>0$ we see that $F(X,n_1)\equiv 0\modd{q}$ has
$O_{\ep}(B^{\ep}A^{k(k-1)/2})$ roots when $a\le A$, so that the
range $|X|\le N$ produces $O_{\ep}(B^{\ep}A^{k(k-1)/2}(N/q+1))$
solutions.  The contribution to Case III from the $\ll B^\ep$
appropriate $q$ and non-zero integers $n_1$ with $|n_1| \leq N$ such that
$a(n_1,q)\le A$ is therefore
\[\ll_{\ep}B^{2\ep}A^{k(k-1)/2}(N/q+1)N \ll_{\ep}
B^{3\ep-1/(32k)}A^{k(k-1)/2}N^2,\]
since $q\gg_{\ep}B^{-\ep+1/(32k)}$. In the remaining case either
$n_1=0$ or $r|n_1$ for some integer $r|q^{\infty}$ satisfying $A<r\le
N$.  The number of possibilities for $r$ is $O_{\ep}(B^{\ep})$, by the
same argument that bounded the number of possibilities for $q$.  Thus
there are $O_{\ep}(1+B^{\ep}N/A)$ choices for $n_1$, with a total
contribution $O_{\ep}((1+B^{\ep}N/A)N)$ for Case III. We will choose
$A=B^{\phi}$ with $\phi=(16k(2+k(k-1)))^{-1}$, in which case we may verify that $B^{3\ep-1/(32k)}A^{k(k-1)/2}\ll B^{-1/(16k^3)}.$
This shows that Case III makes a satisfactory contribution
of $O_{\ep}(N^2B^{-1/(16k^3)})$ to Theorem \ref{thm_gen}. This
completes the proof of the theorem.

We turn now to the proof of Theorem \ref{thm_prime}. By (\ref{x1}) and
(\ref{x2}) we see that there is a constant $p_1$ depending on $Q_1$
and $Q_2$ alone, such that $\sigma_p>0$ for all good primes $p\ge
p_1$.  We may of course choose $p_1$ sufficiently large that we have
$p<p_1$ for all bad primes. Whenever $\sigma_p>0$ we see from
(\ref{sig_N}) that $N(\u{n};p^e)>0$ for all large enough $e$, from
which a compactness argument shows that there is at least one solution
of $\u{Q}(\xbf)=\u{n}$ with $\xbf\in\Z_p^k$.  Thus this local
condition holds for all $\u{n}\in\Z^2$ as soon as $p\ge p_1$.
\label{note_compactp}

For each of the finitely many primes $p<p_1$ we will show that there
is an exponent $e=e(p)$ and a congruence class $\u{n}^{(p)}\modd{p^e}$ with $p\nmid n_1^{(p)}n_2^{(p)}$, such that $\u{Q}(\xbf)=\u{n}$ has a
solution $\xbf\in\Z_p^k$ whenever $\u{n}\equiv\u{n}^{(p)}\modd{p^e}$.
To do this we consider vectors $\ybf_p$ of the form $\xbf_p+p\b{m}$
where $\xbf_p$ is as in the statement of Theorem \ref{thm_prime} and
$\b{m}$ runs over $\Z^k$.  Since the determinants $\Delta_{ij}(\xbf)$
do not all vanish identically we can find an integer vector $\b{m}$
for which some $\Delta_{ij}(\ybf_p)$ is non-zero.  We then set
$\u{n}^{(p)}=\u{Q}(\ybf_p)$ so that \[n_i^{(p)}=Q_i(\ybf_p)\equiv Q_i(\xbf_p)\not\equiv 0\modd{p}\]
for $i=1,2$.  Suppose now that $e=2f+1$ where
$p^f||\Delta_{ij}(\ybf_p)$ and let
$\u{n}\equiv\u{n}^{(p)}\modd{p^e}$.  Then
$\u{Q}(\xbf)\equiv\u{n}\modd{p^e}$ has a solution $\xbf=\ybf_p$
which can be lifted to $\Z_p$ by Hensel's Lemma since $e\ge 2f+1$.
This establishes our claim.

For the real valuation we can produce a completely analogous argument.
There is a neighbourhood of $\xbf_0$ on which the forms $\u{Q}$ are
both positive, and this neighbourhood will include a point $\ybf_0$ at
which some determinant $\Delta_{ij}(\ybf_0)$ is non-vanishing. There
is then a small $\delta>0$ such that the system $\u{Q}(\xbf)=\u{\mu}$
has a real solution whenever $|\u{\mu}-\u{Q}(\ybf_0)|\le \delta$.

We can now use the Chinese Remainder Theorem to produce a modulus
$M=\prod_{p<p_1}p^{e(p)}$ and a residue class $\u{n}^{(M)}\modd{M}$
such that $n_1^{(M)}$ and $n_2^{(M)}$ are both coprime to $M$, and
with the property that if $\u{n}\equiv\u{n}^{(M)}\modd{M}$ then
$\u{Q}(\xbf)=\u{n}$ has a solution in every ring $\Z_p$ for $p<p_1$.
Of course there is also a solution for $p\ge p_1$, by our choice of
$p_1$. It follows that if $r_1$ and $r_2$ are primes such that $r_i
\equiv n_i^{(M)}\modd{M}$ for $i=1,2$ and such that \[|r_i-KQ_i(\ybf_0)|\le \delta K/2,\;\;\; (i=1,2)\] for some dilation factor $K>0$, then the equations $\u{Q}(\xbf)=(r_1,r_2)$ have a solution in $\R$ and in every ring $\Z_p$.  The number of such pairs of primes is
$\gg K^2(\log K)^{-2}$ by the Prime Number Theorem for arithmetic
progressions, but, according to Theorem \ref{thm_gen} at most
$O(K^{2-\varpi})$ pairs can fail to have a representation over $\Z$.
Hence there is at least one representable pair of primes $(r_1,r_2)$ for each
sufficiently large $K$.  This completes the proof of the theorem.

We conclude with an explanation of the conjecture that the result of
Theorem \ref{thm_prime} should continue to hold if the assumption of
Condition \ref{cond1} is replaced by the assumption that neither
$Q_1(\xbf)$ nor $Q_2(\xbf)$ factors over $\Z$. By one of the assumptions of
Theorem \ref{thm_prime}, there exists $\abf \in \Z^k$ such that
$Q_1(\abf)$ and $Q_2(\abf)$ are both positive.  Set
$A:=Q_1(\abf)Q_2(\abf)$ and $M=\prod_{p|A} p$.  By a further
hypothesis of Theorem \ref{thm_prime}, for each $p$ there exists
$\xbf_p \in \Z^k$ such that $p \ndiv Q_1(\xbf_p)Q_2(\xbf_p)$. Apply
the Chinese Remainder Theorem to construct a residue class $\mbf^{(M)}
\modd{M}$ such that $\mbf^{(M)} \con \xbf_p \modd{p}$ for each $p|A$. Fix any $\bbf \in \Z^k$ with $\bbf \con \mbf^{(M)}\modd{M}$ and $\bbf \neq \abf$ (of which
there are infinitely many choices). Define for this choice of $\abf,\bbf$ a pair of quadratic polynomials
in a real variable $t$ given by \[ q_i(t) = Q_i(\bbf + t(\abf - \bbf)),\]
for $i=1,2$.
Then $q_i(t) \maps \infty$ as $t \maps \infty$. By construction, if $p \ndiv A$, then $p \ndiv q_i(t)$ when $t=1$, and
if $p|A$ then $p \ndiv q_i(t)$ when $t=0$. Thus the polynomials
$q_1(t),q_2(t)$ have no fixed prime divisor. If $q_1,q_2$ are irreducible over $\Z$ as polynomials in $t$, then
Schinzel's Hypothesis would imply that $q_1(t)$, $q_2(t)$
simultaneously attain infinitely many prime values, and hence so would
$Q_1,Q_2$. Thus the remaining consideration is to show that if
$Q_1,Q_2$ are irreducible over $\Q$, then there is a choice of $\bbf\con \mbf^{(M)}\modd{M}$ such that $q_1(t)$ and $q_2(t)$ are both irreducible.

The polynomial $q_i(t)$ is reducible if and only if its discriminant
is a square.  However the discriminant will be (up to a factor of 4)
\beq\label{Qdisc}
Q_i(\bbf,\abf)^2 - Q_i(\bbf)Q_i(\abf),
\eeq
where $Q_i(\xbf,\ybf) = \xbf^t Q_i \ybf$ is the associated bilinear
form. As a function of $\bbf$, (\ref{Qdisc}) is a quadratic form, $R_i(\bbf)$,
say.  Moreover, if $R_i$ were the square $L_i(\xbf)^2$ of a linear
form $L_i$ defined over $\Q$, then we would have
\[Q_i(\xbf)=
Q_i(\abf)^{-1}\{Q_i(\xbf,\abf) + L_i(\xbf)\}\{Q_i(\xbf,\abf) - L_i(\xbf)\},\]
in contradiction to the assumption that $Q_i(\xbf)$ is irreducible
over $\Q$.  Thus neither $R_1(\xbf)$ nor $R_2(\xbf)$ can be squares.

It follows that each of the quadratic forms
\[S_1(\xbf,y)=R_1(\xbf)-y^2\;\;\;\mbox{and}\;\;\;
S_2(\xbf,y)=R_2(\xbf)-y^2\]
is irreducible over $\Q$. For any $\bbf$ such that $q_1(t)$ is
reducible, there will be a corresponding integer $y$ such that
$S_1(\bbf,y)=0$.  Moreover if $|\bbf|\le B $, say, then $y\ll B$.
Since an irreducible quadratic form $S$ in $n$ variables has $O_{S,\varepsilon}(B^{n-2+\varepsilon})$ integral zeros of size $O(B)$
we deduce, on taking $\varepsilon=1/2$, that there are $O(B^{k-1/2})$
admissible values of $\bbf$ with $|\bbf|\le B$, such that $q_1(t)$ is reducible.

There is a similar estimate for $q_2(t)$.  However there are $\gg B^k$
vectors $\bbf\con \mbf^{(M)}\modd{M}$ such that $|\bbf|\le B$, whence
if $B$ is large enough there must be some such value for which both
$q_1(t)$ and $q_2(t)$ are irreducible.  This provides the final step
in our argument. We conclude with the observation that the above
argument is essentially a proof of a case of the Hilbert
irreducibility theorem.

\section*{Acknowledgements}

Pierce was partially supported during this work by a Marie Curie
Fellowship funded by the European Commission and NSF DMS-0902658. Gratitude is extended to Damaris Schindler for comments on the manuscript, and to the Max-Planck-Institut f\"{u}r Bioanorganische Chemie, which graciously provided a very pleasant working
environment to the second author during the preparation of an early
version of this manuscript.

\bibliographystyle{amsplain}
\bibliography{NoThBibliography}

\end{document}

%% file: format.tex
\newtheorem{thm}{Theorem}[section]
\newtheorem{prop}{Proposition}[section]
\newtheorem{lemma}{Lemma}[section]
\newtheorem{cor}{Corollary}[thm]

\newtheorem{condition}{Condition}

\newtheorem{dfn}{Definition}

\newcommand\nonumfootnote[1]{%
  \begingroup
  \renewcommand\thefootnote{}\footnote{#1}%
  \addtocounter{footnote}{-1}%
  \endgroup
}

\renewcommand{\b}[1]{\mathbf{#1}}
\renewcommand{\u}[1]{\underline{#1}}

\newcommand{\twosum}[2]{\sum_{\substack{#1\\#2}}}
\newcommand{\ep}{\varepsilon}

\newcommand{\con}{\equiv}

\newcommand{\ndiv}{\nmid}
\newcommand{\modd}[1]{\; ( \text{mod} \; #1)}

\newcommand{\maps}{\rightarrow}
\newcommand{\intersect}{\cap}

\newcommand{\Union}{\bigcup}

\newcommand{\diag}{\text{diag}\,}

\newcommand{\disc}{{\rm disc \;}}

\newcommand{\rk}{{\rm rk }}

\newcommand{\al}{\alpha}
\newcommand{\be}{\beta}

\newcommand{\del}{\delta}
\newcommand{\Del}{\Delta}

\newcommand{\Sig}{\Sigma}
\newcommand{\sig}{\sigma}
\newcommand{\lam}{\lambda}

\newcommand{\Ecal}{\mathcal{E}}

\newcommand{\Ical}{\mathcal{I}}
\newcommand{\Jcal}{\mathcal{J}}

\newcommand{\Mcal}{\mathcal{M}}

\newcommand{\Qcal}{\mathcal{Q}}
\newcommand{\Scal}{\mathcal{S}}

\newcommand{\abf}{{\bf a}}
\newcommand{\bbf}{{\bf b}}
\newcommand{\cbf}{{\bf c}}

\newcommand{\ebf}{{\bf{e}}}

\newcommand{\lbf}{{\bf l}}
\newcommand{\mbf}{{\bf m}}

\newcommand{\rbf}{{\bf r}}

\newcommand{\ubf}{{\bf u}}
\newcommand{\vbf}{{\bf v}}
\newcommand{\wbf}{{\bf{w}}}
\newcommand{\x}{{\bf x}}
\newcommand{\xbf}{{\bf x}}
\newcommand{\ybf}{{\bf y}}
\newcommand{\zbf}{{\bf z}}

\newcommand{\lambf}{\boldsymbol\lambda}

\newcommand{\unu}{\u{\nu}}
\newcommand{\uQ}{\u{Q}}

\newcommand{\Qbar}{\overline{\Q}}

\newcommand{\C}{\mathbb{C}}
\newcommand{\F}{\mathbb{F}}

\newcommand{\N}{\mathbb{N}}
\newcommand{\Proj}{\mathbb{P}}

\newcommand{\Q}{\mathbb{Q}}
\newcommand{\R}{\mathbb{R}}

\newcommand{\Z}{\mathbb{Z}}

\newcommand{\Mf}{\mathfrak{M}}

\newcommand{\mf}{\mathfrak{m}}
\newcommand{\Sf}{\mathfrak{S}}

\newcommand{\rank}{\text{rank}\,}

\newcommand{\beq}{\begin{equation}}
\newcommand{\eeq}{\end{equation}}

%% file: QuadForms_V13.bbl
\providecommand{\bysame}{\leavevmode\hbox to3em{\hrulefill}\thinspace}
\providecommand{\MR}{\relax\ifhmode\unskip\space\fi MR }
\providecommand{\MRhref}[2]{%
  \href{http://www.ams.org/mathscinet-getitem?mr=#1}{#2}
}
\providecommand{\href}[2]{#2}
\begin{thebibliography}{10}

\bibitem{Bir61}
B.~J. Birch, \emph{Forms in many variables}, Proc. Roy. Soc. Ser. A
  \textbf{265} (1961/62), 245--263.

\bibitem{BHB05}
T.~D. Browning and D.~R. Heath-Brown, \emph{Counting rational points on
  hypersurfaces}, J. Reine Angew. Math. \textbf{584} (20045), 83--115.

\bibitem{BrMu13}
T.~D. Browning and R.~Munshi, \emph{Rational points on singular intersections
  of quadrics}, Compositio Math., in press.

\bibitem{CTSSD}
J.~L. Colliot-Th\'el\`ene, J.~J. Sansuc, and P.~Swinnerton-Dyer,
  \emph{Intersections of two quadrics and châtelet surfaces. i}, J. Reine
  Angew. Math. \textbf{373} (1987), 37--107.

\bibitem{FI}
E.~Fouvry and H.~Iwaniec, \emph{Gaussian primes}, Acta Arith. \textbf{79}
  (1997), 249--287.

\bibitem{HL}
G.~H. Hardy and J.~E. Littlewood, \emph{Some problems of ``{Partitio}
  {Numerorum}'' {VI}: {Further} researches in {Waring's} problem}, Math.
  Zeitschrift \textbf{23} (1925), no.~1, 1--37.

\bibitem{HB96}
D.~R. Heath-Brown, \emph{A new form of the circle method, and its application
  to quadratic forms}, J. Reine Angew. Math. \textbf{481} (1996), 149--206.

\bibitem{Hoo}
C.~Hooley, \emph{On {H}ypothesis ${K}^*$ in {W}aring's problem}, Sieve Methods,
  Exponential Sums, and their Applications in Number Theory, London Math. Soc.
  Lecture Notes No. 237, Cambridge University Press, 1997, pp.~175--185.

\bibitem{Hux79}
M.~N. Huxley, \emph{A note on polynomial congruences}, Recent progress in
  analytic number theory, Durham, 1979, vol.~1, Academic Press, London--New
  York, 1981, pp.~193--196.

\bibitem{Kloos}
H.~D. Kloosterman, \emph{On the representation of numbers in the form $ax^2
  +by^2 +cz^2 + dt^2$}, Acta Math. \textbf{49} (1926), 407--464.

\bibitem{Mun13}
R.~Munshi, \emph{Pairs of quadrics in 11 variables}, arXiv:1305.1461v1 (2013).

\bibitem{Reid}
M.~Reid, \emph{The complete intersection of two or more quadrics}, PhD Thesis,
  Trinity College, Cambridge 1972.

\end{thebibliography}
